%% file: main.tex
\documentclass[reqno,twoside]{amsart}
\usepackage{amsfonts}
\usepackage{amsmath,amssymb}
\usepackage{epic}
\usepackage{pstricks}
\usepackage{graphicx}
\usepackage{xcolor} % color, soul
\usepackage{upgreek}
\usepackage[ruled,vlined,linesnumbered]{algorithm2e}
\usepackage{enumitem}
\usepackage{booktabs}
\usepackage{subcaption}
\usepackage[colorlinks=true]{hyperref}
\usepackage{tikz}
\usetikzlibrary{arrows.meta, positioning}
\usepackage{a4wide,amsmath,amssymb,latexsym,amsthm}
\usepackage{graphicx}
\usepackage{verbatim}

\newtheorem{hypothesis}{} %Hypothesis}
\newtheorem{theorem}{Theorem}[]

\newtheorem{remark}[]{Remark}
\newtheorem{lemma}[]{Lemma}

\newtheorem{Problem}[]{Problem}

\newcommand{\be}{\begin{equation}}
\newcommand{\ee}{\end{equation}}
\newcommand{\ba}{\begin{eqnarray}}
\newcommand{\ea}{\end{eqnarray}}
\newcommand{\beq}{\begin{equation}}
\newcommand{\eeq}{\end{equation}}

\usepackage{color}

\numberwithin{equation}{section}

\usepackage{color}
\usepackage{stmaryrd}
\usepackage{comment}

\newcommand{\diag}{\operatorname{diag}}
\subjclass[2010]{}

\begin{document}
\title[Source reconstruction algorithms for coupled 
	parabolic systems]{Source reconstruction algorithms for coupled 
	parabolic systems from internal measurements of one scalar state}

\author{Cristhian  Montoya$^{*}$}
\address{C. Montoya,  Escuela de Ciencias Aplicadas e Ingenier\'{i}a, EAFIT, Medell\'{i}n Colombia}
\email{cdmontoyaz@eafit.edu.co}
\thanks{$^{*}$ Corresponding author - Cristhian Montoya:  cdmontoyaz@eafit.edu.co}

\author{Ignacio Brevis}
\address{I. Brevis,  School of Mathematical Sciences, University of Nottingham, UK}
\email{ignacio.brevis1@nottingham.ac.uk}

\author{David Bolivar}
\address{D. Bolivar,  Escuela de Ciencias Aplicadas e Ingenier\'{i}a, EAFIT, Medell\'{i}n Colombia}
\email{dboliva1@eafit.edu.co}
%%%%%%%%%%%%%%%%%%%%%%%%%%%%%%%%%%%%%%%%%%%%%%%%%%
% ABSTRACT %%%%%%%%%%%%%%%%%%%%%%%%%%%%%%%%%%%%%%%
\begin{abstract}
% This paper is devoted to the study of source reconstruction algorithms for coupled systems of heat equations, with either constant or spatially dependent coupling terms, where internal measurements are available from a reduced number of observed states. Two classes of systems are considered. The first comprises parabolic equations with constant zero-order coupling terms (through a matrix potential term). The second type considers parabolic equations coupled by a matrix potential that depends on spatial variables, which leads to the analysis of a non-self-adjoint operator. In all configurations, the source is assumed to be of separate variables, the temporal part is a known scalar function, and the spatial dependence is an unknown vector field. Several numerical examples using the finite element method in 1D and 2D are presented to show the reconstruction of space-dependent sources.
This paper is devoted to the study of source reconstruction algorithms for coupled systems of heat equations, with either constant or spatially dependent coupling terms, where internal measurements are available from a reduced number of observed states. Two classes of systems are considered. The first comprises parabolic equations with constant zero-order coupling terms (through a matrix potential term or via the diffusion matrix). The second type considers parabolic equations coupled by a matrix potential that depends on spatial variables, which leads to the analysis of a non-self-adjoint operator. In all configurations, the source is assumed to be of separate variables, the temporal part is a known scalar function, and the spatial dependence is an unknown vector field. Several numerical examples using the finite element method in 1D and 2D are presented to show the reconstruction of space-dependent sources.
		

\end{abstract}
% ABSTRACT %%%%%%%%%%%%%%%%%%%%%%%%%%%%%%%%%%%%%%%
%%%%%%%%%%%%%%%%%%%%%%%%%%%%%%%%%%%%%%%%%%%%%%%%%%
\maketitle
\tableofcontents

%%%%%%%%%%%%%%%%%%%%%%%%%%%%%%%%%%%%%%%%%%%%%%%%%%
% THEORY %%%%%%%%%%%%%%%%%%%%%%%%%%%%%%%%%%%%%%%%%
\input{01_theory}
% THEORY %%%%%%%%%%%%%%%%%%%%%%%%%%%%%%%%%%%%%%%%%
%%%%%%%%%%%%%%%%%%%%%%%%%%%%%%%%%%%%%%%%%%%%%%%%%%
%%%%%%%%%%%%%%%%%%%%%%%%%%%%%%%%%%%%%%%%%%%%%%%%%%
% NUMERICAL RESULTS %%%%%%%%%%%%%%%%%%%%%%%%%%%%%%
\input{02_numerics}
% NUMERICAL RESULTS %%%%%%%%%%%%%%%%%%%%%%%%%%%%%%
%%%%%%%%%%%%%%%%%%%%%%%%%%%%%%%%%%%%%%%%%%%%%%%%%%
%%%%%%%%%%%%%%%%%%%%%%%%%%%%%%%%%%%%%%%%%%%%%%%%%%
% CONCLUSIONS %%%%%%%%%%%%%%%%%%%%%%%%%%%%%%%%%%%%
\input{03_conclusions}
% CONCLUSIONS %%%%%%%%%%%%%%%%%%%%%%%%%%%%%%%%%%%%
%%%%%%%%%%%%%%%%%%%%%%%%%%%%%%%%%%%%%%%%%%%%%%%%%%
% APPENDIX %%%%%%%%%%%%%%%%%%%%%%%%%%%%%%%%%%%%%%%
%%%%%%%%%%%%%%%%%%%%%%%%%%%%%%%%%%%%%%%%%%%%%%%%%%
\input{04_appendix}
%%%%%%%%%%%%%%%%%%%%%%%%%%%%%%%%%%%%%%%%%%%%%%%%%%
% ACKNOWLEDGEMENTS %%%%%%%%%%%%%%%%%%%%%%%%%%%%%%%
\section*{Acknowledgements}\label{sec:acknow}
The first and third authors acknowledge the internal project 1109--11090202021 EAFIT University. 
Ignacio Brevis was supported by the Engineering and Physical Sciences Research Council (EPSRC), UK, under Grant EP/W010011/1. This work was also supported by the London Mathematical Society (LMS), under the Scheme 5 Research Grant 52402.
% ACKNOWLEDGEMENTS %%%%%%%%%%%%%%%%%%%%%%%%%%%%%%%
%%%%%%%%%%%%%%%%%%%%%%%%%%%%%%%%%%%%%%%%%%%%%%%%%%

%\IB{this reference O. A. Veliev. On nonselfadjoint Sturm-Liouville operators with matrix potentials. Mat. Zametki, 81(4):496–506, 2007. is in Russian, should we cite this instead \href{https://link.springer.com/article/10.1134/S0001434607030200}{Math Notes}. Also the DOI is from the English version}
%\bibliographystyle{plain}
\bibliographystyle{abbrv}
\bibliography{references}

\end{document}

%% file: 01_theory.tex
% MAIN PROBLEMS %%%%%%%%%%%%%%%%%%%%%%%%%%%%%%%%%% %%%%%%%%%%%%%%%%%%%%%%%%%%%%%%%%%%
\section{Main problems}\label{section.problems}
\subsection{Introduction}
Inverse problems that involve determining coefficients or sources in coupled systems of partial differential equations have attracted increasing interest in the last decade, especially in the case of parabolic or hyperbolic systems. Such problems naturally arise in diverse scientific and engineering fields, including fluid mechanics, biology, and medicine, among others. A challenging task involving inverse problems for coupled systems is whether it is possible to determine all sources (or coefficients) from a reduced number of measurements. Typically, the measurements are given for all state variables on local subsets (either from boundary local subsets or internal local subsets). Thus, a fundamental question is whether measuring fewer states than the number of sources or coefficients is sufficient for a unique determination. This question holds significant theoretical and practical interest.
	
From a theoretical point of view, this question leads to different guidelines depending on the type of coupling. For instance, linear coupling in low-order terms may exhibit matrix potentials with non-self-adjoint operators, where tools such as perturbation theory, semigroup theory, and spectral analysis are frequently employed to analyze these systems. Alternatively, linear coupling in the main operator can also lead to degenerate operators. On the other hand, the analysis of systems with nonlinear coupling terms is more complex than the previous one, which in turn requires knowledge of approximation theory as well as fixed-point arguments. The reconstruction of a general external source, either from internal or boundary measurements, is not determined uniquely (see, e.g.,~\cite{bookIsakov}), thus, by adding external forces into coupled systems, the inverse source problem becomes solvable if some a priori knowledge is assumed, i.e., if the unknown source is a characteristic function \cite{bookIsakov}, a point source \cite{2005Badia}, or a separated variables function \cite{1998Badia}. 
	
Regarding practical applications, there exist models in real life involving partial data of physical quantities, for example, pressure estimation from velocity phase-contrast MRI \cite{brown2014magnetic}, wireless communication where only some components of the electric fields are measured \cite{1943Dorn}, elastography \cite{2012doyleymodel}, molecular multi-photon transitions in laser fields \cite{bandrauk1993molecules}, heat transfer \cite{alifanov2012inverse}, and hybrid inverse problems \cite{2013Bal}.	
	
In this paper, we consider sources that separate into temporal and spatial components, where the temporal part is known and scalar; meanwhile, the spatial dependence is an unknown vector field. In fact, these external forces are associated with coupled systems of heat equations. Our primary goal is to recover the spatial distribution of these external forces in coupled parabolic systems from internal measurements from a reduced number of states. More precisely, our purpose is to provide source reconstruction algorithms, which show that it is possible to determine sources in coupled systems of heat equations from a limited number of components of the state in small subdomains. Our work is inspired by a previous methodology for determining sources in scalar equations (heat equation \cite{2013GOT} and wave equation \cite{1995Yamamoto}) and Stokes fluids \cite{2017GMO}. In that sense, we extend the existing results on scalar heat equations to coupled systems of heat equations. Furthermore, we present a Lipschitz-type stability result derived from Carleman inequalities. To preserve the flow of the main text, the full proof is presented in the appendix, where interested readers can find all necessary details.
		
To be more precise, we define our inverse source problems within an abstract framework.
Throughout the paper, $\Omega$ denotes a nonempty bounded domain of $\mathbb{R}^d$ ($d\in\mathbb{N}$) with smooth boundary $\partial\Omega$, and $\mathcal{O}\subset\Omega$ is a nonempty open subset, denoting the spatial subset of measurements. We denote by $n\in\mathbb{N}$ the number of equations and by $m\in\mathbb{N}$ the number of observed components, where $m<n$. Let  $A$ be an appropriate linear partial differential operator with domain $D(A)\subset L^2(\Omega,\mathbb{K})$, time-independent, and (possibly) with space-dependent coefficients, where $\mathbb{K}$ represents $\mathbb{R}$ or $\mathbb{C}$. The source term is of the form $\sigma(t)F(x)$, where $\sigma$ is a non-zero known scalar function, whereas $F=(f_1,\dots, f_n)^t$ is unknown ($(\cdot )^t$ denotes the transpose), and both are in suitable spaces to be defined later on. In the present paper, we will focus on the following problems:
	\begin{Problem}\label{Problem1}{\bf{Systems with constant coefficients}}
    
    The diagonal case with the same operator $A$ on each line, that is, \begin{equation}\label{sys.case1.matrizQ.ctes}
   		\left\{
    	\begin{array}{llll}
        \partial_t Y+I_nAY+QY={\mathbf{\sigma}}(t)F(x), &\mbox{in} &\Omega\times(0,T),\\
         Y(\, \cdot \, ,0)=0, &\mbox{in}&\Omega,        
    	\end{array}\right.
		\end{equation}
		where $Q\in \mathcal{M}_n(\mathbb{K})$ is a coupling matrix with constant coefficients (so-called constant potential term) and 
		$I_n$ is the identity matrix of size $n$. Here, we are interested in solving the following question: can we recover the source term $F=(f_1,\dots, f_n)^t$ in system \eqref{sys.case1.matrizQ.ctes} from a reduced number of measurements of the solution in $\mathcal{O}\times(0,T)$?
	\end{Problem}
	\begin{Problem}\label{Problem2}{\bf{Systems with space-dependent coefficients}}
    
		In this case, the coupling occurs on a zero-order term (so-called variable potential term in space) with space-dependent coefficients, i.e., 
		\begin{equation}\label{sys.Q.variable}
   		\left\{
    	\begin{array}{llll}
            \partial_t Y+AY+Q(x)Y=\sigma(t)F(x), &\mbox{in} &\Omega\times(0,T),\\
            Y(\,\cdot\, ,0)=0, &\mbox{in}&\Omega,    
    	\end{array}\right.
		\end{equation}
		where $Q\in \mathcal{M}_n(\mathbb{K})$ is a coupling matrix with space-dependent coefficients. Now, the question is: can we recover the source term $F=(f_1,\dots, f_n)^t$ in system \eqref{sys.Q.variable} from a reduced number of measurements of the solution in $\mathcal{O}\times(0,T)$? 	
	\end{Problem}

It is important to highlight that explicit formulas for source reconstruction and their algorithms for coupled parabolic systems have not been reported in the literature. Therefore, the main goal of this paper is to fill that gap with theoretical formulas for systems \eqref{sys.case1.matrizQ.ctes}--\eqref{sys.Q.variable}, where $A$ is the Laplace operator.    
	
\subsection{State of the art}
In the following review, we briefly outline the main classes of inverse problems related to coupled systems, particularly inverse coefficient problems and inverse source problems, highlighting works that specifically address coupled parabolic systems.
\begin{enumerate}
\item [a)]{\textit{Inverse coefficient problems in coupled systems.}} Under the structure of systems \eqref{sys.Q.variable} with $A=-\Delta$, the Bukhgeim-Klibanov method~\cite{BK81}, which is based on Carleman estimates, has been widely employed to obtain results for inverse coefficient problems \cite{2021KlibanovLi}. For instance, the articles \cite{2006CristofolGaitan,2009Yamamoto} proved Lipschitz-type stability inequalities for $2\times 2$ reaction-diffusion systems using a single observation in a subdomain. These results were further extended to a nonlinear case in~\cite{2012Cristofol}. In a different context, the recent work~\cite{2018NicoRobertoAxel} addressed the uniqueness and Lipschitz-type stability for cascade hyperbolic systems composed of $n$ equations, where internal measurements were available for all components of the solution except the last one. 
Another relevant work on inverse problems is \cite{yamamoto2009carleman}, where the author presents self-contained derivations of Carleman estimates for parabolic equations and demonstrates their powerful applications in proving uniqueness and stability for inverse problems and control theory.
Another recent work~\cite{2019Yamamoto} considered two coupled Schr\"{o}dinger equations (with $A=i\Delta$ in \eqref{sys.Q.variable}), in which the authors investigated the logarithmic stability for determining two potentials from internal observations of a single component of the solution. In contrast to the previous results, \cite{2019Yamamoto} combined Carleman estimates along with the Fourier-Bros-Iagolnitzer transform in order to prove the logarithmic stability. The paper \cite{2013CristofolGaitan} tackled the problem of identification for two discontinuous coefficients in a one-dimensional coupled parabolic system, observing only one component. The authors proved theoretical results based on Carleman-type inequalities, while numerical results were obtained using the finite difference method for both temporal and spatial discretization, jointly with the interior-point method (meaning that the numerical part is treated as an optimization problem).
To summarize, these articles show identifiability and stability results for hyperbolic and parabolic systems coupled through potential matrices (i.e., systems as \eqref{sys.Q.variable}). 
    
\item [b)]{\textit{Inverse source problems in coupled systems.}}  
    In relation to inverse source problems in coupled systems, the existing literature remains somewhat limited. The work \cite{2020lipschitz} addressed an identification problem for a cascade-type system of $n$ degenerate parabolic equations coupled through zero-order terms, with a general external source $G=G(x,t)$ instead of the separated form $\sigma(t)F(x)$. Under certain assumptions on $G$ and its temporal derivative, \cite{2020lipschitz} proved a Lipschitz-type stability result for determining the source $G$. The proof utilized Carleman inequalities, relying on the observation of a single component over a subdomain, along with simultaneous data of the $n$ components at a fixed positive time over the entire spatial domain. Similarly, the article \cite{2016Alabau} considered two wave equations coupled in cascade, providing identifiability and stability results for space-dependent sources from incomplete boundary observations. Nevertheless, the article did not provide insights into reconstruction algorithms. 
    
    The study most closely related to our research is presented in \cite{2017GMO}, which established theoretical and numerical source reconstruction algorithms for the Stokes system. The methodology in~\cite{2017GMO} combines spectral analysis, null controllability properties, and Volterra equations.     
    Similarly, our work integrates these techniques into a comprehensive framework (see Figure~\ref{fig:diagram} for an overview of the structure of this methodology).
 \end{enumerate}
Finally, we refer to the articles \cite{2019Wuetal, 2012Wu, 2009Jishan} and the references therein, which are closely linked to inverse problems for coupled parabolic systems.
	
As mentioned, our approach follows the guidelines proposed in \cite{1995Yamamoto, 2013GOT} for the wave equation, in~\cite{2013GOT} for the heat equation, and in~\cite{2017GMO} for the Stokes equations (see Figure~\ref{fig:diagram}). Accordingly, in Section~\ref{section.preli}, we present the necessary preliminary results concerning spectral analysis (where Riesz bases are now required instead of the standard Fourier bases) as well as controllability properties for coupled parabolic systems and essential aspects of Volterra equations.
	 
The remaining of the paper is organized as follows: Section~\ref{section.coeff.ctes} provides an explicit source reconstruction formula for coupled parabolic systems with constant coefficient matrices (as a potential term into the system), thus solving Problem~\ref{Problem1}. Section~\ref{section.coeff.variables} addresses Problem~\ref{Problem2}, presenting a reconstruction formula for a one-dimensional $2\times 2$ parabolic system coupled through a non-constant zero-order matrix potential. Here, the main challenge arises from the non-self-adjoint operator, necessitating careful use of eigenfunctions associated with non-self-adjoint Sturm–Liouville operators to achieve reconstruction based solely on interior measurements and a reduced number of states. Section~\ref{section.numerics} validates these theoretical findings through numerical experiments implemented with finite element methods and gradient descent algorithms. These experiments include reconstructions for Problem~\ref{Problem2} in 1D, 
%analyses of sensitivity to regularization parameters and observational noise, 
and reconstructions addressing Problem~\ref{Problem1} in 2D. Finally, the Lipschitz-type stability result is placed in the appendix.

\begin{figure}
    \centering
    \input{diagram}
    \caption{Diagram depicting the source reconstruction methodology for coupled parabolic systems in the form~\eqref{sys.case1.matrizQ.ctes}-\eqref{sys.Q.variable}.}
    \label{fig:diagram}
\end{figure}        
% PRELIMINARIES %%%%%%%%%%%%%%%%%%%%%%%%%%%%%%%%%%
 %%%%%%%%%%%%%%%%%%%%%%%%%%%%%%%%%%
%%%%%%%%%%%%%%%%%%%%%%%%%%%%%%%%%%%%%%%%%%%%%%%%%%
%%%%%%%%%%%%%%%%%%%%%%%%%%%%%%%%%%%%%%%%%%%%%%%%%%
\section{Preliminaries}\label{section.preli} 
	As mentioned, in order to provide source reconstruction formulas for coupled systems of heat equations, we shall mainly combine existing null controllability results for coupled parabolic systems with distributed control, 
	spectral properties for linear operators, and also integro-differential equations. 
	Thus, this section is devoted to presenting these topics. Henceforth, Problems \ref{Problem1} and \ref{Problem2} will be
	analyzed by considering the operator $A=-\Delta$ with domain $\mathcal{D}(A)=H^2(\Omega)\cap H_0^1(\Omega)$.

%--------------SPECTRAL PROPERTIES .--------------
\subsection{Spectral analysis}\label{subsection.spectral} 
We describe the eigenvalues and eigenfunctions of the non-self-adjoint operators $L,L^*:H^2(0,\pi)^2\cap H_0^1(0,\pi)^2\subset L^2(0,\pi)^2\rightarrow L^2(0,\pi)^2$ ($L^*$ denotes the adjoint operator) related to the following Sturm-Liouville problem
	\begin{equation}\label{eq.sturm.problem}
   	\left\{
    \begin{array}{lll}
        L\Phi:=-\Delta \Phi+V(x)\Phi=\lambda\Phi, &\mbox{in } (0,\pi),\\
		\Phi(0,t)=\Phi(\pi,t)= 0,  & \\
    \end{array}\right.
	\end{equation}
	where 
	\begin{equation}\label{Qmatriz2por2}
		V(x)=\begin{pmatrix}
			 0 & 0\\
			 q(x) & 0
		\end{pmatrix},
		\quad\mbox{and}\quad q\in L^\infty(0,\pi).
	\end{equation}
	Additionally, for every $k\in \mathbb{N}$, we consider the normalized eigenfunctions of
	the Laplace operator with Dirichlet boundary conditions over $(0,\pi)$, i.e., 
	\begin{equation}\label{eigen.laplacian.I}
		\varphi_k(x)=\sqrt{\displaystyle\frac{2}{\pi}}\sin(kx),\quad\mbox{and }\quad
		I_k(q):=\int\limits_0^\pi q(x)\varphi_k(x)\,dx,\quad \forall\,k\in \mathbb{N}.
	\end{equation}
	The next Lemma establishes biorthogonal Riesz bases associated to the operators $L$ and $L^*$. 
	Further details can be found in \cite[Proposition 2.1]{2017duprez} and \cite[Proposition 2.1]{2016Khdjaetal}. 
	\begin{lemma}\label{lemma.duprez}
		Consider the families {\scriptsize{
		\[\mathcal{B}=\Biggl\{\Phi_{1,k}=\begin{pmatrix} 0\\ \varphi_k\end{pmatrix}, 
		\Phi_{2,k}=\begin{pmatrix} \varphi_k\\ \psi_k\end{pmatrix}: k\in\mathbb{N}\Biggr\}\,\,
		\mbox{and}\,\,
		\mathcal{B}^*=\Biggl\{\Phi_{1,k}^*=\begin{pmatrix} \psi_k\\ \varphi_k\end{pmatrix},
		\Phi_{2,k}^*=\begin{pmatrix} \varphi_k\\ 0\end{pmatrix}, : k\in\mathbb{N}\Biggr\},\]}}
		where $\psi_k$ is defined for all $x\in (0,\pi)$ by
		\begin{equation}\label{eq.basis.heat}
   		\left\{
    	\begin{aligned}
        	\psi_k(x)&=\alpha_k\varphi_k(x)-\displaystyle\frac{1}{k}\int\limits_{0}^x\sin(k(x-\zeta))
        	(I_k(q)\varphi_k(\zeta)-q(\zeta)\varphi_k(\zeta))\,d\zeta,\\
        	\alpha_k &=\displaystyle\frac{1}{k}\int\limits_{0}^\pi\int\limits_0^x
        	\sin(k(x-\zeta))((I_k(q)\varphi_k(\zeta)-q(\zeta)\varphi_k(\zeta)))\varphi_k(x)\,d\zeta \,dx.\\
		\end{aligned}\right.
		\end{equation}
		Then, we have
		\begin{enumerate}
		\item [a)] The spectrum of $L^*$ and $L$ are given by $\rho(L^*)=\rho(L)=\{k^2: k\in \mathbb{N}\}.$
		\item [b)] For every $k\in \mathbb{N}$, the eigenvalue $k^2$ of $L^*$ has algebraic multiplicity $1$. Moreover, in this case,
			\begin{equation}\label{spectral.decomLstar.heat}
			\left\{
			\begin{array}{l}
			{\left(L^{*}-k^{2}Id\right)\Phi_{1, k}^{*}=I_{k}(q) \Phi_{2, k}^{*}}, \\ 
			{\left(L^{*}-k^{2}Id\right)\Phi_{2, k}^{*}=0}.
			\end{array}\right.
			\end{equation}
		\item [c)] For every $k\in \mathbb{N}$, the eigenvalue $k^2$ of $L$ has algebraic multiplicity $1$. Moreover, in this case,
			\begin{equation}\label{spectral.decomL.heat}
			\left\{
			\begin{array}{l}
			{\left(L-k^{2}Id\right)\Phi_{1, k}=0},\\ 
			{\left(L-k^{2}Id\right)\Phi_{2, k}=I_{k}(q) \Phi_{1, k}}.
			\end{array}\right.
			\end{equation}
		\item [d)] The sequences $\mathcal{B}$ and $\mathcal{B}^*$ are biorthogonal Riesz bases of $L^2(0,\pi)^2$.
		\item [e)] The sequence  $\mathcal{B}^*$ is a Schauder basis of $H_0^1(0,\pi)^2$ and $\mathcal{B}$ is its 
			biorthogonal basis in $H^{-1}(0,\pi)^2$.	
		\end{enumerate}
	\end{lemma}
	\begin{remark}\label{obs.spectrumgeneral}
		Spectral theory for regular and singular Sturm–Liouville problems encompasses a broad range of topics concerning eigenvalues and eigenfunctions (see e.g., \cite{bookLevitan,2016Shkalikov,bookKato}). In our case, 
        Lemma \ref{lemma.duprez} provides conditions on the potential $V(x)$ that ensure the eigenfunctions associated with the operators $L$ and $L^*$ form Riesz bases. In this work, we restrict our analysis to the one-dimensional case of a $2\times 2$ system featuring a non-self-adjoint matrix potential. This specific choice highlights the role of spectral properties within our reconstruction strategy and identifies critical points for future research. However, we emphasize that extending the spectral analysis to general $n\times n$ systems is significantly more challenging and involves a deeper and more delicate investigation beyond the scope of this paper. For further discussions on these advanced topics, we refer readers to~\cite{2007Veliev,2014Seref,2019Mityagin}.
	\end{remark}

%--------------CONTROLLABILITY.--------------
\subsection{Controllability}\label{subsection.control} 
	Depending on the structure of the coupling matrix $Q^t$  and the control matrix $B$, different results on 
	null controllability for the adjoint system of \eqref{sys.case1.matrizQ.ctes} can be derived, see for 
	instance \cite{2010Burgos,2014Benabdallah,2015EFCetal}. To our purpose, let us assume that 
	$Q^t\in L^\infty(\Omega)^{n^2}$ and $B\in\mathcal{M}_n(\mathbb{R})$ have the structure: 
	\begin{equation}\label{def.matrixQandB}
	Q^t=
	\begin{pmatrix}
		q_{11} & 0& 0&\cdots &0\\
		q_{21} & q_{22}& 0&\cdots  &0\\
		q_{31}&q_{32}& q_{33}&\cdots & 0\\
		\vdots &\vdots &\ddots & \ddots & \vdots\\
		q_{n1}&q_{n2}& \cdots &q_{n,n-1} & q_{nn}\\
	\end{pmatrix}\quad\mbox{and}\quad 
		B=\diag(0,0,\dots,0,1),	
	\end{equation}
	where 
	\begin{equation}\label{cond.matrixQ}
	q_{ij}\geq q_0>0\,\,\quad\mbox{in an open set}\,\,\mathcal{O}_0\subset\mathcal{O},\quad \forall\, i>j,\text{ and }  i,j=1,\dots,n.
	\end{equation}
	
	The following result holds from \cite{2010Burgos}.
	\begin{lemma}\label{lema.controlQ}
		Assume that $Q^t\in L^\infty(\Omega)^{n^2}$ and $B\in\mathcal{M}_n(\mathbb{R})$ are given by 
		\eqref{def.matrixQandB} and satisfy \eqref{cond.matrixQ}. Let $\tau\in (0,T]$ and $\Psi^{0}\in L^2(\Omega)^n$. Then, there exists a control function 
        $U^{(\tau)}=U^{(\tau)}(\Psi^{0})=(u_1^{(\tau)},\dots,u_n^{(\tau)})^t$ belonging to $L^2(0,T;L^2(\mathcal{O})^n)$
		such that the solution $\Psi$ of the problem 
		\begin{equation}\label{sys.heatctes2.adjoint}
   			\left\{
    		\begin{array}{llll}
        	-\partial_t \Psi-\Delta \Psi+Q^t\Psi=1_{\mathcal{O}}BU^{(\tau)}, &\mbox{in} &\Omega\times(0,\tau),\\
        	 \Psi=0 ,&\mbox{on}&\partial\Omega\times(0,\tau),\\        
			 \Psi(\,\cdot\, ,\tau)=\Psi^{0}, &\mbox{in}&\Omega,        
    		\end{array}\right.
		\end{equation}
		satisfies 
		\begin{equation}\label{identity.nullcontrol}
			\Psi(\cdot, 0)=0\quad\mbox{in}\,\,\Omega.	
		\end{equation}
		Moreover, there exists a positive constant $C_0$ depending only on $\Omega$ and $\mathcal{O}$ such that
		\begin{equation}\label{ine.control1}
		\|u_{n}^{(\tau)}\|_{L^2(0,T;L^2(\mathcal{O}))}\leq C_{0}e^{C(\tau)}\|\Psi^{0}\|_{L^2(\Omega)^n},	
		\end{equation}
		where
		\[C(\tau)=\tau+\frac{1}{\tau}+\max\limits_{j\leq i}\Bigl(|q_{ij}|^{2/(3(i-j)+3)}+\tau|q_{ij}|\Bigr).\]
	\end{lemma}
	\begin{remark}
		It is worth mentioning that Lemma \ref{lema.controlQ} was proved in \cite[Theorem 1.2]{2010Burgos}
		for a forward system instead of a backward system. However, the above configuration is more appropriate 
		for solving our inverse source problem (Problem \ref{Problem1}). In addition, notice that the null controllability property with one scalar control for general complete matrices is not possible
		\cite{2010Burgos}.
	\end{remark}  

%--------------VOLTERRA EQUATIONS.--------------
\subsection{Volterra equations}\label{subsection.volterra} 
    We recall some technical results related to scalar Volterra equations of the first and second kind, which will be useful in subsequent sections. For further details, we refer the reader to \cite{2013GOT} and \cite[Chapter 3]{GripenbergLondenStaffans}.
     
	\begin{lemma}\label{lemma_volterra1} 
		Assume $0<t<\tau\leq T$, $\sigma\in W^{1,\infty}(0,\tau)$,
        $\sigma(0)\neq 0$, $\eta\in L^{2}(0,\tau;L^{2}(\Omega))$. Then,     
    	there exists a unique $\theta \in H^1(0,\tau;L^2(\Omega))$ satisfying the Volterra equation 
    	\begin{equation}\label{volterrask}
    	\begin{array}{ll}
        	&\sigma(0)\partial_{t}\theta(x,t)+\displaystyle\int\limits_{t}^{\tau}( \sigma(s-t)\theta(x,s)
        	+\partial_t\sigma(s-t)\partial_t\theta(x,s))\,ds=\eta(x,t),\\
        	&\theta(x,\tau)=0.
        \end{array}
    \end{equation}
    Furthermore, there exists a constant $C>0$ depending on $\|\sigma\|_{W^{1,\infty}(0,\tau)}$ such that
    \begin{equation}\label{inevolterra}
      \|\theta\|_{H^1(0,\tau;L^{2}(\Omega))}\leq C\|\eta\|_{L^{2}(0,\tau;L^{2}(\Omega))}.
    \end{equation}
	\end{lemma}
	
	Now, let us  consider the operator $K: L^2(0,T;L^2(\Omega))\to H^1(0,T;L^2(\Omega))$ defined by
    \begin{equation}\label{def_Kg}
        (Kv)(x,t):=\int\limits_{0}^{t}\sigma(s)v(x,t-s)\,ds.
    \end{equation}

	The following Lemma provides properties related to the operators $K$ and $K^*$.  
	\begin{lemma}\label{lemma_volterra2} Assume $\sigma\in W^{1,\infty}(0,T)$ and $\sigma(0)\neq 0$. Then, there exist  positive constants 
	$C_1$ and $C_2$ depending only on $\Omega, T$ and
    $\|\sigma\|_{W^{1,\infty}(0,T)}$ such that
    \begin{equation}\label{inequality_kg}
        C_1\|Kv\|_{ H^{1}(0,T;L^2(\Omega))}\leq \|v\|_{L^{2}(0,T;L^2(\Omega))}\leq C_2\|Kv\|_{H^{1}(0,T;L^2(\Omega))}.
    \end{equation}
    Furthermore, the adjoint operator $K^{*}: H^{1}(0,T;L^2(\Omega))\to L^{2}(0,T;L^2(\Omega))$ is 
    given by
    \begin{equation}\label{adjK}
        (K^*\theta)(x,t)=\sigma(0)\partial_{t}\theta(x,t)+\int\limits_{t}^{T}(\sigma(s-t)\theta(x,s)    
        +\partial_t\sigma(s-t)\partial_{t}\theta(x,t))\,ds.
    \end{equation}
	\end{lemma}
    
% FIRST MAIN RESULT %%%%%%%%%%%%%%%%%%%%%%%%%%%%%%
%\input{Systems_with_constant_coefficients}
% FIRST MAIN RESULT %%%%%%%%%%%%%%%%%%%%%%%%%%%%%%
%%%%%%%%%%%%%%%%%%%%%%%%%%%%%%%%%%%%%%%%%%%%%%%%%%
%%%%%%%%%%%%%%%%%%%%%%%%%%%%%%%%%%%%%%%%%%%%%%%%%%
\section{First main result: coupled systems with constant coefficients}\label{section.coeff.ctes} 
	
    In this section, we present an explicit formula for reconstructing the spectral coefficients of the vector source $F(x)$ in terms of a family of null controls $U^{(\tau)}$, which are linked to the systems \eqref{sys.heatctes2.adjoint}, indexed by $\tau\in(0,T]$. 
    In other words, we solve Problem \ref{Problem1}.  To achieve this, we begin by considering the following system:
	\begin{equation}\label{sys.heatctes1}
   	\left\{
    	\begin{array}{llll}
        \partial_t Y-\Delta Y+QY=\sigma(t)F(x), &\mbox{in} &\Omega\times(0,T),\\
         Y=0, &\mbox{on}&\partial\Omega\times(0,T),\\        
		 Y(\,\cdot\, ,0)=0, &\mbox{in}&\Omega,        
    	\end{array}\right.
	\end{equation}
	where $Q\in L^\infty(\Omega)^{n^2}$ satisfies \eqref{def.matrixQandB}, \eqref{cond.matrixQ}.
    %, and $D$ satisfies ~\ref{itm:A1} and 
    %~\ref{itm:A2}. 
	\begin{remark}
	   Let us observe that, thanks to the assumptions on 
       %the diffusion matrix $D$ and 
       the potential matrix $Q$, for every 
	   source $\sigma(t)F(x)\in L^2(0,T;L^2(\Omega)^n)$, system \eqref{sys.heatctes1} %(resp. system \eqref{sys.heatctes2}) 
       admits a unique weak 
	   solution $Y\in C([0,T];L^2(\Omega)^n)\cap L^2(0,T;H^1_0(\Omega)^n)$. Additionally, by considering $\sigma\in W^{1,\infty}(0,T)$ and 
	   following, for example, the procedure given in  \cite[Chapter 3, Section 6]{bookladyvzenskaja1968}, solvability in the space 
	   $$W_2^{2,1}(\Omega\times(0,T)):=L^2(0,T;H^2(\Omega)^n\cap H^1_0(\Omega)^n)\cap H^1(0,T;L^2(\Omega)^n)$$ also holds for
	   system \eqref{sys.heatctes1}.
       %and \eqref{sys.heatctes2}. 
	\end{remark}
	\begin{remark}
		Before presenting the main results of this section, we present integral representation to the solution for system  \eqref{sys.heatctes1}, 
        %a common denominator to the systems \eqref{sys.heatctes1} and 
		%\eqref{sys.heatctes2} is the integral representation of their solutions, 
        which are possible thanks to the 
		linearity of this system. More precisely, from the Duhamel principle, the solution $Y$ from either of those systems can be written by 
		\begin{equation}\label{id.def.op.K}
			Y(x,t)=\int\limits_0^t\mathbf{\sigma}(s)W(x,t-s)\,ds,\quad (x,t)\in \Omega\times(0,T),
		\end{equation}
		where $W$ satisfies
		\begin{equation}\label{eq.linear.primal.forw}\small{
   			\left\{
   		 	\begin{array}{llll}
        		\partial_t W-\Delta W+QW=0, &\mbox{in} &\Omega\times(0,T),\\
        		 W=0, &\mbox{on}&\partial\Omega\times(0,T),\\        
				W(\,\cdot\, ,0)=\sigma(0)F(\cdot), &\mbox{in}& \Omega,        
    		\end{array}\right.}
    		%\quad\mbox{or}\quad
    		%\left\{
    		%\begin{array}{llll}
        	%	\partial_t W-D\Delta W=0,&\mbox{in} &\Omega\times(0,T),\\
        	%	W=0, &\mbox{on}&\partial\Omega\times(0,T),\\        
        	%	W(\,\cdot\, ,0)=\sigma(0)F(\cdot), &\mbox{in}& \Omega.    
            %    \end{array}\right.}
		\end{equation}
		Furthermore, since   $\partial_t Y(x,t)=\sigma(0)W(x,t)+\int\limits_{0}^{t} \partial_t
		\sigma(t-s)W(x,s)\,ds$,
		by evaluating at $t=T$ the main equations of \eqref{sys.heatctes1} 
        %and \eqref{sys.heatctes2}, 
        and using \eqref{id.def.op.K}, we obtain the following identity:
		\begin{equation}\label{eq.identityglobal}\small{
			\sigma(0)W(x,T)+\int\limits_0^T\partial_t\sigma(T-s)W(x,s)\,ds-\Delta Y(x,T)
			+\int\limits_0^T\mathbf{\sigma}(s)QW(x,T-s)\,ds=\sigma(T)F(x)}
		\end{equation}
		%and
		%\begin{equation}\label{eq.identityglobal2}
		%	\sigma(0)W(x,T)+\int\limits_0^T\partial_t\sigma(T-s)W(x,s)\,ds-D\Delta Y(x,T)=\sigma(T)F(x).
		%\end{equation}
	\end{remark}

	As previously mentioned, our inverse source problems closely depend on the null controllability properties of the adjoint system associated with \eqref{sys.heatctes1}, %and \eqref{sys.heatctes2}, 
    as well as on the spectral properties of the main operator. Consequently, the remainder of this section is dedicated to establishing explicit connections between the identity \eqref{eq.identityglobal} and 
    %\eqref{eq.identityglobal2} and 
    this essential theoretical component. Specifically, two key points must be addressed. First, the matrix $Q^t\in\mathcal{M}_n(\mathbb{R})$, corresponding to the coupling matrices $Q,\in\mathcal{M}_n(\mathbb{R})$ from system \eqref{sys.heatctes1},
    %and \eqref{sys.heatctes2}, 
    must satisfy the controllability conditions outlined in \eqref{cond.matrixQ}. 
    %and assumptions~\ref{itm:A1},~\ref{itm:A2}, and~\ref{itm:A3}.
    Second, since the coupling matrix $Q$  contains constant entries, we only need to consider $L^2$-eigenfunctions and eigenvalues of the Laplace operator with Dirichlet boundary conditions in $\Omega$ for the spectral assumptions. These eigenfunctions and eigenvalues will be denoted by $\{\varphi_k\}_{k\in\mathbb{N}}$ and $\{\lambda_k\}_{k\in\mathbb{N}}$, respectively. Finally, before stating our first result, we introduce some additional hypotheses.
	
	\begin{hypothesis}\label{Hy1}%\tag{H1}
		 Consider $\sigma\in W^{1,\infty}(0,T)$ with $\sigma(T)\neq 0$. Furthermore, for some $k\in\mathbb{N}$ 
		\begin{equation}\label{H1.eq1}
			a^{Q}_{{j,k}}(T):=\Biggl(1-\frac{\lambda_k}{\sigma(T)}\sum\limits_{i=1}^nm_{ij}(T)\Biggr)\neq 0,\quad \forall\, 
			j=1,\dots, n,
		\end{equation}		
		 where $M=(m_{ij}(t))=\displaystyle\int\limits_{0}^t\tilde\Phi_k(t)\tilde\Phi_k^{-1}(s)\sigma(s)\,ds$ and 
		$\tilde\Phi_k$ is a fundamental matrix associated with the linear ordinary differential system: $Z'+(\lambda_kI_n+Q)Z=0$.
	\end{hypothesis}
	
	\begin{hypothesis}\label{Hy2} 
		For every $k\in\mathbb{N}$ and $\Psi^0_{k}=(\varphi_k,\dots,\varphi_k)^t\in L^2(\Omega)^n$, consider Lemma \ref{lema.controlQ}, where $U=(u_1^{(\tau)},\dots,u_n^{(\tau)})^t$ is the 
		control function associated to the system \eqref{sys.heatctes2.adjoint} extended by zero in $(\tau,T)$, and  
		\begin{equation}\label{H2.eq1}
			BU=(0,0,\dots,u_n^{(\tau)})^t\quad\mbox{with}\,\,\, 
            (BU)_n=u_n^{(\tau)}\neq 0. 
		\end{equation}
	\end{hypothesis}
	
    Our first source reconstruction result is given in the following theorem.
	\begin{theorem}\label{th1} 
		Let hypotheses \ref{Hy1} and \ref{Hy2} be satisfied. Let $\theta_n$ and 
        $\hat{\theta}_n$ be solutions of 
        \eqref{volterrask} for $\eta=(BU)_n$ and $\hat{\eta}=(Q^tBU)_n$, respectively. 
        Then, for every solution $Y\in W_2^{2,1}(\Omega\times(0,T))$ of 
		\eqref{sys.heatctes1}, $F=(f_1,\dots, f_n)^t\in L^2(\Omega)^n$ satisfies the local reconstruction 
		identity
		\begin{equation}\label{formula1Qctes}
			\sum\limits_{j=1}^n a^{Q}_{j,k}(T)(f_j,\varphi_k)_{L^2(\Omega)}
			=\mathcal{C}_1(\Psi^0_{k})+\mathcal{C}_2(\Psi^0_{k})+\mathcal{C}_3(Q^t\Psi^0_{k}),
		\end{equation}
        where 
        \begin{equation*}
		\begin{aligned}
			\mathcal{C}_1(\Psi^0_{k})=&-\frac{\sigma(0)}{\sigma(T)}(y_n,\theta_n)
            _{H^1(0,T;L^2(\mathcal{O}))},\\
			\mathcal{C}_2(\Psi^0_{k})=&-\frac{1}{\sigma(T)}\int\limits_{0}^{T}\partial_t\sigma(T-s)
			(y_n,\theta_n)_{H^1(0,T;L^2(\mathcal{O}))}\,ds,\\
			\mathcal{C}_3(Q^t\Psi^0_{k})&=-\frac{1}{\sigma(T)}\int\limits_0^T\sigma(T-s)(y_n,
            \hat{\theta}_n)_{H^1(0,T;L^2(\mathcal{O}))}\,ds.
		\end{aligned}
		\end{equation*}
        
	\end{theorem}
%--------------PROOF OF FIRST PART OF THEOREM 1--------------
\begin{proof} Due to the fact that the proof essentially combines three different topics, we divide this proof into three steps: spectral representation, null controllability, and Volterra equations.	
	  
	\noindent{\textit{Step 1. Spectral representation.}} First, from the $L^2$-eigenfunctions 
	$\{\varphi_k\}_{k\in\mathbb{N}}$ of the Laplace operator with homogeneous Dirichlet conditions, note that the solution of \eqref{sys.heatctes1} can also be written as
	\begin{equation}\label{eq.recovery.heat0}
		Y(x,t)=\sum\limits_{k\in\mathbb{N}}Y_k(t)\varphi_k(x),
	\end{equation}
	where $Y_k(t)=(y^k_1(t),\dots, y^k_n(t))^t$ is the unique solution of the ordinary differential system 
	\begin{equation}\label{sys.odes.case.Q}
   	\left\{
    \begin{array}{llll}
        Y'_k(t)+(\lambda_kI_n+Q)Y_k(t)=\sigma(t)F_k,& &\\
		Y_k(0)=0, &&        
    \end{array}\right.
	\end{equation}
	and 
	$F_k=((f_1,\varphi_k)_{L^2(\Omega)},\dots, (f_n,\varphi_k)_{L^2(\Omega)})^t=:(f_1^k,\dots, f_n^k)^t$.
	
	By solving \eqref{sys.odes.case.Q}, for every $k\in\mathbb{N}$, we obtain
	\begin{equation}\label{eq.solutionY.caseQ}\small{
		Y_k(t)=\overbrace{\Biggl(\int\limits_{0}^t\tilde\Phi_k(t)\tilde\Phi_k^{-1}(s)\sigma(s)\,ds\Biggr)}^{M=(m_{ij}
		(t))_{i,j=1}^n}F_k
		=\Biggl(\sum\limits_{j=1}^n m_{1j}(t)f_j^k,\,\,\sum\limits_{j=1}^n m_{2j}(t)f_j^k,\,\, \dots,
		\sum\limits_{j=1}^n m_{nj}(t)f_j^k\Biggr)^t,}
	\end{equation}
	where $\tilde\Phi_k$ is a fundamental matrix associated with the ordinary differential system: $Y'_{k}+(\lambda_kI_n+Q)Y_{k}=0$ (see assumption \ref{Hy1}).
		
	Now, let us consider the sequence 
	$\mathcal{B}=\{\Psi^0_{k}\}_{k\in\mathbb{N}}$. Multiplying the identity \eqref{eq.identityglobal} by elements of 
	$\mathcal{B}$ and integrating in space,  we get 
	\begin{equation}\label{eq.recoveryheat.aux1}\small{
	\begin{aligned}
		\sigma(T)(F,\Psi^0_{k})_{L^2(\Omega)^n}
		=&\sigma(0)(W(T),\Psi^0_{k})_{L^2(\Omega)^n}
		+\int\limits_{0}^{T}\partial_t\sigma(T-s)(W(s),\Psi^0_{k})_{L^2(\Omega)^n}\,ds\\
		&-(\Delta Y(T),\Psi^0_{k})_{L^2(\Omega)^n}+\int\limits_0^T\sigma(T-s)(QW(s),\Psi^0_{k})_{L^2(\Omega)^n}\,ds.
	\end{aligned}}
	\end{equation}
	Using \eqref{eq.recovery.heat0} and the fact that $\{\varphi_k\}_{k\in\mathbb{N}}$ are eigenfunctions for the Laplace operator, the third term at the right-hand side of \eqref{eq.recoveryheat.aux1} can be transformed as follows:
	\begin{equation}\label{eq.recoveryheat.aux2}\small{
		\begin{aligned}
			-(\Delta Y(T),\Psi^0_{k})_{L^2(\Omega)^n}
            %=-(Y(T),\Delta\Psi^0_{k})_{L^2(\Omega)^n}
			=\lambda_k (Y(T),\Psi^0_{k})_{L^2(\Omega)^n}
            =\lambda_k\sum\limits_{i=1}^n (y_i(T),\varphi_k)_{L^2(\Omega)}
			=\lambda_k\sum\limits_{i=1}^n y^k_i(T).
		\end{aligned}}
	\end{equation}
	From \eqref{eq.recoveryheat.aux1} and \eqref{eq.recoveryheat.aux2}, at this moment our reconstruction formula is given by
	\begin{equation*}\small{
	\begin{aligned}
		(F,\Psi^0_{k})_{L^2(\Omega)^n}
		=&\frac{\sigma(0)}{\sigma(T)}(W(T),\Psi^0_{k})_{L^2(\Omega)^n}
		+\frac{1}{\sigma(T)}\int\limits_{0}^{T}\partial_t\sigma(T-s)(W(s),\Psi^0_{k})_{L^2(\Omega)^n}\,ds\\
		&+\frac{\lambda_k }{\sigma(T)}\sum\limits_{i=1}^n y^k_i(T)
		+\frac{1}{\sigma(T)}\int\limits_0^T\sigma(T-s)(QW(s),\Psi^0_{k})_{L^2(\Omega)^n}\,ds,
        \quad \sigma(T)\neq 0,
	\end{aligned}}
	\end{equation*}
	which is equivalent to (using \eqref{eq.solutionY.caseQ}): 
	\begin{equation}\label{eq.recoveryheat.aux3}\small{
	\begin{aligned}
		&\sum\limits_{j=1}^n\Biggl(1-\frac{\lambda_k}{\sigma(T)}\sum\limits_{i=1}^nm_{ij}(T)\Biggr)f_j^k
		=\frac{\sigma(0)}{\sigma(T)}(W(T),\Psi^0_{k})_{L^2(\Omega)^n}\\
		&+\frac{1}{\sigma(T)}\int\limits_{0}^{T}\partial_t\sigma(T-s)(W(s),\Psi^0_{k})_{L^2(\Omega)^n}\,ds
		+\frac{1}{\sigma(T)}\int\limits_0^T\sigma(T-s)(QW(s),\Psi^0_{k})_{L^2(\Omega)^n}\,ds,
	\end{aligned}}
	\end{equation}
	where $\quad \Biggl(1-\frac{\lambda_k}{\sigma(T)}\sum\limits_{i=1}^nm_{ij}(T)\Biggr)\neq 0$ for all $i,j=1,\dots n$. 

    At this point, \eqref{eq.recoveryheat.aux3} shows a reconstruction formula for the coefficients of the source $F$, but using global terms in $L^2(\Omega)^n$.
    
	\noindent{\textit{Step 2. Controllability.}} Now, in order to replace the global terms 
	$(W(s),\Psi^0_{k})_{L^2(\Omega)^n}$ and $(QW(s),\Psi^0_{k})_{L^2(\Omega)^n}$ by local terms in 
	$L^2(0,s;L^2(\mathcal{O})^n)$, for every $s\in(0,T]$, we use the null controllability property for adjoint systems
	associated to \eqref{sys.heatctes1}. In other words, we apply the hypothesis \ref{Hy2}.
	Thus, if $\Psi$ denotes the adjoint state of $Y$ and $\overline\Psi:=Q^t\Psi$, Lemma \ref{lema.controlQ} 
	guarantees the existence of a control function $U:=U^{(s)}(\Psi^0_{k})\in L^2(0,s;L^2(\mathcal{O})^n)$ such that the systems
	\begin{equation}\label{sys.control.coupled1.ctes}
   	\left\{
    	\begin{array}{llll}
        -\partial_t \Psi-\Delta\Psi+Q^t\Psi=1_{\mathcal{O}}BU, &\mbox{in} &\Omega\times(0,s),\\
        \Psi=0, &\mbox{on}&\partial\Omega\times(0,s),\\        
		 \Psi(\,\cdot\, ,s)=\Psi^0_{k}(\cdot), &\mbox{in}&\Omega        
    \end{array}\right.
    \end{equation}
    and
	\begin{equation}\label{sys.control.coupled22.ctes}
    \left\{
    	\begin{array}{llll}
        -\partial_t \overline\Psi-\Delta\overline\Psi+Q^t\overline\Psi=1_{\mathcal{O}}Q^{t}BU, &\mbox{in} 
        &\Omega\times(0,s),\\
     	\overline \Psi=0, &\mbox{on}&\partial\Omega\times(0,s),\\        
		\overline\Psi(\,\cdot\, ,s)=Q^t\Psi^0_{k}(\cdot), &\mbox{in}&\Omega,        
    \end{array}\right.
	\end{equation}
	satisfy 
	\begin{equation}\label{cond.controlcero}
		\Psi(x,0)=\overline\Psi(x,0)=0,\quad \forall x\in\Omega.
	\end{equation}
	On the other hand, multiplying the first system of \eqref{eq.linear.primal.forw}  by $\Psi$, with $\Psi$ solution of
	\eqref{sys.control.coupled1.ctes} and \eqref{cond.controlcero}, and integrating by parts 
	in $L^2(0,s;L^2(\Omega)^n)$, we obtain (after extending $U$ by zero at $(s,T)$)
	\begin{equation}\label{id.aux.globaltolocal1}
		(W(s),\Psi^0_{k})_{L^2(\Omega)^n}=-(W,1_{\mathcal{O}}BU)_{L^2(0,s;L^2(\Omega)^n)}
		=-(W,BU)_{L^2(0,T;L^2(\mathcal{O})^n)}.
	\end{equation}
	Analogously, by multiplying the first system of \eqref{eq.linear.primal.forw}  by
    $\overline\Psi$, with $\overline\Psi$ solution of \eqref{sys.control.coupled22.ctes} and \eqref{cond.controlcero}, we have
	\begin{equation}\label{id.aux.globaltolocal2}\small{ 
		(QW(s),\Psi^0_{k})_{L^2(\Omega)^n}=-(W,1_{\mathcal{O}}Q^{t}BU)_{L^2(0,s;L^2(\Omega)^n)}
		=-(W,Q^{t}BU)_{L^2(0,T;L^2(\mathcal{O})^n)}.}
	\end{equation}
    Note that, thanks to the identities \eqref{id.aux.globaltolocal1} and \eqref{id.aux.globaltolocal2}, the global terms (in norm $L^2(\Omega)^n$) in the reconstruction formula \eqref{eq.recoveryheat.aux3} are now local terms in the subdomain $\mathcal{O}\times(0,T)$. 

    The rest of the proof establishes the relationship among the control functions $BU$ and $Q^{t}BU$ with systems of Volterra equations. 
    
	\noindent{\textit{Step 3. Volterra equations.}} In order to ensure an appropriate connection with \eqref{id.aux.globaltolocal1} and
	\eqref{id.aux.globaltolocal2}, in this step we  adapt 
	subsection \ref{subsection.volterra} to the case of systems of Volterra equations.  
    For every $k\in\mathbb{N}$, let $\Theta=(\theta_1,\theta_2,\dots,\theta_n)^t$ be a solution to $n$ copies of \eqref{volterrask} with right-hand side given by \eqref{H2.eq1}, i.e., for every $k\in\mathbb{N}$, Lemma \ref{lemma_volterra1} and Lemma \ref{lemma_volterra2} guarantees the existence of a unique function 
    $\Theta=(\theta_1,\theta_2,\dots,\theta_n)^t$ satisfying the following system (see \eqref{volterrask} and \eqref{adjK})
    \begin{equation}\label{exp.ññ1}
        \overline{K^*}(\Theta)(x,t)= 1_{\mathcal{O}}BU(x,t)\quad\mbox{ and}\quad 
	    \Theta(x,s)=0,\quad \forall k\in\mathbb{N},
    \end{equation} 
    where $\overline{K^*}(\Theta):=(K^*\theta_1, K^*\theta_2,\dots, K^*\theta_n)^t$, and  $K^*\theta_{j}$ is the adjoint operator given in Lemma \ref{lemma_volterra2}, for $j=1,\dots, n$. 
    
    Now, replacing the right-hand sides of the system \eqref{exp.ññ1} into \eqref{id.aux.globaltolocal1}, we obtain
	\begin{equation}\label{eq.globallocal1}
    \small{
    \begin{array}{ll}
         (W(s),\Psi^0_{k})_{L^2(\Omega)^n}
         &=-(W,BU)_{L^2(0,T;L^2(\mathcal{O})^n)}\vspace{0.1cm}\\
         &=-(W,\overline{K^*}(\Theta))_{L^2(0,T;L^2(\mathcal{O})^n)}\vspace{0.1cm}\\
         & =-(w_n,K^*\theta_n)_{L^2(0,T;L^2(\mathcal{O}))},
    \end{array}}   
	\end{equation}
    where the last identity comes from  \eqref{exp.ññ1} and \eqref{H2.eq1} (consequence of having null controllability with only one scalar component in the control function $BU$).

    From \eqref{def_Kg}, it is easy to deduce that $Y=KW$, in particular $y_n=Kw_n$. Therefore, using \eqref{eq.globallocal1} follows
	\begin{equation}\label{eqs.final.proof1}
		(w_n,K^*\theta_n)_{L^2(0,T;L^2(\mathcal{O}))}
		=(y_n,\theta_n)_{H^1(0,T;L^2(\mathcal{O}))}.	
	\end{equation}

    Now, applying the same arguments above, it is possible to find a function $\hat\Theta=(\hat\theta_1,\hat\theta_2,\dots,\hat\theta_n)^t$ solution to 
    \begin{equation}\label{exp.ññ2}
        \overline{K^*}(\hat\Theta)(x,t)= 1_{\mathcal{O}}Q^{t}BU(x,t)\quad\mbox{ and}\quad 
	   \hat\Theta(x,s)=0,\quad \forall k\in\mathbb{N}. 
    \end{equation}
    Replacing the right-hand sides of \eqref{exp.ññ2} into \eqref{id.aux.globaltolocal2}, we have
	\begin{equation}\label{eq.globallocal2}
    \small{
    \begin{array}{ll}
         (QW(s),\Psi^0_{k})_{L^2(\Omega)^n}&=-(W,Q^{t}BU)_{L^2(0,T;L^2(\mathcal{O})^n)}\vspace{0.1cm}\\
         &=-(W,\overline{K^*}(\hat\Theta))_{L^2(0,T;L^2(\mathcal{O})^n)}\vspace{0.1cm}\\
         & =-(w_n,K^*\hat\theta_n)_{L^2(0,T;L^2(\mathcal{O}))}.
    \end{array}}   
	\end{equation}
     Again, from \eqref{def_Kg}, it is easy to deduce that $Y=KW$, in particular $y_n=Kw_n$. Therefore, using \eqref{eq.globallocal2} it follows
	\begin{equation}\label{eqs.final.proof2}
		(w_n,K^*\hat\theta_n)_{L^2(0,T;L^2(\mathcal{O}))}
		=(y_n,\hat\theta_n)_{H^1(0,T;L^2(\mathcal{O}))}.
	\end{equation}  
	Finally, putting together \eqref{eq.globallocal1},\eqref{eqs.final.proof1},\eqref{eq.globallocal2}, \eqref{eqs.final.proof2}, and replacing in 
    \eqref{eq.recoveryheat.aux3}, we obtain the desired reconstruction formula \eqref{formula1Qctes}. This completes the proof of Theorem \ref{th1}.
	\end{proof}  
    
	\begin{remark}
        Under hypothesis~\ref{Hy1}, the reconstruction formula~\eqref{formula1Qctes} is valid for all $k\in\mathbb{N}$ in several specific cases of the time-dependent function $\sigma(t)$ (see~\cite{2013GOT,2017GMO}). This includes the following scenarios:
		\begin{enumerate}
			\item [a)] $\sigma(t)=\sigma_0$ constant.
			\item [b)] $\sigma$ a non-negative and increasing function.
			\item [b)] $\sigma(t)=1+\frac{1}{2}\cos\Bigl(\frac{4t}{T-\varepsilon}\Bigr)$  for $t<T-\varepsilon$, and $\sigma(t)=\frac{3}{2}$ for $t>T-\varepsilon$.
		\end{enumerate}	
	\end{remark}

% SECOND MAIN RESULT %%%%%%%%%%%%%%%%%%%%%%%%%%%%%
%\input{Systems_with_space_dependent_coefficients}
% SECOND MAIN RESULT %%%%%%%%%%%%%%%%%%%%%%%%%%%%%
%%%%%%%%%%%%%%%%%%%%%%%%%%%%%%%%%%%%%%%%%%%%%%%%%%
\section{Second main result: coupled systems with space dependent coefficients}\label{section.coeff.variables} 
    In this section, we focus on the problem of recovering the spatial component of a source term in the coupled system of second-order parabolic equations corresponding to Problem \ref{Problem2}. Specifically, we consider the 1D system:
	\begin{equation}\label{sys.2por2heatprimal}
   	\left\{
    \begin{array}{lll}
        \partial_t Y+(\overbrace{-\Delta+Q(x)}^{L})Y=\sigma(t)F(x), &\mbox{in } (0,\pi)\times(0,T),\\
		Y(0,t)=Y(\pi,t)= 0,  &\mbox{in } (0,T),\\
        Y(\,\cdot\, ,0)=0, &\mbox{in } (0,\pi),        
    \end{array}\right.
	\end{equation}	
	where $L:H^2(0,\pi)^2\cap H_0^1(0,\pi)^2\subset L^2(0,\pi)^2\rightarrow L^2(0,\pi)^2$ and $Q$ is given  by
	\begin{equation*}
		Q(x)=\begin{pmatrix}
			 0 & 0\\
			 q(x) & 0
		\end{pmatrix}
		\quad\mbox{and}\,\, q(x)\geq C>0,\quad  q\in L^\infty(0,\pi)\cap W^{1,\infty}(\tilde{\mathcal{O}}),\,\,
		\tilde{\mathcal{O}}\subset\mathcal{O}\subset (0,\pi).
	\end{equation*}
	
Unlike the previous section, where the coupling terms were constant and the self-adjoint structure of the Laplace operator $A=-\Delta$ enabled the construction of auxiliary control systems to estimate global terms (see Step 2 in the proof of Theorem~\ref{th1}), the strategy for recovering the source $F(x)$ in system~\eqref{sys.2por2heatprimal} requires a different approach. The key difficulty lies in the spectral analysis: since the operators
\[L:=-\Delta+Q(x)\quad \mbox{and}\quad L^*:=-\Delta+Q^{t}(x)\]
are not self-adjoint, and therefore one cannot rely on classical Fourier-type representations for spectral decomposition. Instead, the analysis must be carried out in terms of Riesz bases. This shift is fundamental and represents the main novelty of this section compared to the previous one.
         
As discussed in Remark \ref{obs.spectrumgeneral}, the spectral study of systems of $n\times n$ equations (with matrix operator $L$) is complex and beyond the scope of this work. Furthermore, for higher dimensions in space, an exhaustive revision from spectral theory might be useful in order to establish a source reconstruction formula as in the theorem below (Theorem \ref{th3}) for more general systems than the one considered here in \eqref{sys.2por2heatprimal} (see~\cite{bookKato}).   
	
	Before presenting the main theorem of this section, we will first consider some hypotheses.  
	
	\begin{hypothesis}\label{Hy7}
		Consider $\sigma\in W^{1,\infty}(0,T)$ with $\sigma(T)\neq 0$. Furthermore, for some $k\in\mathbb{N}$:
		\begin{equation}\label{H7.eq}
		\begin{array}{lll}
			a_k^L(T)&:=\sigma(T)\Biggl(1-\displaystyle\frac{k^2}{\sigma(T)}\displaystyle\int\limits_{0}^Te^{-k^2(T-s)}\sigma(s)\,ds\Biggr)\neq 0,\\
			b_k^L(T)&:=-I_k(q)\Biggl(1-k^2\displaystyle\int\limits_{0}^T(T-s)e^{-k^2(T-s)}\sigma(s)\,ds\Biggr).
		\end{array}	
		\end{equation}
	\end{hypothesis}
 	\begin{hypothesis}\label{Hy8} 
 		Consider Lemma \ref{lemma.duprez} (spectral analysis for the operators $L$ and $L^*$). 
	\end{hypothesis}
	
 	\begin{hypothesis}\label{Hy9}
 		For every $k\in\mathbb{N}$, any $\tau\in (0,T]$, and using Lemma \ref{lemma.duprez}, consider $\Phi_{1,k}^*+\Phi_{2,k}^*$ as the initial datum for the following distributed control system
	\begin{equation}\label{eq.dd}
   		\left\{
    	\begin{array}{llll}
    	-\partial_t \Psi+L^*\Psi=(0, 1_{\mathcal{O}}u_2),  &\mbox{in } (0,\pi)\times(0,\tau),\\
		\Psi(0,t)=\Psi(\pi,t)=0,  &\mbox{in } (0,\tau),\\
        \Psi(\,\cdot\, ,\tau)=\Psi^{0}(\cdot), &\mbox{in } (0,\pi),        
   		\end{array}\right.
	\end{equation}
        where $\Psi$ satisfies $\Psi(\,\cdot\, ,0)=0\quad\mbox{in}\,\,(0,\pi)$ (i.e., system \eqref{eq.dd} satisfies the null controllability property, 
        see \cite[Theorem 1.1]{2017duprez}). 
 	\end{hypothesis}

	\begin{theorem}\label{th3}
		Let hypotheses \ref{Hy7}, \ref{Hy8}, and \ref{Hy9} be satisfied. 
        Let $\theta$ be a solution to \eqref{volterrask} with right-hand side 
        given by $1_{\mathcal{O}}u_2$. Then, for any solution $Y\in W_2^{2,1}((0,\pi)\times(0,T))$ of \eqref{sys.2por2heatprimal}, 
		the source $F=(f_1,f_2)^{t}\in L^2(0,\pi)^2$ satisfies 
		\begin{equation}\label{formula.recoveryth3}
			a^{L}_{k}(T)\Bigl(f_1^{\varphi_k}+f_1^{\psi_k}+f_2^{\varphi_k}\Bigr)
			+b^{L}_{k}(T)f_1^{\varphi_k}
			=\mathcal{C}_1(\mathcal{B},\mathcal{B}^*)+
            \mathcal{C}_2(\mathcal{B},\mathcal{B}^*),
		\end{equation} 
        where 
        \[\mathcal{C}_1(\mathcal{B},\mathcal{B}^*)=-\sigma(0)(y_2,\theta)_{H^1(0,T;L^2(\mathcal{O}))},\]
        \[\mathcal{C}_2(\mathcal{B},\mathcal{B}^*)=-\displaystyle\int\limits_{0}^{T}\partial_t\sigma(T-s)(y_2,\theta)_{H^1(0,T;L^2(\mathcal{O}))}\,ds,\]
		and $f_1^{\varphi_k}:=(f_1,\varphi_k)_{L^2(0,\pi)},\,\,f_1^{\psi_k}:=(f_1,\psi_k)_{L^2(0,\pi)}$ and 
		$f_2^{\varphi_k}:=(f_2,\varphi_k)_{L^2(0,\pi)}$ are Riesz coefficients, with 
        $\varphi_k$ and $\psi_k$ defined in \eqref{eigen.laplacian.I} and \eqref{eq.basis.heat}, respectively. 	
	\end{theorem}
	\begin{proof} The proof is divided into four steps: spectral property, integral representation, controllability, and Volterra equations.
	
	\noindent{\textit{Step 1. Spectral property.}} The main ingredient in this step is the spectrum of the operators $L$ and $L^*$. First, considering \eqref{eigen.laplacian.I}, for every $k\in \mathbb{N}$, 
	we define the normalized eigenfunctions $\varphi_k(x)=\sqrt{\frac{2}{\pi}}\sin(kx)$ of the Laplace operator 
	$\partial_{xx}$ with Dirichlet boundary conditions over $(0,\pi)$.  
	
	Using Lemma \ref{lemma.duprez}, the solution to \eqref{sys.2por2heatprimal} can be represented in 
	the form:
	\begin{equation}\label{amamaama}
		Y(x,t)=\sum\limits_{k\in\mathbb{N}}\alpha_k(t)\Phi_{1,k}(x)+\beta_k(t)\Phi_{2,k}(x),
	\end{equation}
	where $\alpha_k(t)=(Y(t),\Phi_{1,k})_{L^2(0,\pi)^2}, \beta_k(t)=(Y(t),\Phi_{2,k})_{L^2(0,\pi)^2}$, and 
	$\Phi_{1,k},\Phi_{2,k}$ belong to the family $\mathcal{B}$, which constitutes a Riesz basis of $L^2(0,\pi)^2$. 
	Moreover, the sequences $\{\alpha_{k}(T)\}_{k\in\mathbb{N}}$ and $\{\beta_k(T)\}_{k\in\mathbb{N}}$ can be determined 
	in explicit form. In fact, after replacing \eqref{amamaama} into \eqref{sys.2por2heatprimal} and multiplying twice, one by $\Phi_{1,k}^*$ and another one by $\Phi_{2,k}^*$ (see Lemma \ref{lemma.duprez}), we can deduce the coupled system of ordinary differential equations
    (a consequence of the biorthogonal property between the families $\mathcal{B}$ and $\mathcal{B}^*$ as well as the relations \eqref{spectral.decomLstar.heat} and \eqref{spectral.decomL.heat}):
	\[\frac{d}{dt}\begin{pmatrix} \alpha_{k}(t)\\ \beta_{k}(t)\end{pmatrix}
		+\begin{pmatrix} k^2&I_k(q)\\0 & k^2\end{pmatrix}\begin{pmatrix} \alpha_{k}(t)\\ \beta_{k}(t)\end{pmatrix}
		=\sigma(t)\begin{pmatrix} (F,\Phi_{1,k}^*)_{L^2(0,\pi)^2}\\ (F,\Phi_{2,k}^*)_{L^2(0,\pi)^2}
	\end{pmatrix},\quad \forall\,k\in \mathbb{N},\]
	where $I_k(q):=\displaystyle\int\limits_0^\pi q(x)\varphi_k(x)\,dx$.
	
	From the definitions of $\Phi_{1,k}^*$ and $\Phi_{2,k}^*$ (see Lemma \ref{lemma.duprez}), the explicit solution
	to the previous systems corresponds to: for every $k\in \mathbb{N}$,
	\begin{equation}\label{eq.sol.coeffs.heat}
	\left\{
	\begin{aligned}
		\beta_k(t)&=f_1^{\varphi_k}\int\limits_{0}^te^{-k^2(t-s)}\sigma(s)\,ds.\\
		\alpha_k(t)&=(f_1^{\psi_k}+f_2^{\varphi_k})\int\limits_{0}^te^{-k^2(t-s)}\sigma(s)\,ds
		-I_k(q)\int\limits_{0}^te^{-k^2(t-s)}\beta_k(s)\,ds,\\
	\end{aligned}\right.
	\end{equation}
	where by simplicity we define $f_1^{\varphi_k}:=(f_1,\varphi_k)_{L^2(0,\pi)},\,\,f_1^{\psi_k}:=(f_1,\psi_k)_{L^2(0,\pi)}$
	and $f_2^{\varphi_k}:=(f_2,\varphi_k)_{L^2(0,\pi)}$.
	
	\noindent{\textit{Step 2. Integral representation.}} 
    Using the Duhamel principle, the solution to the system 
    \eqref{sys.2por2heatprimal} can be written by
	\begin{equation}\label{sys.2por2heatw}\small{
		Y(x,t)=\int\limits_0^t\mathbf{\sigma}(s)W(x,t-s)\,ds\quad\mbox{and}\quad
		\left\{
    	\begin{array}{lll}
        	\partial_t W+LW=0, &\mbox{in } (0,\pi)\times(0,T),\\
			W(0,t)=W(\pi,t)= 0,  &\mbox{in } (0,T),\\
        	W(\,\cdot\, ,0)=\sigma(0)F(\cdot), &\mbox{in } (0,\pi).        
    	\end{array}\right.}
	\end{equation} 
	
	Since $\partial_tY(x,t)=\mathbf{\sigma}(0)W(x,t)+\int\limits_{0}^{t}\partial_t\sigma(t-s)W(x,s)\,ds$, then, by 
	evaluating the main equations of \eqref{sys.2por2heatprimal} in $t=T$, and multiplying by elements of the family 
	$\mathcal{B^*}$ and integrating over $(0,\pi)$, we obtain
	\begin{equation}\label{eq.a}
	\begin{aligned}
		(\sigma(T)F,\Phi_{1,k}^*+\Phi_{2,k}^*)_{L^2(0,\pi)^2}
		=&\sigma(0)(W(T),\Phi_{1,k}^*+\Phi_{2,k}^*)_{L^2(0,\pi)^2}\\
		&+\int\limits_{0}^{T}\partial_t\sigma(T-s)(W(s), \Phi_{1,k}^*+\Phi_{2,k}^*)_{L^2(0,\pi)^2}\,ds\\
		&+(LY(T),\Phi_{1,k}^*+\Phi_{2,k}^*)_{L^2(0,\pi)^2}.
	\end{aligned}
	\end{equation}
	Using Lemma \ref{lemma.duprez} and \eqref{amamaama}, the last term in the right-hand side of \eqref{eq.a} can 
	be written as follows:
	\begin{equation}\label{eq.b}
		\begin{aligned}
			(LY(T),\Phi_{1,k}^*+\Phi_{2,k}^*)_{L^2(0,\pi)^2}=&(Y(T),I_k(q)\Phi_{2,k}^*)_{L^2(0,\pi)^2}
			+k^2(\alpha_k(T)+\beta_k(T))\\
			=& (I_k(q)+k^2)\beta_k(T)+k^2\alpha_{k}(T).
		\end{aligned}
	\end{equation}
	Thus, at this stage, our reconstruction formula (with global terms in 
    $L^2(0,\pi)^2$) is given by
	\begin{equation}\label{eq.c}
	\begin{aligned}
		\sigma(T)(F,\Phi_{1,k}^*+\Phi_{2,k}^*)_{L^2(0,\pi)^2}=
		&\sigma(0)(W(T),\Phi_{1,k}^*+\Phi_{2,k}^*)_{L^2(0,\pi)^2}\\
		&+\int\limits_{0}^{T}\partial_t\sigma(T-s)(W(s), \Phi_{1,k}^*+\Phi_{2,k}^*)_{L^2(0,\pi)^2}\,ds\\
		&+(I_k(q)+k^2)\beta_k(T)+k^2\alpha_{k}(T).
	\end{aligned}
	\end{equation}
	where $\alpha_k(T)$ and $\beta_k(T)$ have been obtained in \eqref{eq.sol.coeffs.heat}.
	
	\noindent{\textit{Step 3. Controllability.}} As in the proof of Theorem~\ref{th1}, where global terms were replaced using two adjoint systems and their corresponding controllability results, we follow a similar strategy here. Specifically, the global terms $(\mathbf{w}(\,\cdot\, ,\tau),\Phi_{1,k}^*+\Phi_{2,k}^*)_{L^2(0,\pi)^2}$ for $\tau\in (0,T]$, are replaced by local terms in 
    $L^2(0,\tau;L^2(\mathcal{O}))^2$, using the controllability assumption stated in Hypothesis~\ref{Hy9}. 
	
	Thus, reasoning as in the previous section, for every $\tau\in(0,T]$ follows
	\begin{equation}\label{eq.ff}\small{
		(W(\,\cdot\, ,\tau),\Phi_{1,k}^*+\Phi_{2,k}^*)_{L^2(0,\pi)^2}
		=-(W,(0, 1_{\mathcal{O}}u_2))_{L^2(0,\tau;L^2(\mathcal{O})^2)}
		=-(w_2,u_2)_{L^2(0,T;L^2(\mathcal{O}))}.}
	\end{equation}
	
	\noindent{\textit{Step 4. Volterra equation.}} Let $\theta$ be a solution to \eqref{volterrask} 
	with right-hand side given by $1_{\mathcal{O}}u_2$.  
    
    Taking into account \eqref{volterrask} and Lemma \ref{lemma_volterra2}, 
	we get $K^*\theta=1_{\mathcal{O}}u_2$. Replacing this into \eqref{eq.ff}, we have
	 \begin{equation}\label{eq.ggg}\small{
		(W(\,\cdot\, ,\tau),\Phi_{1,k}^*+\Phi_{2,k}^*)_{L^2(0,\pi)^2}
		=-(w_2,u_2)_{L^2(0,T;L^2(\mathcal{O}))}
		=-(w_2,K^*\theta)_{L^2(0,T;L^2(\mathcal{O}))}.}
	\end{equation}

	From \eqref{def_Kg} and \eqref{sys.2por2heatw}, $y_2=Kw_2$. Therefore,
	\begin{equation}\label{eq.ii}\small{
		(W(\,\cdot\, ,\tau),\Phi_{1,k}^*+\Phi_{2,k}^*)_{L^2(0,\pi)^2}
		=-(w_2,K^*\theta)_{L^2(0,T;L^2(\mathcal{O}))}
		=-(y_2,\theta)_{H^1(0,T;L^2(\mathcal{O}))}}.
	\end{equation}
	
	Finally, putting together \eqref{eq.sol.coeffs.heat}, \eqref{eq.c}, \eqref{eq.ff}, and \eqref{eq.ii}, our reconstruction formula is given by
    \begin{equation*}
        \begin{aligned}
            a(T)(f_1^{\psi_k}+f_2^{\varphi_k}+f_1^{\varphi})+b(T)f_1^{\varphi}
		&=-\sigma(0)(y_2,\theta)_{H^1(0,T;L^2(\mathcal{O}))}\\ 
        &\quad-\displaystyle\int\limits_{0}^{T}\partial_t\sigma(T-s)(y_2,\theta)_{H^1(0,T;L^2(\mathcal{O}))}\,ds,
        \end{aligned}
	\end{equation*}
	where 
		\[a(T)=\sigma(T)\Biggl(1-\frac{k^2}{\sigma(T)}\int\limits_{0}^Te^{-k^2(T-s)}\sigma(s)\,ds\Biggr)\]
	and 
		\[b(T)=-I_k(q)\Biggl(1-k^2\int\limits_{0}^Te^{-k^2(T-s)}\Biggl(\int\limits_{0}^se^{-k^2(s-\tau)}
		\sigma(\tau)\,d\tau\Biggr)\,ds\Biggr).\]
	This completes the proof of Theorem \ref{th3}.	
	\end{proof}

	\begin{remark}
		Note that, the advantage of using the linear combination $\Phi_{1,k}^*+\Phi_{2,k}^*$ in the above proof 
		lies in the fact that we can recover all coefficients for each term of $F=(f_1,\dots,f_n)^t$ through
		subspaces of $L^2(0,\pi)$, i.e.,  for the source $f_1\in L^2(0,\pi)$, we have $L^2=H_1\oplus H_2$, and it admits 
		the representation $f_1=\sum_{k\in\mathbb{N}}f_1^{\varphi_k}\varphi_k+f_1^{\psi_k}\psi_k$, where 
		$H_1=\langle \{\varphi_k: k\in\mathbb{N}\}\rangle$ and $H_2=\langle \{\psi_k: k\in\mathbb{N}\}\rangle$; meanwhile,
		since the coupling occurs in the second equation, for $f_2\in L^2(0,\pi)$, we only obtain $P_{H_1}f_2$, where 
		$P_{H_1}$ represents the orthogonal projector from $L^2(0,\pi)$ onto $H_1$.  Nevertheless, $H_1$ constitutes an
		orthonormal basis of $L^2(0,\pi)$, and therefore the source reconstruction formula \eqref{formula.recoveryth3} shows
		every term of $F$ in complete form.  
	\end{remark}

%% file: diagram.tex
\scalebox{0.8}{
\begin{tikzpicture}[scale=0.3,
    >=Stealth,
    node distance=3.5cm,
    on grid,
    auto,
    thick,
    boxmain/.style={
        draw,
        rectangle,
        align=left,
        text width=7.4cm,
        minimum height=1.5cm,
        %fill=blue!20,
        rounded corners
    },
    box/.style={
        draw,
        rectangle,
        align=left, %center,
        text width=6.4cm,
        minimum height=1.2cm,
        rounded corners
    },
    boxintegral/.style={
        draw,
        rectangle,
        align=left,
        text width=6.1cm,
        minimum height=1.2cm,
        rounded corners
    },
    boxshort/.style={
        draw,
        rectangle,
        align=left, %center,
        text width=5.cm,
        minimum height=1.2cm,
        rounded corners
    }
]

% Main system (top)
\node[boxmain] (main) {
    \hspace{23mm}\textbf{MAIN SYSTEM}\\
    %\flushleft
    $Y(x,t) = (y_{1}(x,t),y_{2}(x,t),\dots, y_{n}(x,t))^{t}$: states \\
    $\sigma(t)$: known source term\\
    $F(x) = (f_1(x),f_2(x),\dots,f_n(x))^{t}:$ unknown
};

% Spectral analysis (right)
\node[boxshort, below =of main] (spec) {
    \hspace{2mm}\textbf{Step 1B. Spectral analysis}\\
    $Y(x,t) = \sum_{k\in\mathbb{N}} \alpha_k(t)\phi_k(x)$ \\
    $\cdot$ $\alpha_k$: appropriate functions to\\
    \hspace{1mm} be determined \\
    $\cdot$ $\phi_k$: appropriate basis \\ \hspace{1mm} functions to be chosen
};

% Integral representation (left)
\node[boxintegral, left=of spec,xshift=-3.0cm] (intrep) {
    \hspace{2mm}\textbf{Step 1A. Integral representation}\\
    $Y(x,t) = \int_{0}^{t}\sigma(s) W(x,t -s )\,ds = KW$ \\
    $\cdot$ $W$ satisfies a coupled homogeneous\\
    \hspace{1mm} system with initial datum $\sigma(0)F(x)$, \\
     \hspace{1mm} $x\in\Omega$\\
    $\cdot$ $K$ is an appropriate operator
};

% coeff (right)
\node[boxshort, right=of spec, xshift=2.5cm] (coeff) {
    Determine the coefficients $\alpha_k(t)$,$t\in (0,T)$, $k\in\mathbb{N}$, using the main system and properties of $\{\phi_k\}_{k\in\mathbb{N}}$, $x\in\Omega$ for solving the system of ODEs.
};

% reconstruction (below spec)
\node[box, below=of spec] (recons) {
    %\vspace*{-12pt}
    %\begin{itemize}
    \hspace{3mm}\textbf{Reconstruction with global terms}\\
    Reconstruction formula with global terms in $\Omega\times (0,T)$ for each coefficients $f_j^k = (f_j, \phi_k)_{L^{2}(\Omega)}$, $j=1,\dots,n$, $k\in\mathbb{N}$ \\
    $\cdot$ Measurements by $Y|_{\Omega\times (0,T)}$ and \\ \hspace{1mm} $\partial_t Y|_{\Omega\times (0,T)}$
};

% Volterra equations (below spectral)
\node[boxshort, below=of recons] (volt) {
    %\hspace{0mm} 
    \shortstack[c]{\hspace{1mm}\textbf{Step 3. Volterra equations}}\\
    Use the control function $U(x,t)$ in $\mathcal{O}\times (0,T)$ (control domain) given by the solution of the null controllability problem as a datum for the Volterra system $K^{*}\Theta = U$, $\Omega\times(0,T)$\\
    ($K^{*}$ dual operator of $K$) 
};

% Internal null controllability (below integral)
\node[boxshort, left=of volt, xshift=-2.5cm] (null) {
    \hspace{6mm}\shortstack[c]{\textbf{Step 2. Internal null}\\\textbf{controllability}}\\
    Use the null controllability property on the adjoint system associated to the system defined by $W$, with initial datum $\phi_k(x)$, $x\in\Omega$, $k\in\mathbb{N}$.
};

% Internal null controllability (below integral)
\node[boxshort, right=of volt, xshift=2.5cm] (final) {
    \hspace{7mm}\textbf{Final reconstruction}\\
    Reconstruction formula with local terms in $\mathcal{O}\times (0,T)$ for each coefficients $f_j^k = (f_j, \phi_k)_{L^{2}(\Omega)}$, $j=1,\dots, n$, $k\in \mathbb{N}$. \\
    $\cdot$ Measurements by $y_{n}|_{\mathcal{O}\times (0,T)}$ \\
    \hspace{1mm} and $\partial_t y_{n}|_{\mathcal{O}\times (0,T)}$
};
% Arrows
\draw[->, shorten >=1mm] (main) -| (intrep);
\draw[->, shorten >=1mm] (main) -- (spec);
\draw[->, shorten >=1mm] (spec) -- (coeff);
\draw[->, shorten >=1mm] (intrep) |- (recons);
\draw[->, transform canvas={shift={(0,1mm)}}, shorten <=1mm, shorten >=1mm] (coeff) |- (recons);
\draw[->, shorten >=1mm] (null) -- (volt);
\draw[->, shorten >=1mm] (volt) -- (final);
\draw[->, transform canvas={shift={(0,-1mm)}}, shorten >=2mm] (recons.east) -| (final);
\end{tikzpicture}}

%% file: 02_numerics.tex
\section{Numerical results}	\label{section.numerics}

In this section, we present numerical tests related to Problems \ref{Problem1} and \ref{Problem2} proved in Theorems \ref{th1} and \ref{th3}, respectively. Since our theoretical results proved in previous sections are in $\mathbb{R}^{n}$, with $n>1$  when the coupling is given for a constant matrix $Q$ and $n=1$ for the case where the coupling has space-dependent coefficients on the matrix $Q$, we perform 1D and 2D numerical experiments under different configurations.

The following section presents an optimization-based methodology designed for the complete numerical analysis, in contrast to the framework discussed in the previous section, which is grounded in spectral analysis, the null controllability property, and Volterra equations. Further details on this approach are given in subsection ~\ref{sec:impl_details}. In concordance with Problems \ref{Problem1} and \ref{Problem2}, we recover the source term $F(x)=(f_1(x),f_2(x))^t$ as the solution of the following minimization problem: 
\begin{equation}\label{minproblem.functional}
        \min_{F \in L^2(\Omega)}J(Y,F)=\frac{1}{2}\iint\limits_{\Omega\times(0,T)} |\sigma(t)F(x)|^2 %dxdt
        +\frac{k}{2}\sum\limits_{i=1, i\neq j}^2\,\,\iint\limits_{\mathcal{O}\times(0,T)} 
        |y_i - y_{i,\text{obs}}|^2  + |\partial_t y_i - \partial_t y_{i,\text{obs}}|^2 %\,dx \,dt
\end{equation}
    subject to 
    \begin{equation}\label{minproblem.system}
    \left\{
    \begin{array}{llll}
        \partial_t Y- \nu\Delta Y + Q Y=\sigma(t)F(x), &\mbox{in} &\Omega\times(0,T),\\
        Y=0, &\mbox{on}&\partial\Omega\times(0,T),\\        
		Y(\,\cdot\, ,0)=0, &\mbox{in}&\Omega.        
    \end{array}\right.
    \end{equation}
    where $Y_{\text{obs}} := (y_{1,\text{obs}},y_{2,\text{obs}})^t$ are observation data on the subdomain $\mathcal{O}\times(0,T)$ (or also so-called partial measurements of the state $Y=(y_1,y_2)^t$). We also define %$Q$ is the coupling matrix given by
    
    $$Q= \begin{pmatrix} q_{11} & q_{12} \\ q_{21} & q_{22} \end{pmatrix},\quad \text{and} \quad \sigma(t):= \begin{cases} 1+\frac{1}{2}\cos\left(\frac{4\pi t}{T - t_0}\right) & t<T-t_0 \\ \frac{3}{2} & t\geq T-t_0\end{cases},$$ 
    where $Q$ is the coupling matrix and $\sigma$ is the function selected for all experiments (see \cite{2013GOT, 2017GMO}). 
    The functional $J$ is selected in this form based on the theoretical requirements for optimization problems and following ideas from \cite{2017GMO}. Additionally, the Lipschitz-type stability result presented in Appendix \ref{appendix.lipschitz} allows us to clarify the structure of the function $J(Y,F)$ defined in \eqref{minproblem.functional}. The PDE constraining the above minimization problem~\eqref{minproblem.system} is solved using finite elements with conforming piecewise linear Lagrange finite element basis functions in 1D and 2D domains. The finite element schemes are solved using the FEniCS library~\cite{AlnaesEtal2015}. For minimizing the functional~\eqref{minproblem.functional}, we employed a steepest descent method \cite{2011burden,2006Nocedal}, with a penalty term $k=10^{5}$ and a fixed step size $10^{-4}$.
 
\subsection{Source reconstruction in 1D}
We consider $\Omega=(0,1)$, $T=0.5$, and $\nu=0.1$. In this case, we considered 100 elements, a time step $\Delta t=0.001$ for temporal discretization, and $k=10^5$. We show numerical reconstruction for two different sources: 
\begin{enumerate}
    \item [a)] $F_1(x):=(\sin(2\pi x),-\sin(2\pi x))^t$.
    \item [b)] $F_2(x):=(f_{2,1}(x),f_{2,2}(x))^t$ defined as follows:
        \begin{equation}
        \nonumber
        f_{2,1}(x):= 
        \begin{cases}
        8(x - 0.1), & \text{ if } x \in(0.1,0.35), \\
        8(0.6 - x), & \text{ if } x \in(0.35,0.6),\\
        0, & \text{ other case},
        \end{cases}
        \quad\mbox{and}\quad f_{2,2}(x):=f_1(x-0.3)
    \end{equation}
\end{enumerate}

\begin{figure}[htp]
    \centering
    \begin{subfigure}{.3\textwidth}
    \includegraphics[width=1.0\linewidth]{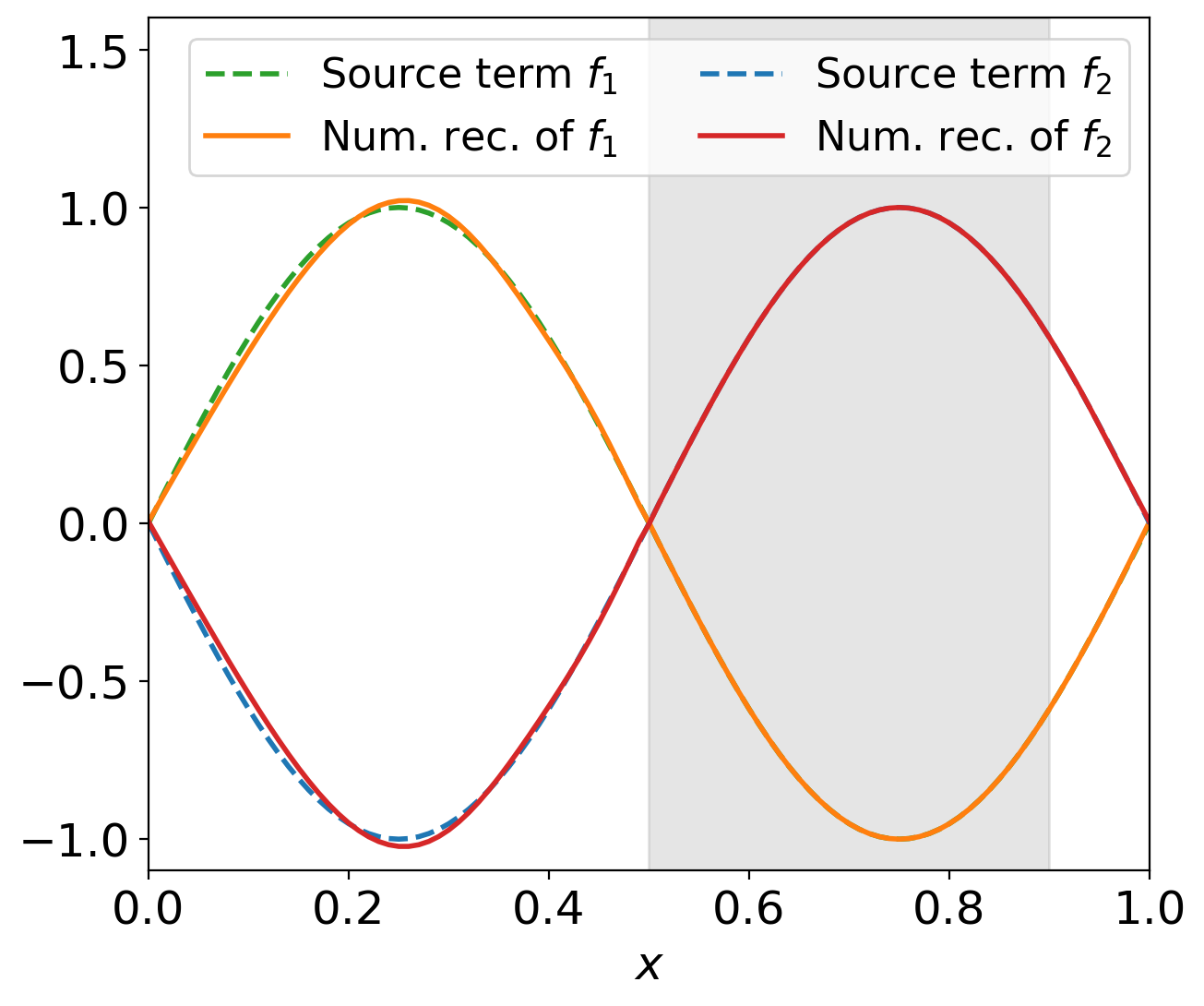}
    \caption{rel. error $2,3\%$}
    \end{subfigure}
    \begin{subfigure}{.3\textwidth}
    \includegraphics[width=1.0\linewidth]{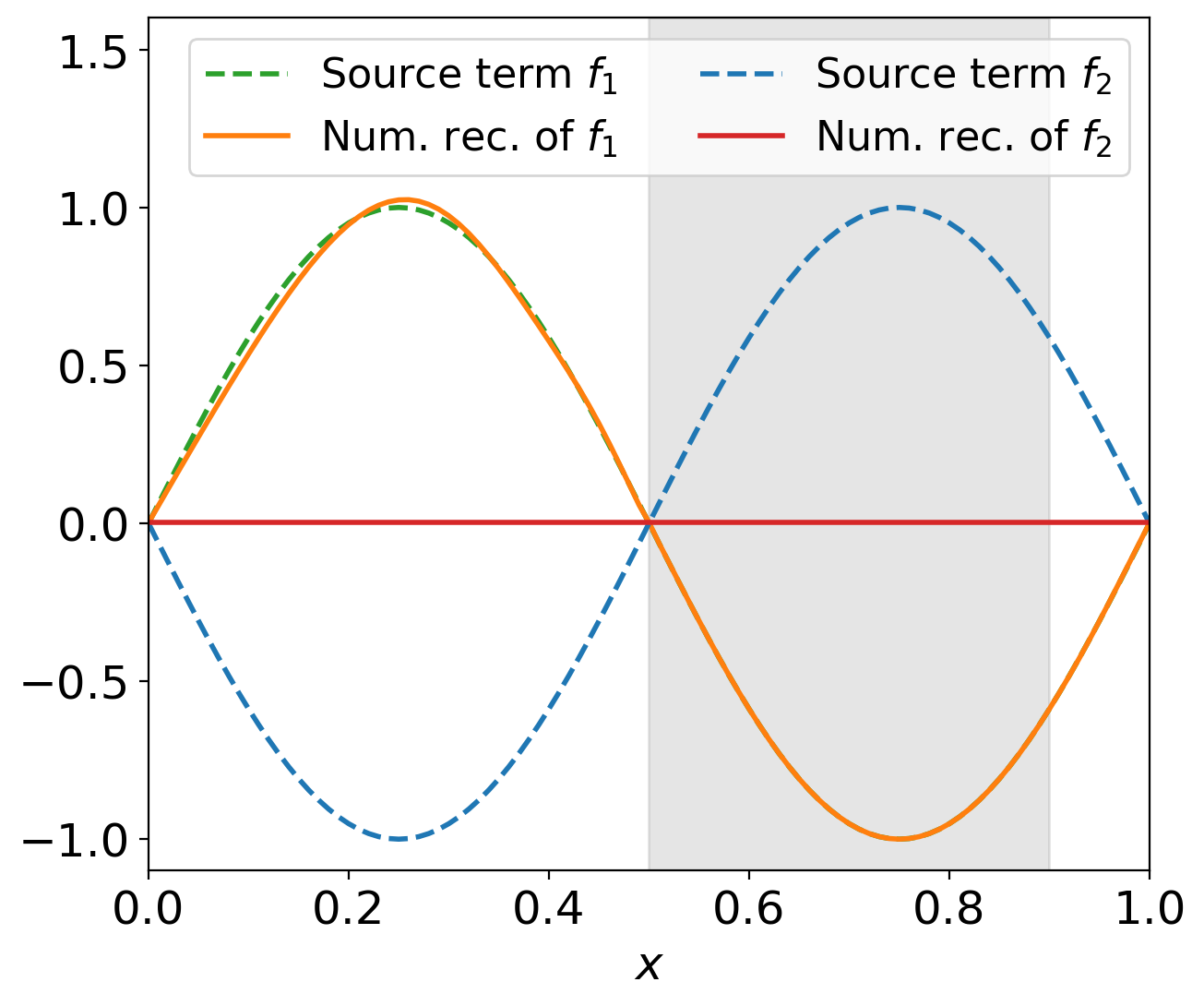}
    \caption{rel. error $70,7\%$}
    \end{subfigure}
    \begin{subfigure}{.3\textwidth}
    \includegraphics[width=1.0\linewidth]{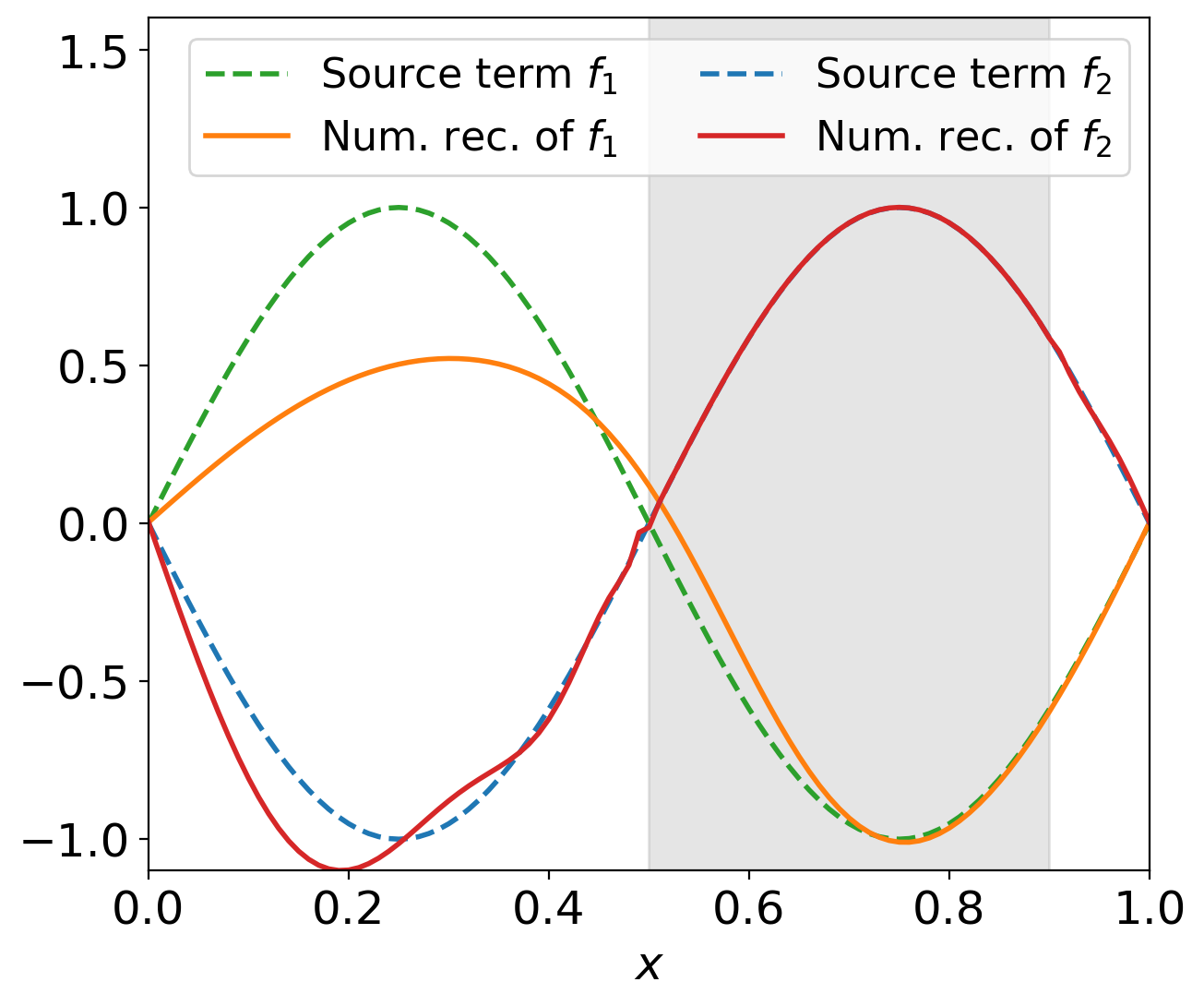}
    \caption{rel. error $25,7\%$}
    \end{subfigure}
    \\
    \begin{subfigure}{.3\textwidth}
    \includegraphics[width=1.0\linewidth]{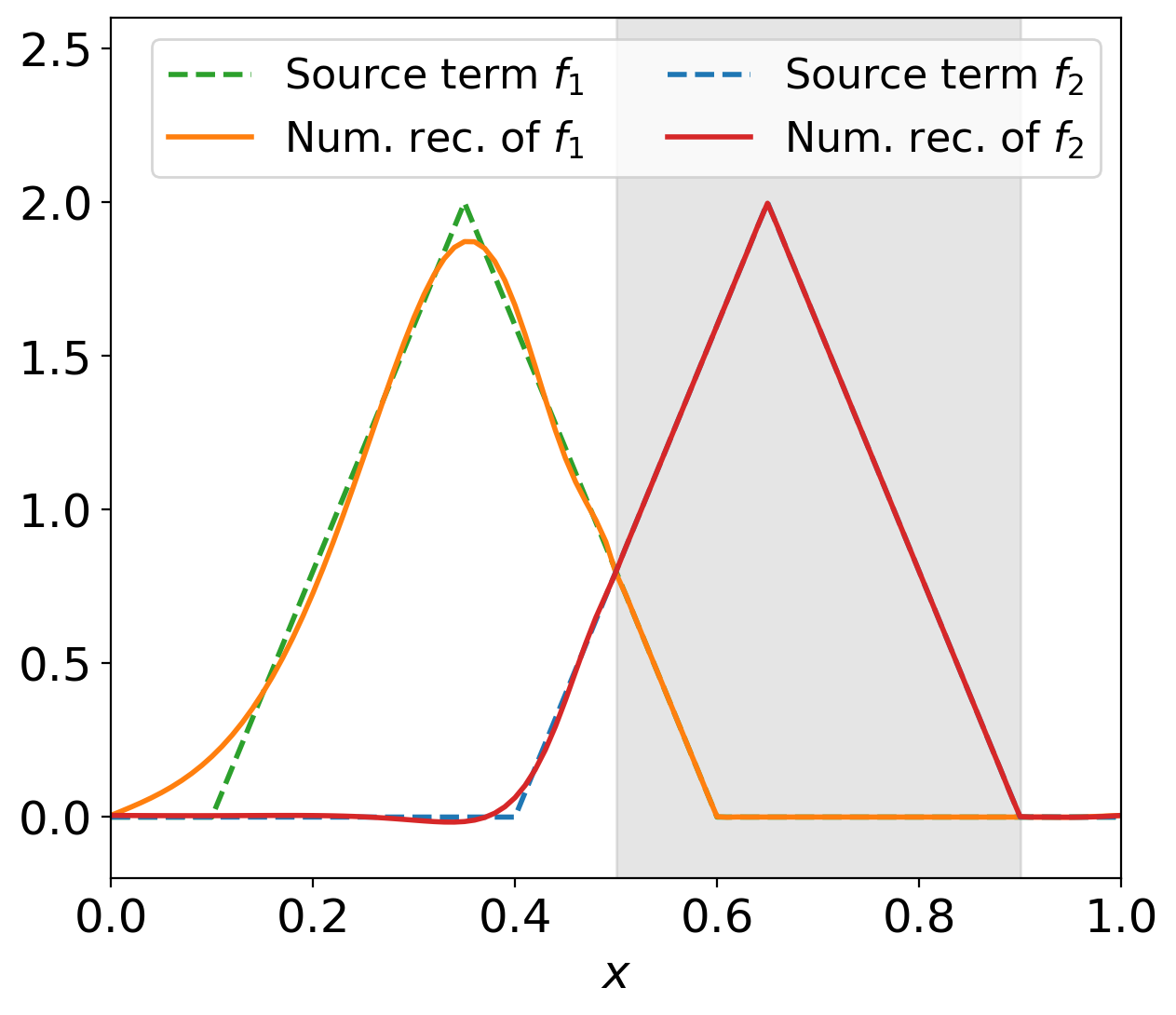}
    \caption{rel. error $4,2\%$}%$,7\%$}
    \end{subfigure}
    \begin{subfigure}{.3\textwidth}
    \includegraphics[width=1.0\linewidth]{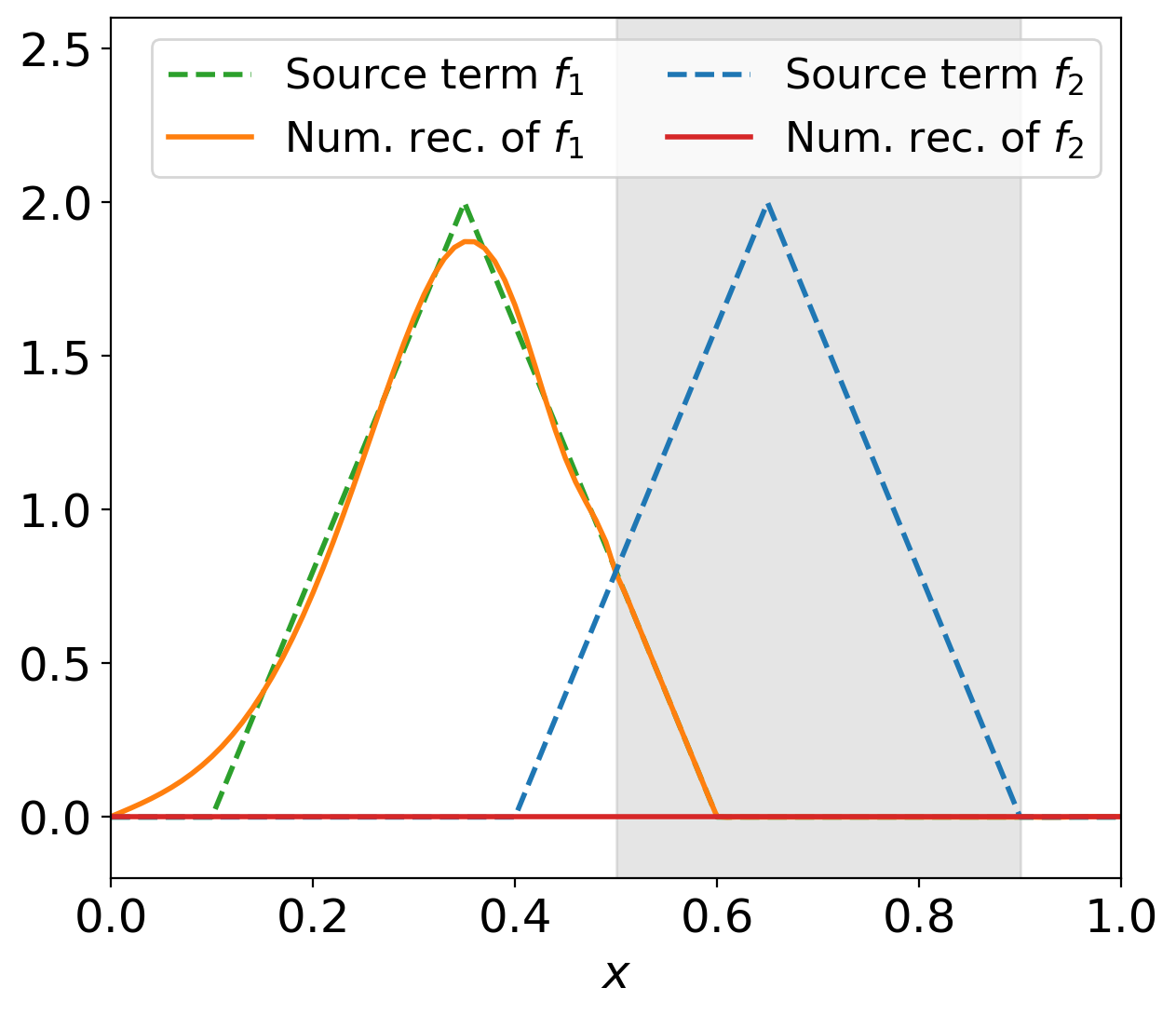}
    \caption{rel. error $70,7\%$}
    \end{subfigure}
    \begin{subfigure}{.3\textwidth}
    \includegraphics[width=1.0\linewidth]{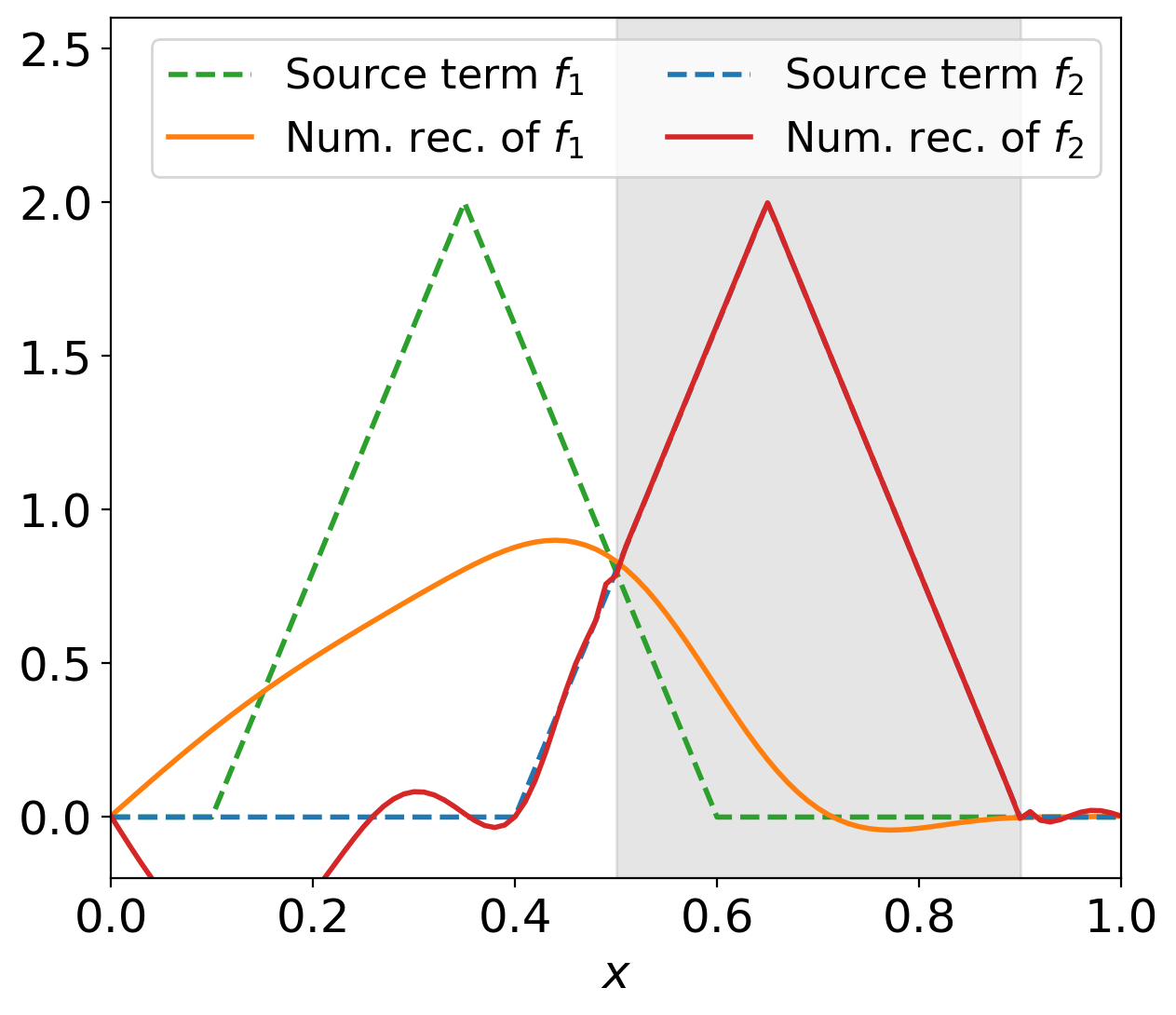}
    \caption{rel. error $38,1\%$}
    \end{subfigure}
    \\
    \begin{subfigure}{.3\textwidth}
    \includegraphics[width=1.0\linewidth]{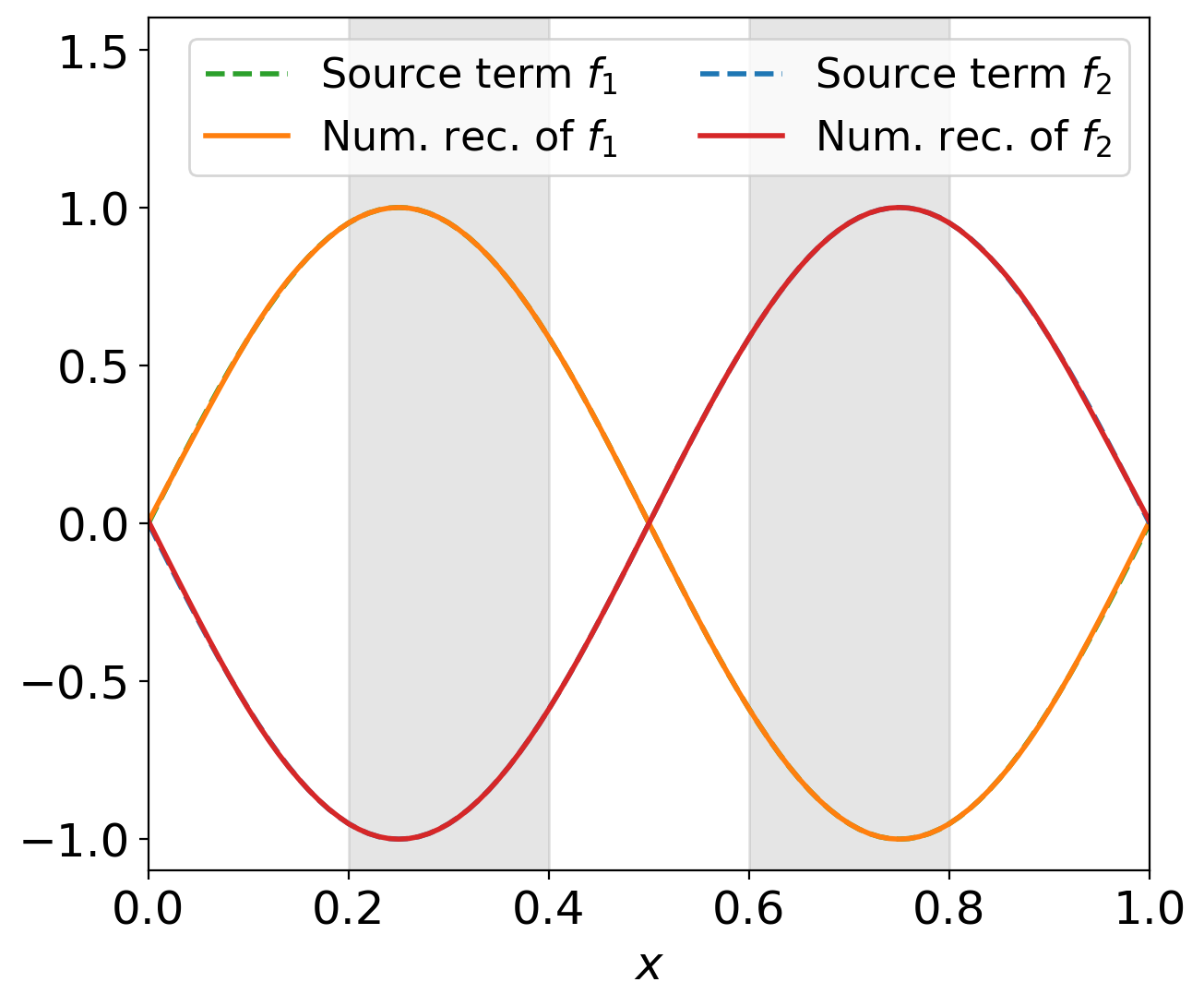}
    \caption{rel. error $0,2\%$}%$2,9\%$}
    \end{subfigure}
    \begin{subfigure}{.3\textwidth}
    \includegraphics[width=1.0\linewidth]{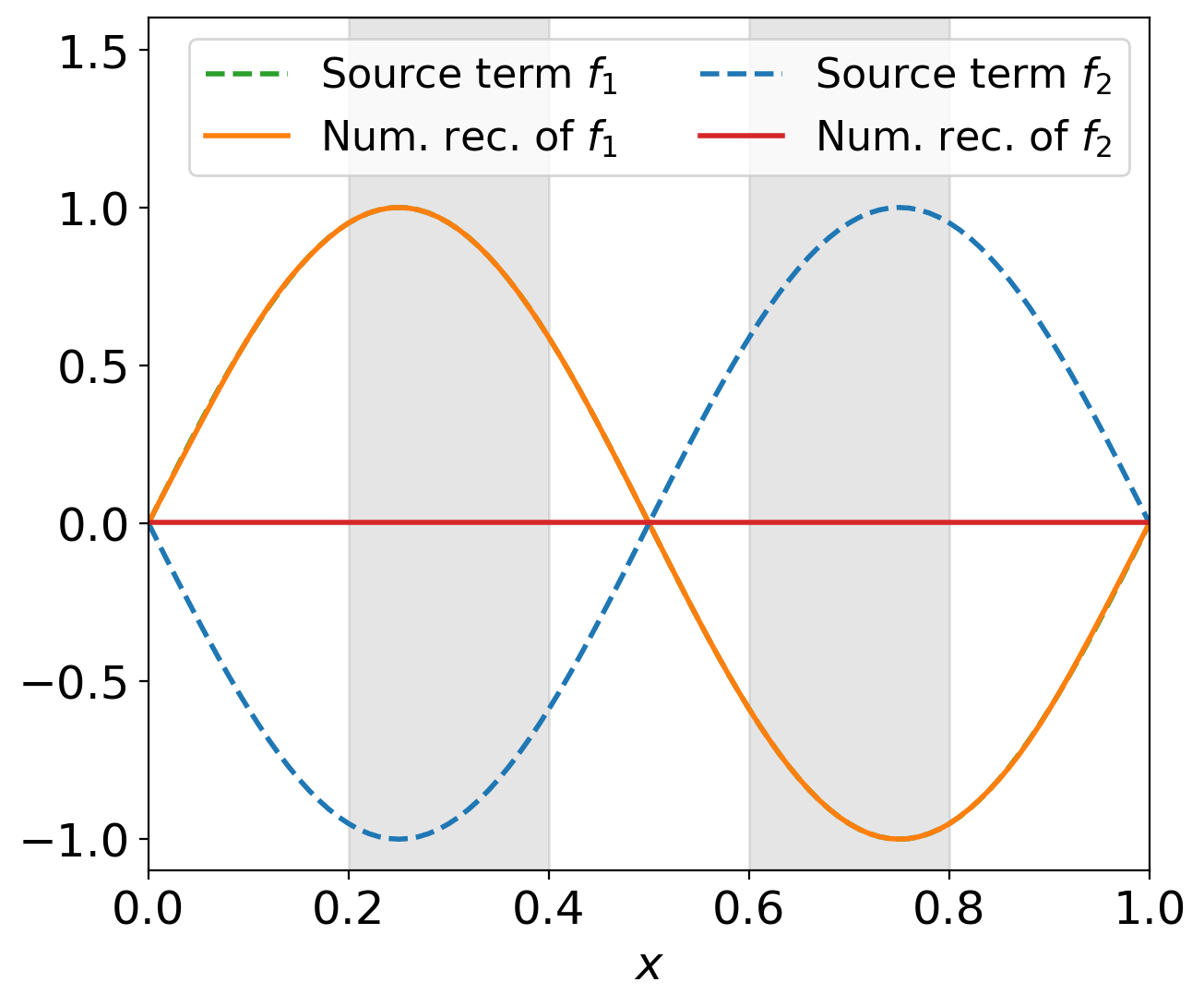}
    \caption{rel. error $70,7\%$}%$24,9\%$}
    \end{subfigure}
    \begin{subfigure}{.3\textwidth}
    \includegraphics[width=1.0\linewidth]{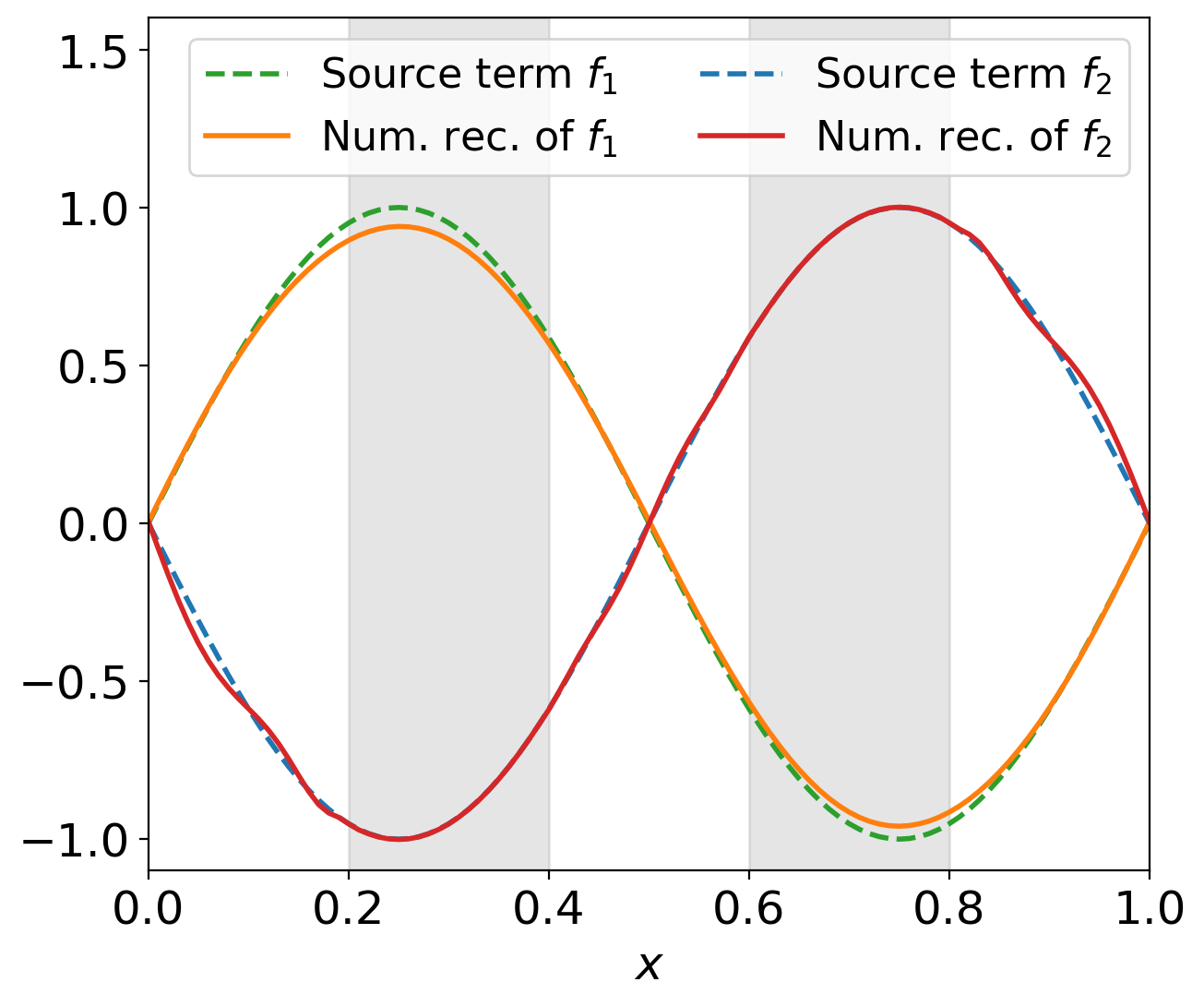}
    \caption{rel. error $3,8\%$}%$41,1\%$} %$42,9\%$}
    \end{subfigure}
    \\
    \begin{subfigure}{.3\textwidth}
    \includegraphics[width=1.0\linewidth]{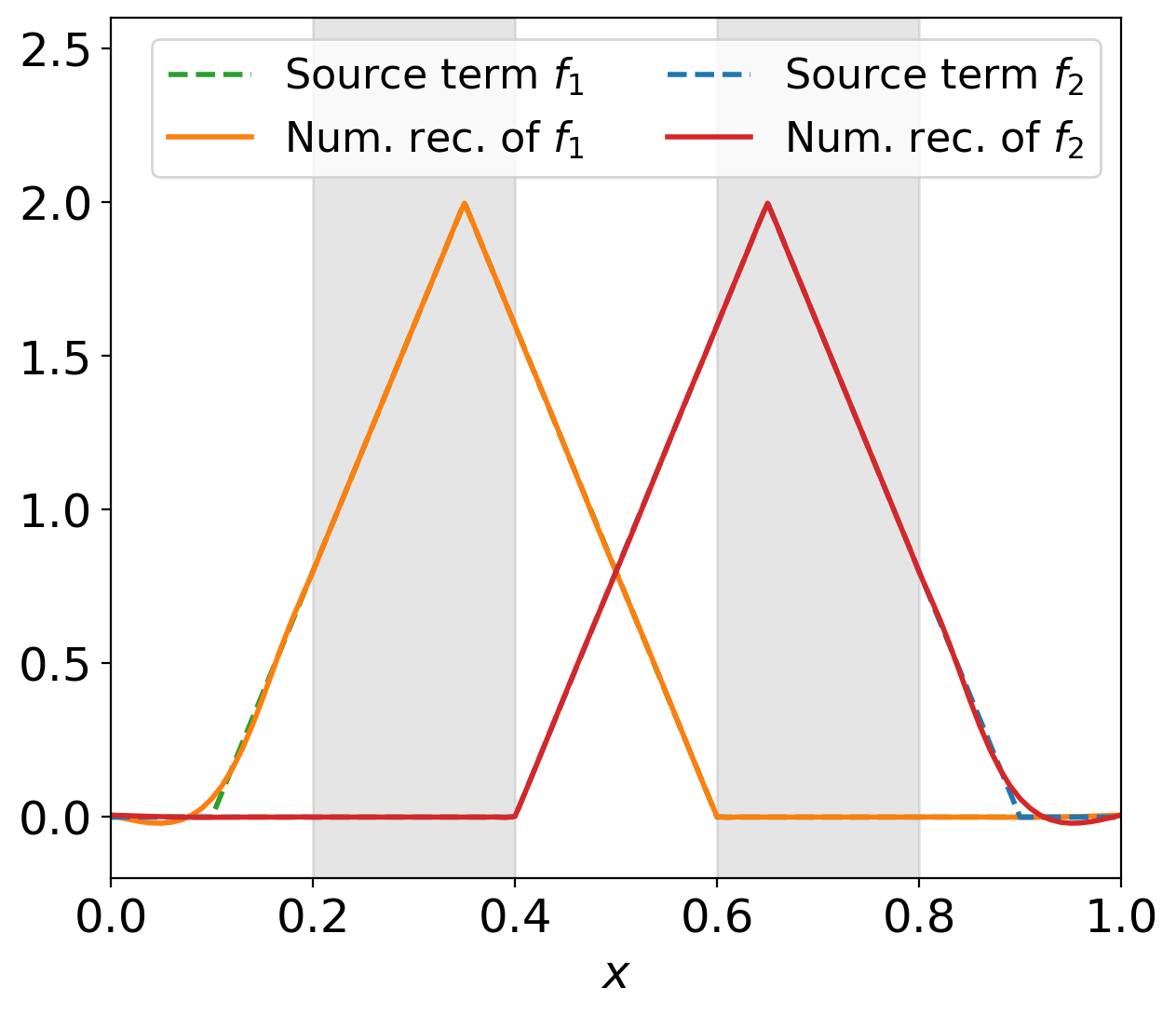}
    \caption{rel. error $1,0\%$}
    \end{subfigure}
    \begin{subfigure}{.3\textwidth}
    \includegraphics[width=1.0\linewidth]{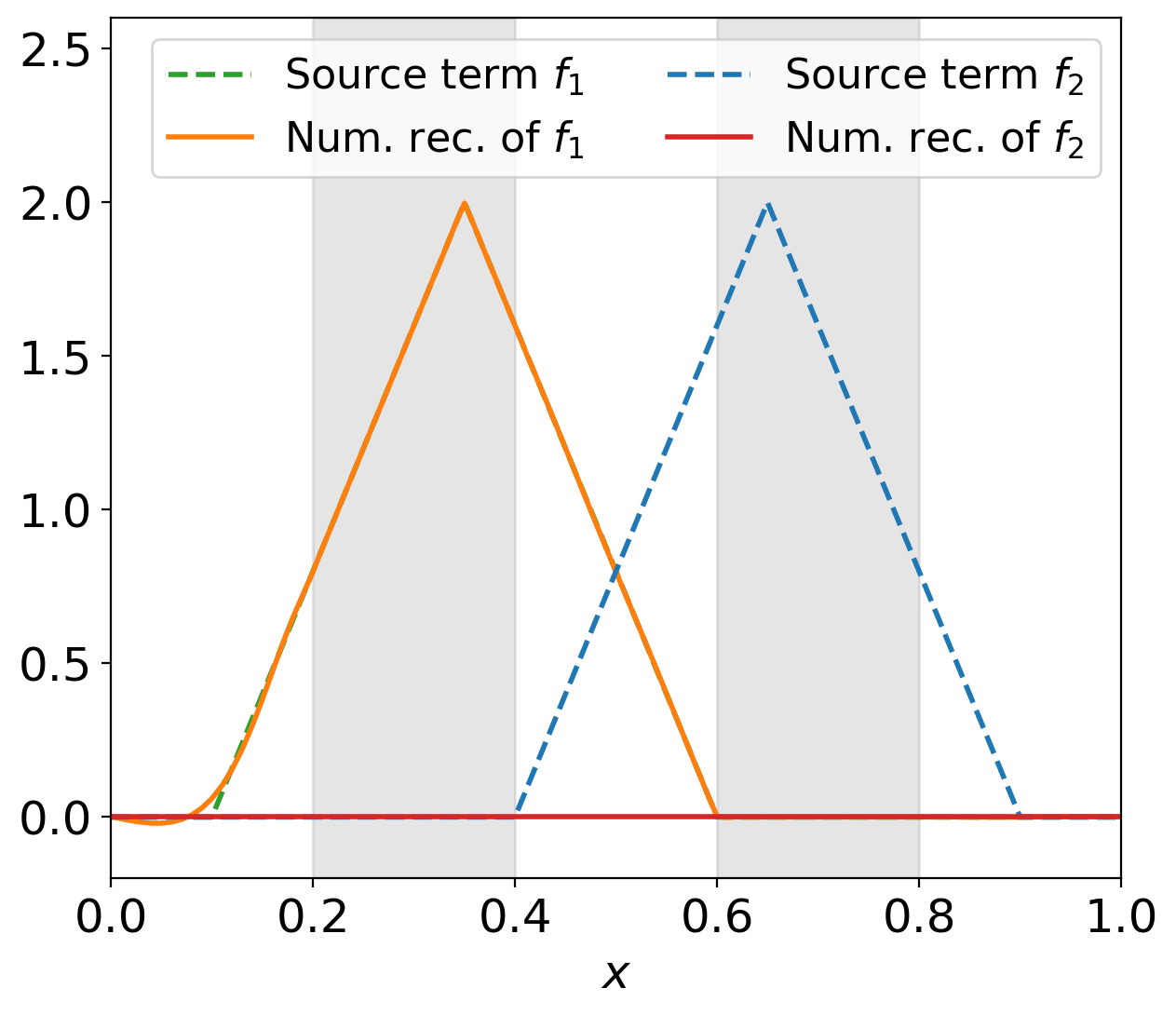}
    \caption{rel. error $70,6\%$}%$21,4\%$}
    \end{subfigure}
    \begin{subfigure}{.3\textwidth}
    \includegraphics[width=1.0\linewidth]{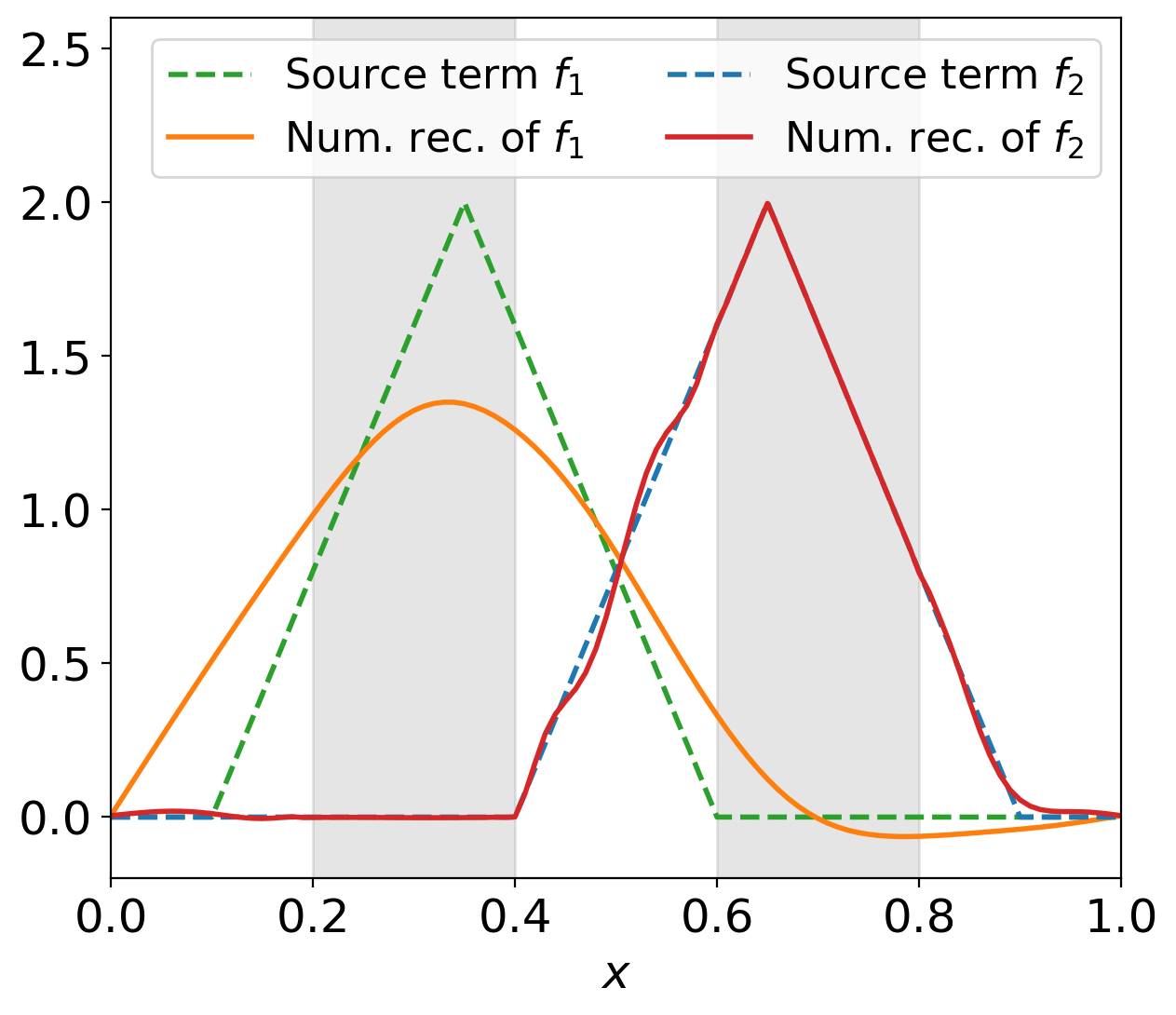}
    \caption{rel. error $20,9\%$}%$44,9\%$}
    \end{subfigure}
    \caption{Source reconstructions when the coupling matrix is given by $q_{11}(x)=q_{12}(x)=q_{22}(x)=0$ and $q_{21}(x)=-x^3 + 4x^2 - 3x + 1$. First column is the reconstruction of $F_1(x)$ and $F_2(x)$ observing both components, i.e., $(y_{1,\text{obs}},y_{2,\text{obs}})^t$. The second column shows the reconstruction of $F_1(x)$ and $F_2(x)$ observing the first component $y_{1,\text{obs}}$, and the last one is the reconstruction observing the second component $y_{2,\text{obs}}$ only.
    The highlighted gray areas correspond to the observation domain $\mathcal{O}$.}
    \label{fig1D_reconstruction_gamma_cubic}
\end{figure}

Figure~\ref{fig1D_reconstruction_gamma_cubic} illustrates reconstructions obtained using a coupling matrix defined by $q_{11}(x)=q_{12}(x)=q_{22}(x)=0$ and $q_{21}(x)=-x^3 + 4x^2 - 3x + 1$ and and two different observation domains, namely $\mathcal{O}_1= (0.5, 0.9)$ and $\mathcal{O}_2= (0.2, 0.4)\cup(0.6,0.8)$. In this setting, reconstructions using measurements from both components consistently exhibit robust performance and relatively low reconstruction errors. On the other hand, reconstructions based on single-component measurements appear more sensitive to both the location and size of the observation subdomain, leading in some cases to larger reconstruction errors, particularly outside the observation region. 
Note that, in the second column, the reconstruction of the source term $f_2$ from the first state does not yield a satisfactory error. 
This is because, given the structure of $Q$, only the second equation contains a coupling term, so the first equation ``does not see'' the second component.
This observation, arising from the numerical experiments, also highlights open problems in the literature on inverse problems, namely, questions regarding the optimal location of the observation subdomain and the possibility of determining an optimal matrix configuration for the coupling.

\begin{figure}[htp]
    \centering
    \begin{subfigure}{.3\textwidth}
    \includegraphics[width=1.0\linewidth]{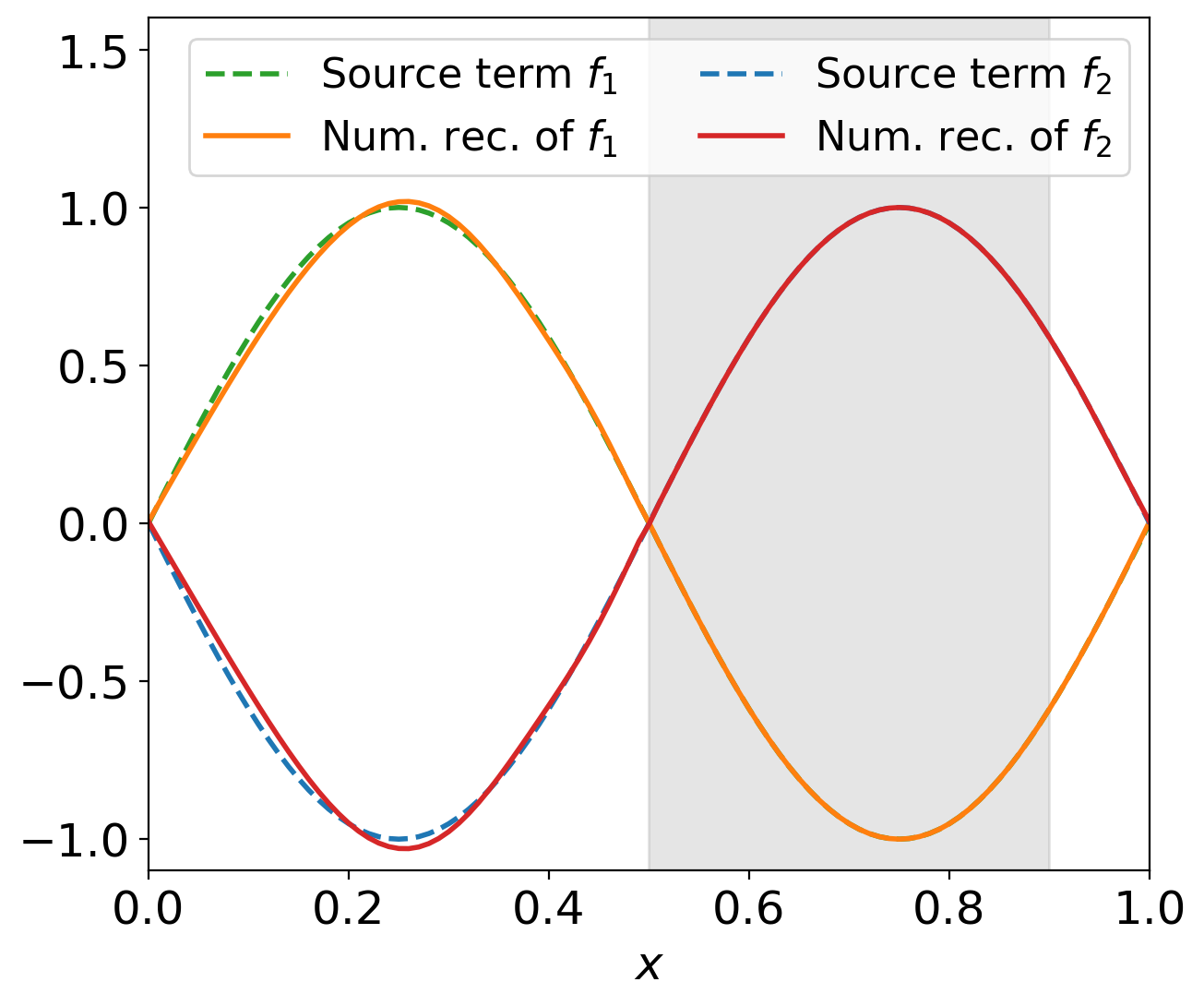}
    \caption{rel. error $2,6\%$}
    \end{subfigure}
    \begin{subfigure}{.3\textwidth}
    \includegraphics[width=1.0\linewidth]{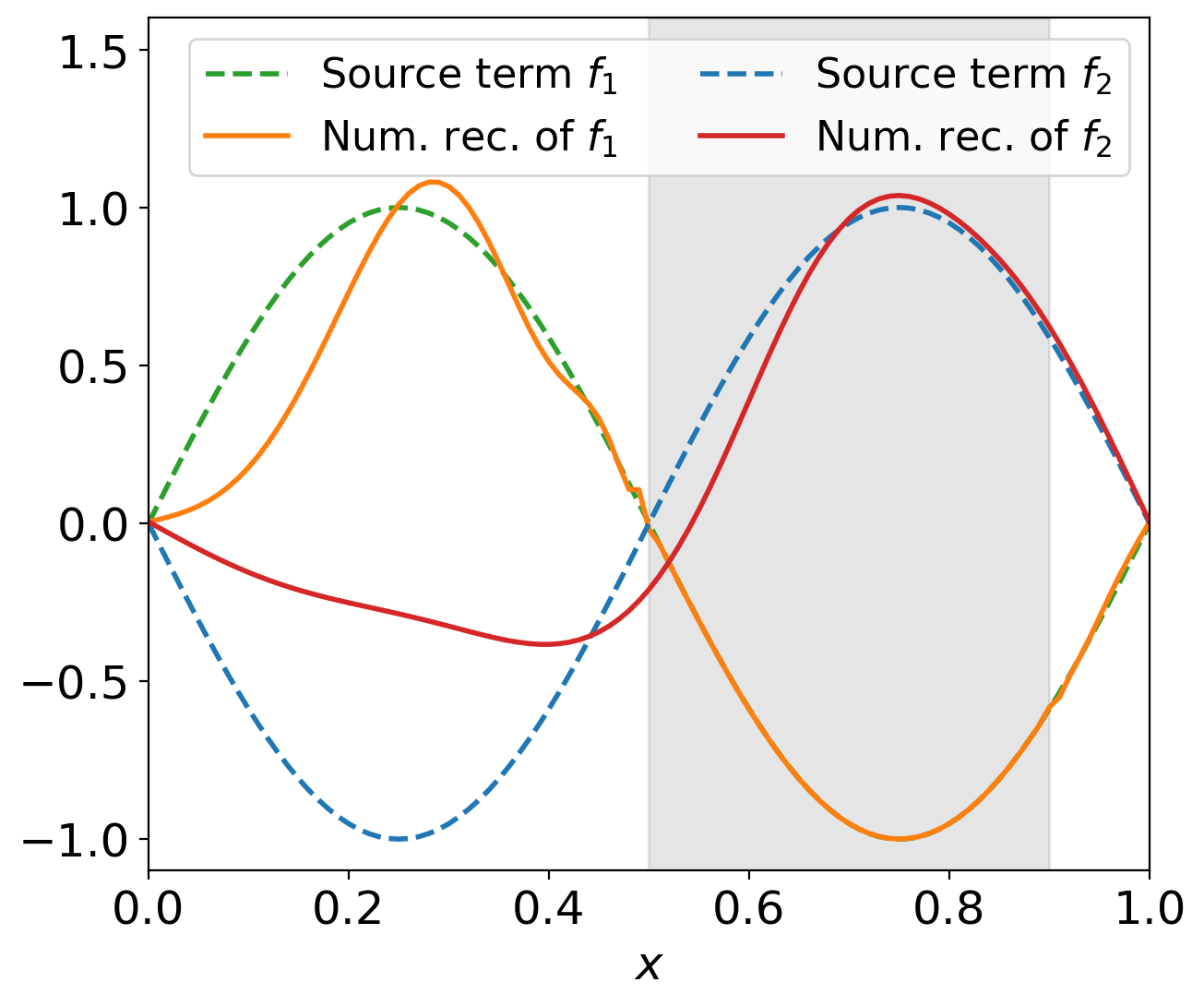}
    \caption{rel. error $37,7\%$}
    \end{subfigure}
    \begin{subfigure}{.3\textwidth}
    \includegraphics[width=1.0\linewidth]{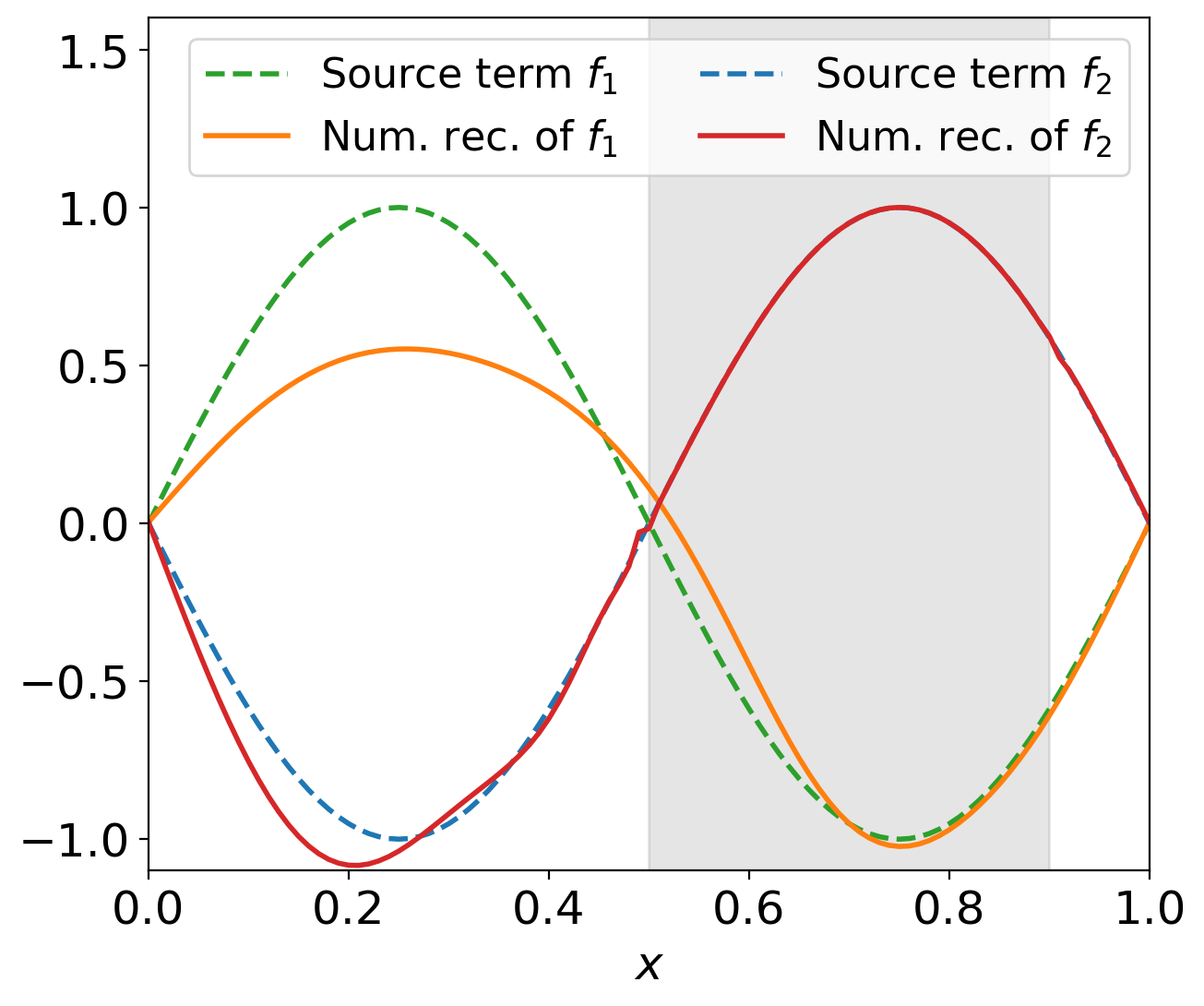}
    \caption{rel. error $22,8\%$}
    \end{subfigure}
    \\
    \begin{subfigure}{.3\textwidth}
    \includegraphics[width=1.0\linewidth]{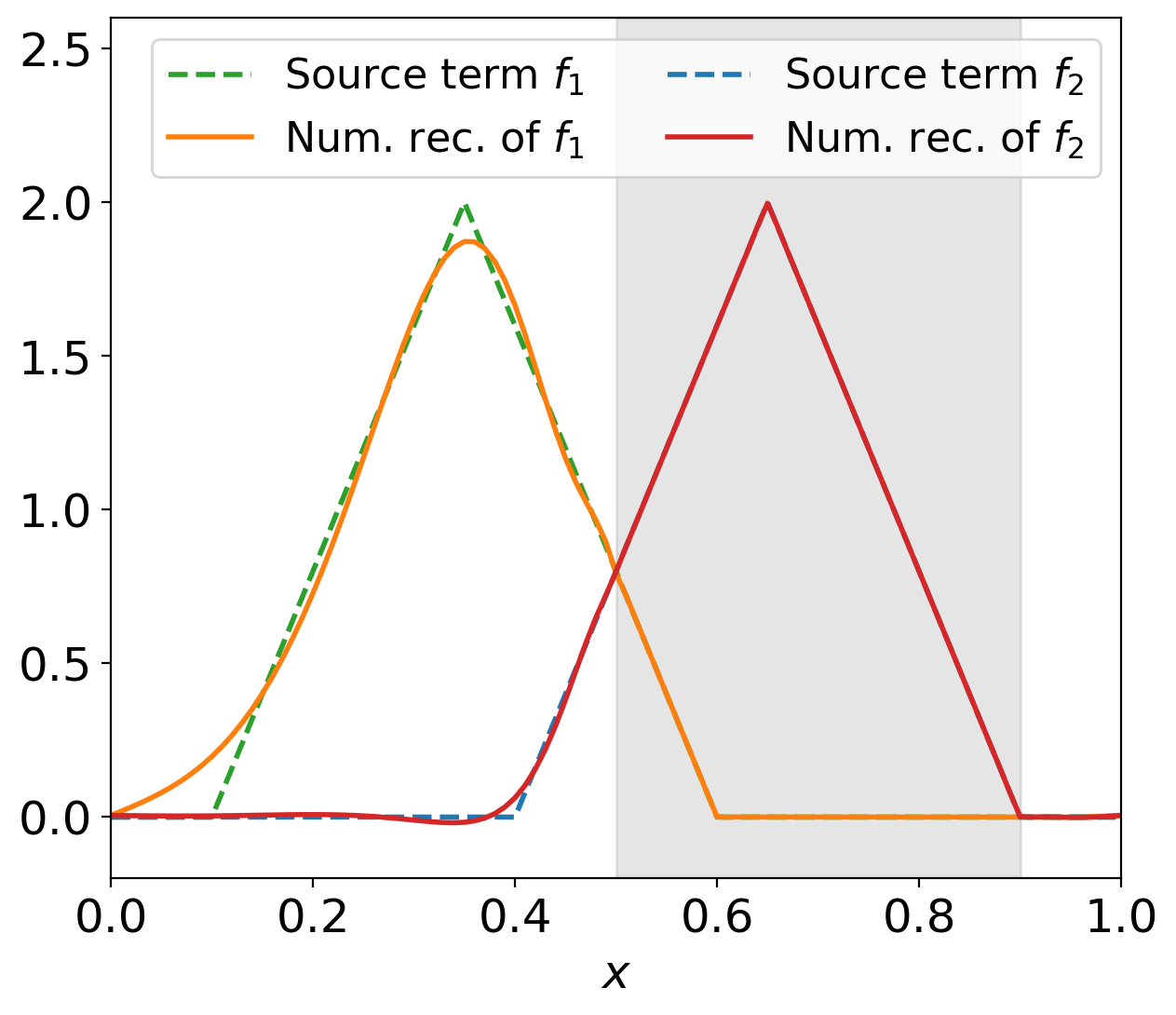}
    \caption{rel. error $4,3\%$}
    \end{subfigure}
    \begin{subfigure}{.3\textwidth}
    \includegraphics[width=1.0\linewidth]{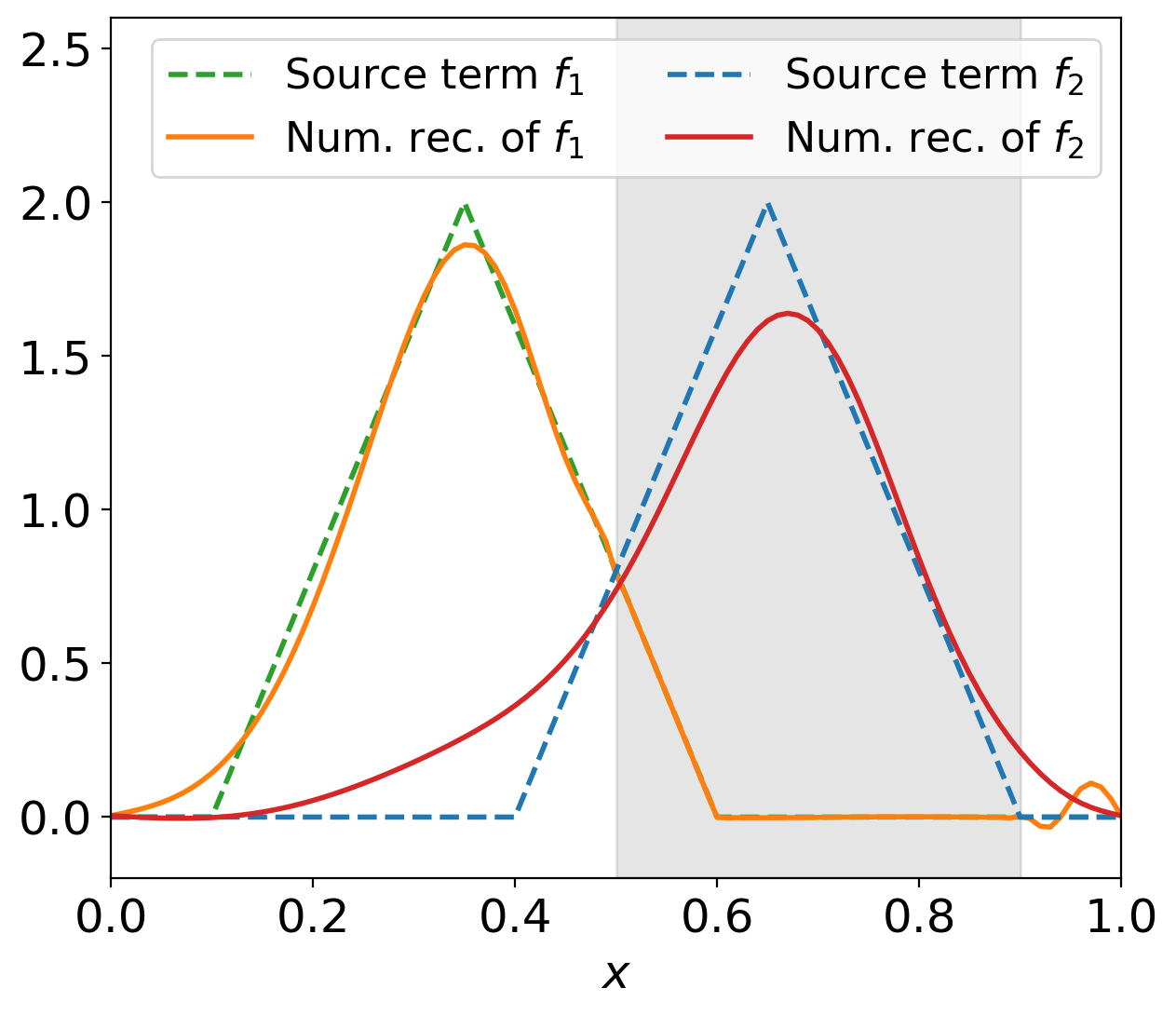}
    \caption{rel. error $13,7\%$}
    \end{subfigure}
    \begin{subfigure}{.3\textwidth}
    \includegraphics[width=1.0\linewidth]{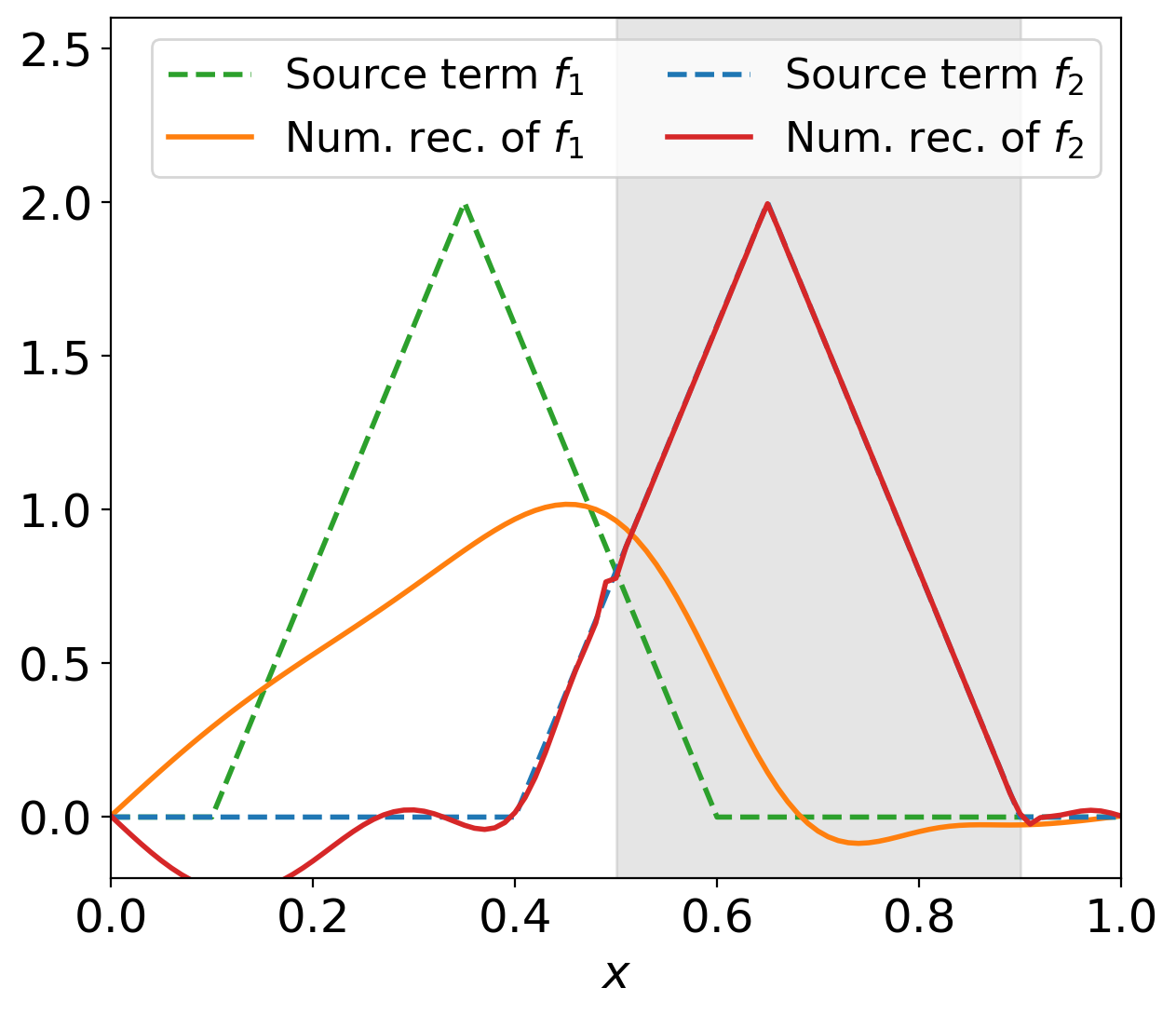}
    \caption{rel. error $35,1\%$}
    \end{subfigure}
    \\
    \begin{subfigure}{.3\textwidth}
    \includegraphics[width=1.0\linewidth]{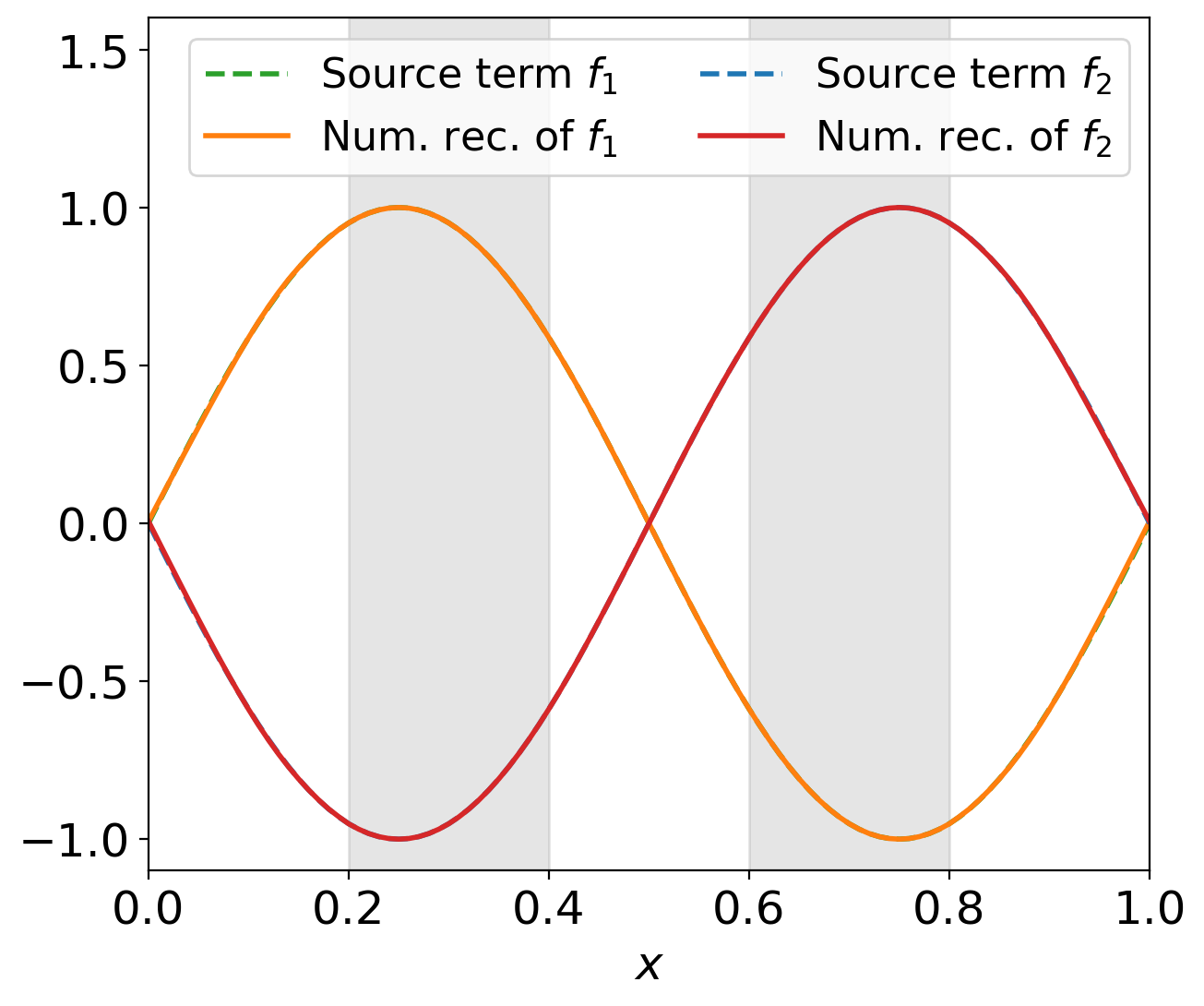}
    \caption{rel. error $0,3\%$}
    \end{subfigure}
    \begin{subfigure}{.3\textwidth}
    \includegraphics[width=1.0\linewidth]{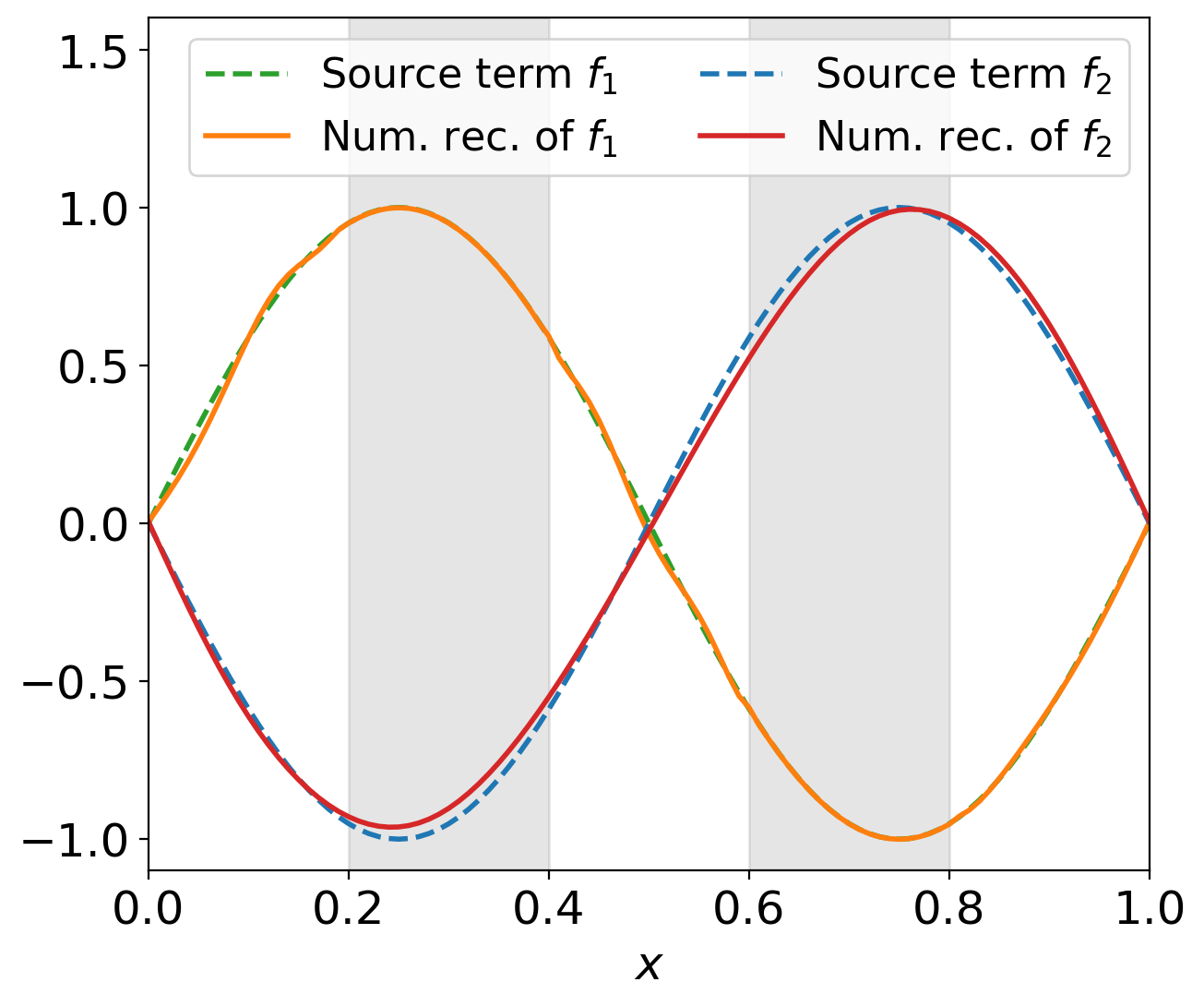}
    \caption{rel. error $3,9\%$}
    \end{subfigure}
    \begin{subfigure}{.3\textwidth}
    \includegraphics[width=1.0\linewidth]{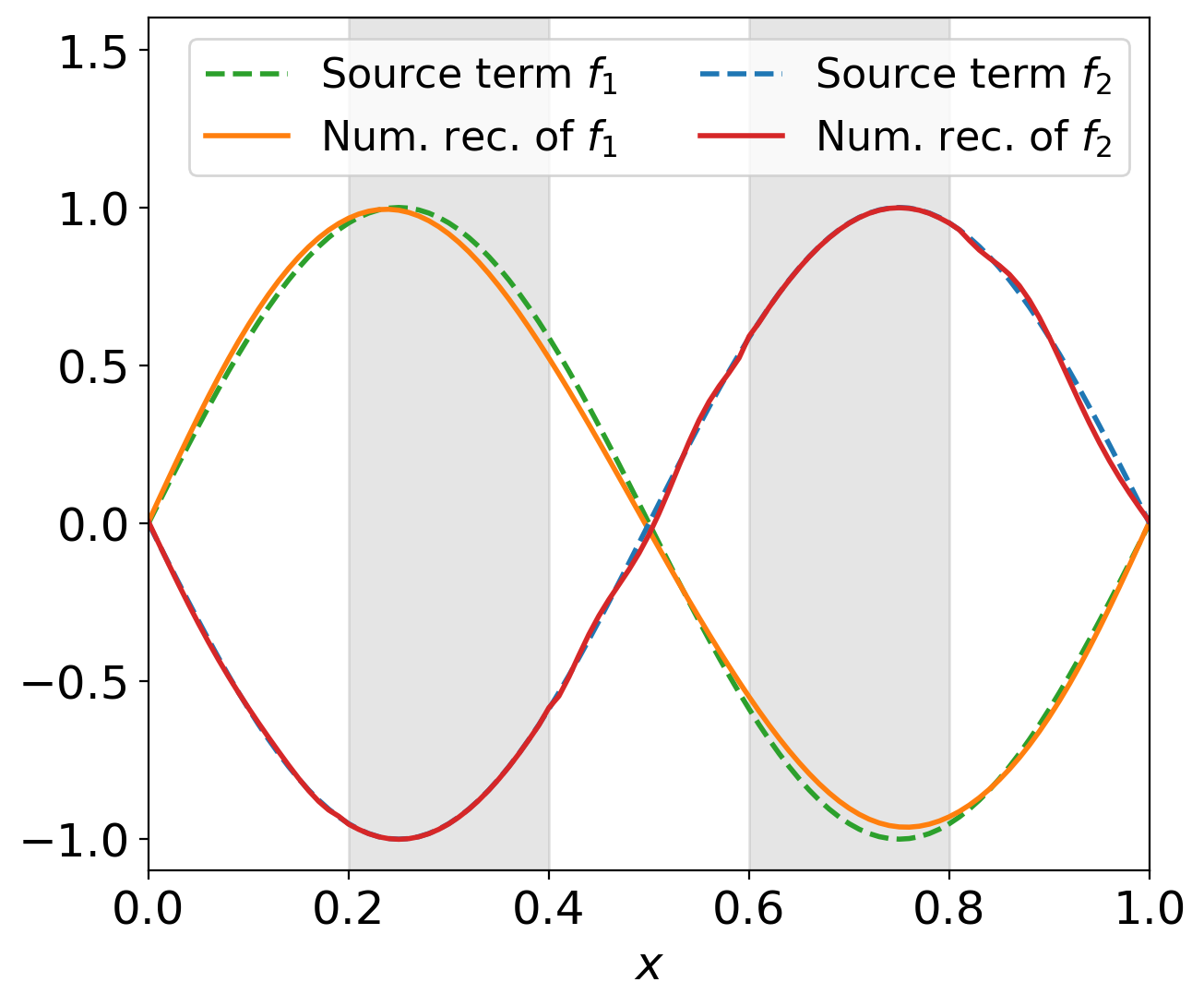}
    \caption{rel. error $3,9\%$}
    \end{subfigure}
    \\
    \begin{subfigure}{.3\textwidth}
    \includegraphics[width=1.0\linewidth]{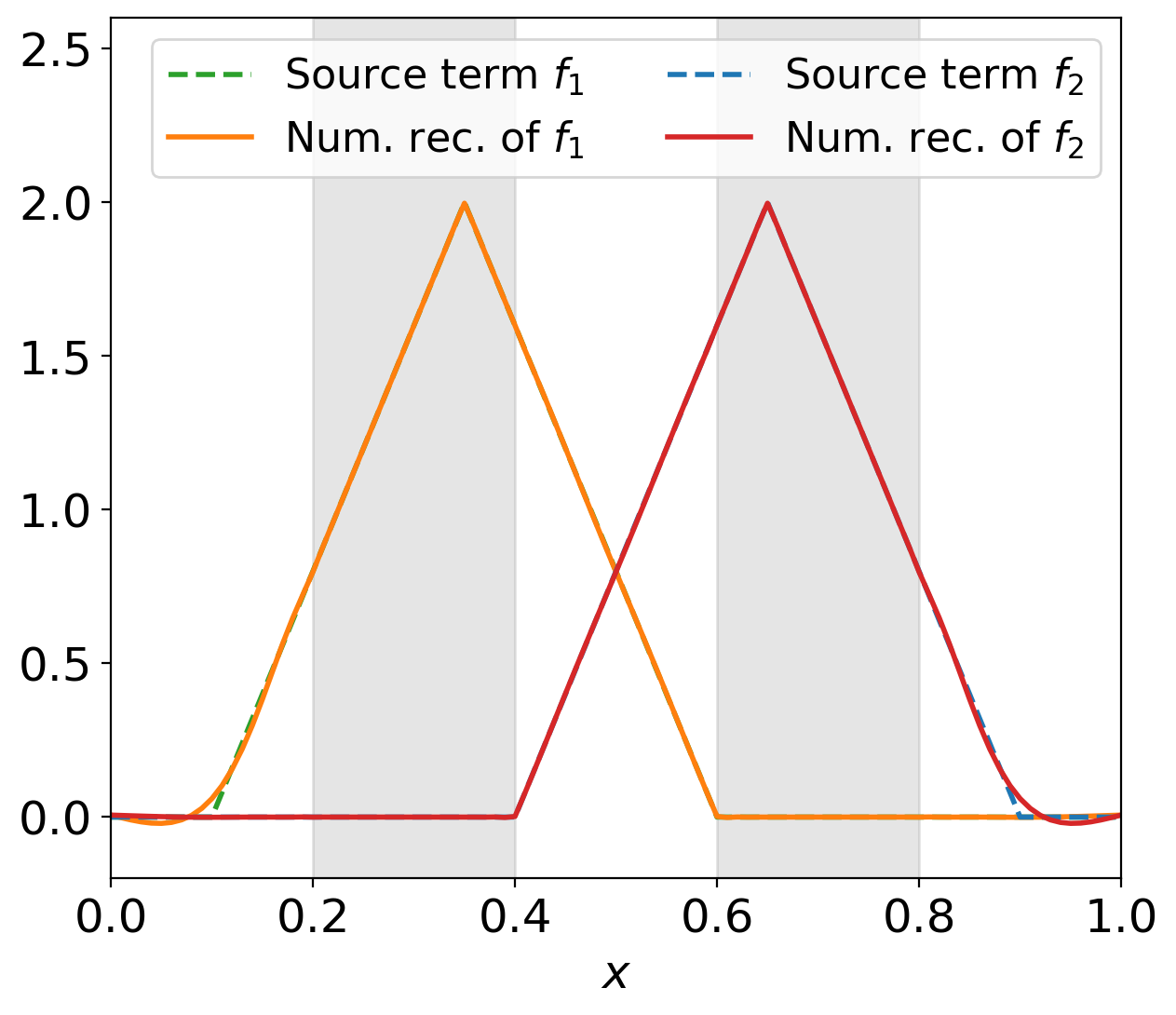}
    \caption{rel. error $1,0\%$}
    \end{subfigure}
    \begin{subfigure}{.3\textwidth}
    \includegraphics[width=1.0\linewidth]{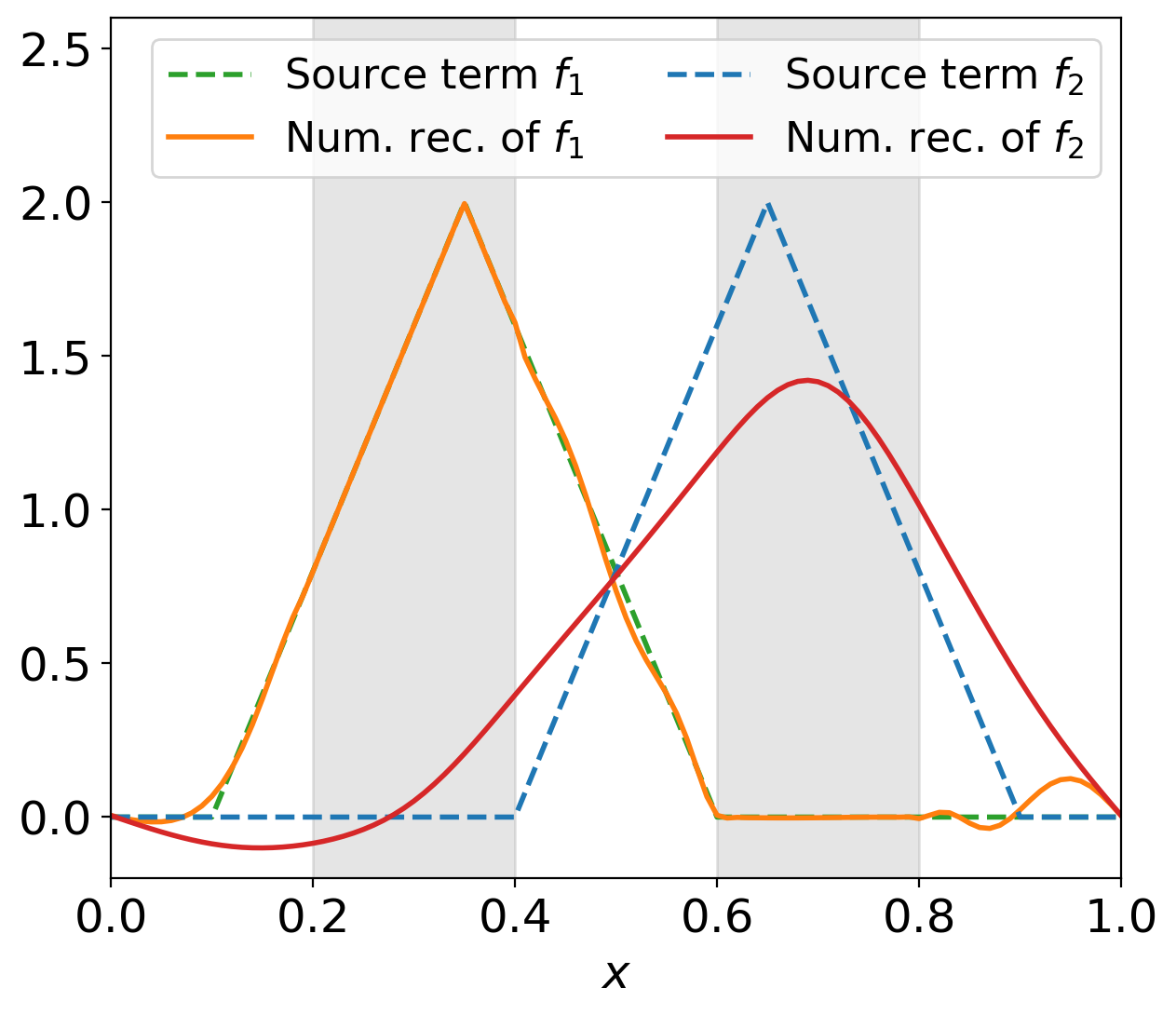}
    \caption{rel. error $21,2\%$}
    \end{subfigure}
    \begin{subfigure}{.3\textwidth}
    \includegraphics[width=1.0\linewidth]{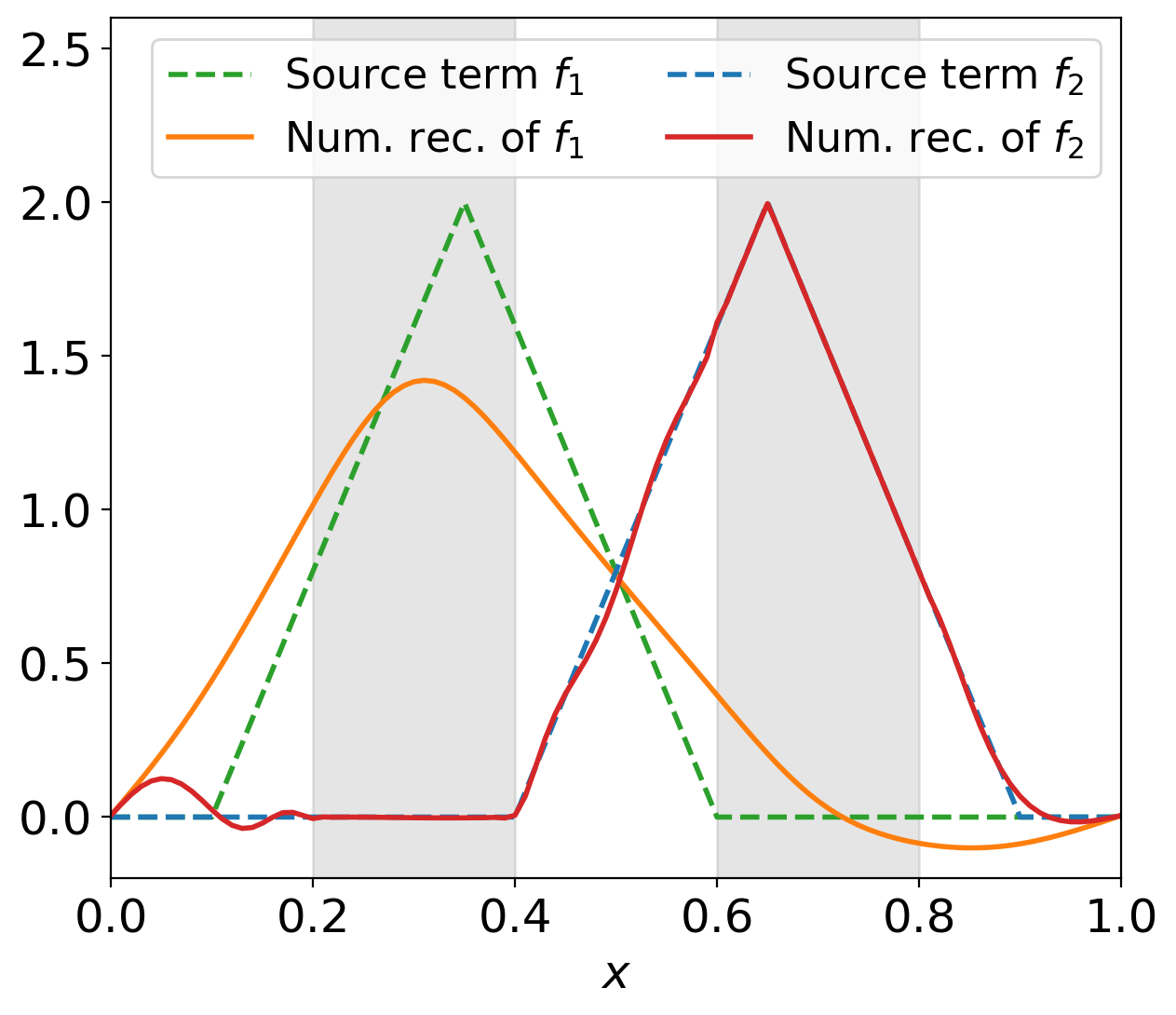}
    \caption{rel. error $21,2\%$}
    \end{subfigure}
    \caption{Source reconstructions when the coupling matrix is given by $q_{11}(x)=q_{22}(x)=0$, $q_{12}(x)=4x-2$, and  $q_{21}(x)=-4x+2$. . First column is the reconstruction of $F_1(x)$ and $F_2(x)$ observing both components, i.e., $(y_{1,\text{obs}},y_{2,\text{obs}})^t$. The second column shows the reconstruction of $F_1(x)$ and $F_2(x)$ observing the first component $y_{1,\text{obs}}$, and the last one is the reconstruction observing the second component $y_{2,\text{obs}}$ only.
    The highlighted gray areas correspond to the observation domain $\mathcal{O}$.}
    \label{fig:1D_reconstruction_gamma_linear}
\end{figure}

Figure~\ref{fig:1D_reconstruction_gamma_linear} presents reconstructions using a different coupling configuration $Q(x)$ given by $q_{11}(x)=q_{22}(x)=0$, $q_{12}(x)=4x-2$, and  $q_{21}(x)=-4x+2$. Note that this choice for the coupling matrix describes a setting more general than presented in Lemma~\ref{lemma.duprez}. %These reconstructions exhibit trends analogous to those observed in the cubic polynomial case, though with slightly improved accuracy overall. 
In this setting, reconstructions from measurements in both components consistently exhibit robust performance and relatively low reconstruction errors. On the other hand, reconstructions based on single-component measurements appear more sensitive to both the location and size of the observation subdomain, leading in some cases to larger reconstruction errors, particularly outside the observation region
Moreover, it can be observed from Figures~\ref{fig1D_reconstruction_gamma_cubic} and~\ref{fig:1D_reconstruction_gamma_linear} that reconstruction quality significantly improves when the observation subdomain covers regions with stronger variations of the source function. This indicates that source reconstruction accuracy is linked not only to the size and location of the observation domain but also to its alignment with the spatial structure of the source.

\subsection{Source reconstruction varying the regularization term}

We analyze the effect of the regularization term $k$ on the reconstruction quality in a one-dimensional setting. Specifically, we consider $F_1(x):=(\sin(2\pi x), - \sin(2\pi x))^t$, the observation domain $\mathcal{O}_1:= (0.5, 0.9)$, and the matrix coupling presented in Figure~\ref{fig:1D_reconstruction_gamma_linear}. All other parameters remain consistent with those from earlier experiments, except for the time discretization, which we set to $\Delta t=0.005$. Our goal here is to investigate how different choices of the regularization parameter $k$ influence the quality of the reconstructed source, particularly in regions inside and outside the observation subdomain.

\begin{table}
    \centering
    \begin{tabular}{c|cccccc}
    $k$ &  $10^{2}$ & $10^{3}$  & $10^{4}$  & $10^{5}$ & $10^{6}$ & $10^{7}$ \\
    \hline     
    % rel. error  &  $44,0\%$ & $21,3\%$ & $8,1\%$ & $10.1\%$ & $29,2\%$ \\
    rel. error  &  $44,0\%$ & $21,3\%$ & $8,1\%$ & $2.9\%$ & $1,3\%$ & $0,9\%$ \\
    \end{tabular}
    \caption{Error reconstruction for $F_1(x)$ under different values of the parameter $k$ for the functional $J$.}
    \label{table:1D_k_analysis}
\end{table}
Table~\ref{table:1D_k_analysis} displays reconstruction results using measurements from both state components, with a regularization parameter $k$ in $\{10^2, 10^3, 10^4, 10^5, 10^6, 10^7\}$. As illustrated, a large $k\sim 10^5$ produces significant improvements in the reconstruction. 
% Excessively large $k$ values may lead to overfitting the observations without fully capturing the PDE problem.  
Excessively large $k$ values yield smaller percentage improvements relative to the computational cost, which can increase by up to 10 times when making $k$ larger by a factor of 10.  
However, an excessively small regularization parameter (e.g., $k=10^2$) results in pronounced reconstruction errors, with oscillations and deviations becoming particularly evident far from the observation domain. Thus, selecting an adequately large regularization parameter is essential for balancing accuracy and stability in source reconstruction tasks.

\subsection{Source reconstruction in 2D}
In this subsection, we analyze source reconstructions on the unit square domain $\Omega=(0,1)^2$, using parameters $T=0.5$ and $\nu = 0.1$. For the numerical discretization, we utilize a regular triangular mesh composed of $3200$ elements, with a temporal discretization step size $\Delta t=0.005$, and $k=10^5$. %We test two sets of observation domains: $\mathcal{O}_1=(0.5,0.9)\times(0.1,0.9)$ and $\mathcal{O}_2=((0.2, 0.4)\cup(0.6,0.8))\times(0.3,0.7)$. 
The reconstruction accuracy is evaluated using the source function 
$$F_3(x,y):=(\sin(2\pi x)\sin(2\pi y),-\sin(2\pi x)\sin(2\pi y))^t.$$

\begin{figure}[htp]
    \centering
    \begin{subfigure}{.24\textwidth}
    \includegraphics[width=1.0\linewidth]{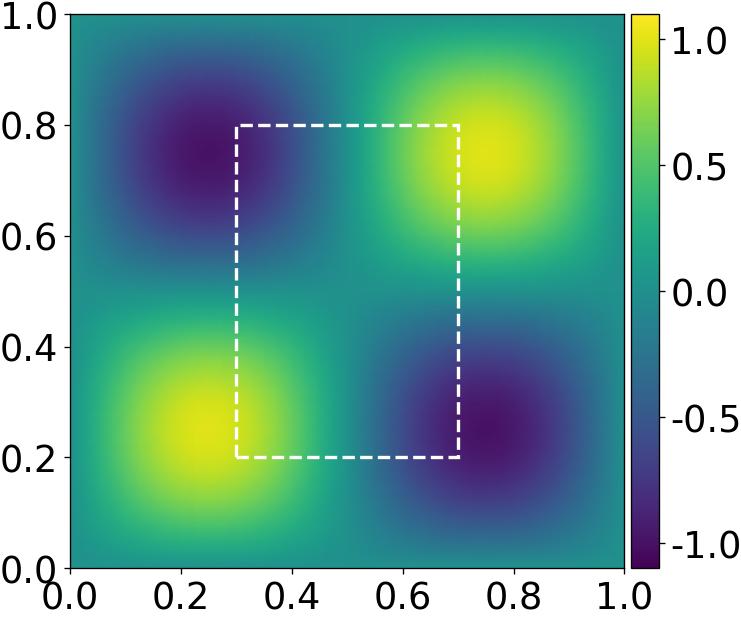}
    \caption{first component}
    \end{subfigure}
    \begin{subfigure}{.24\textwidth}
    \includegraphics[width=1.0\linewidth]{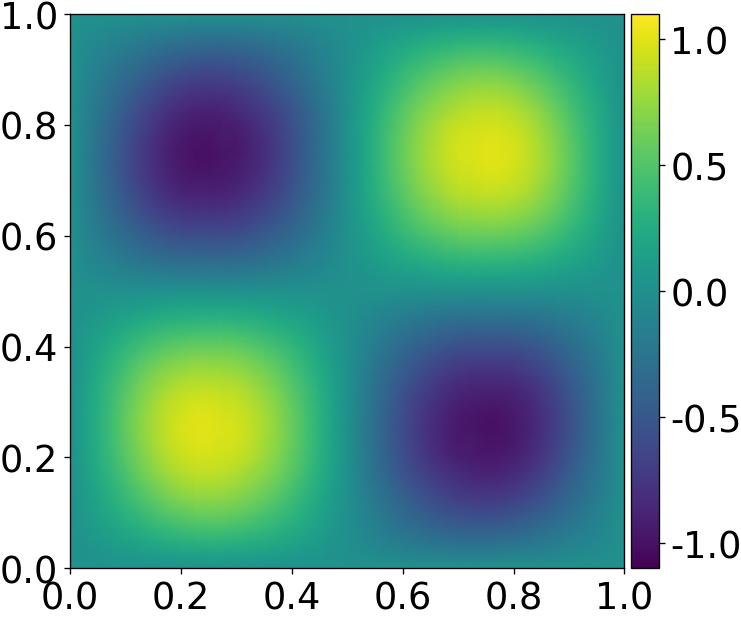}
    \caption{rel. error $2,5\%$}
    \end{subfigure}
    \begin{subfigure}{.24\textwidth}
    \includegraphics[width=1.0\linewidth]{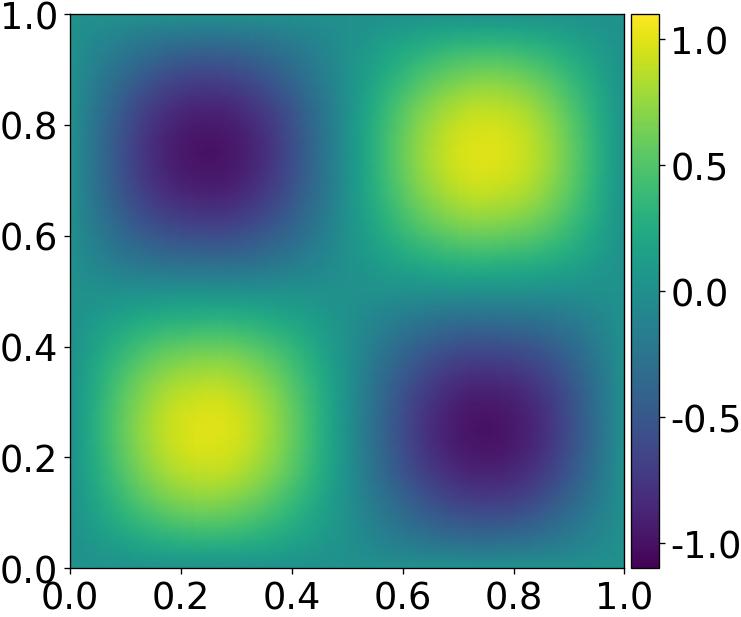}
    \caption{rel. error $1,1\%$}
    \end{subfigure}
    % \begin{subfigure}{.24\textwidth}
    % \includegraphics[width=1.0\linewidth]{numerical-result/2D/source_rec_1_comp_ex01_2D_f1_sigma3_obs_1_2nd_comp.jpg}
    % \caption{rel. error $100,0\%$}
    % \end{subfigure}
    \\
    \begin{subfigure}{.24\textwidth}
    \includegraphics[width=1.0\linewidth]{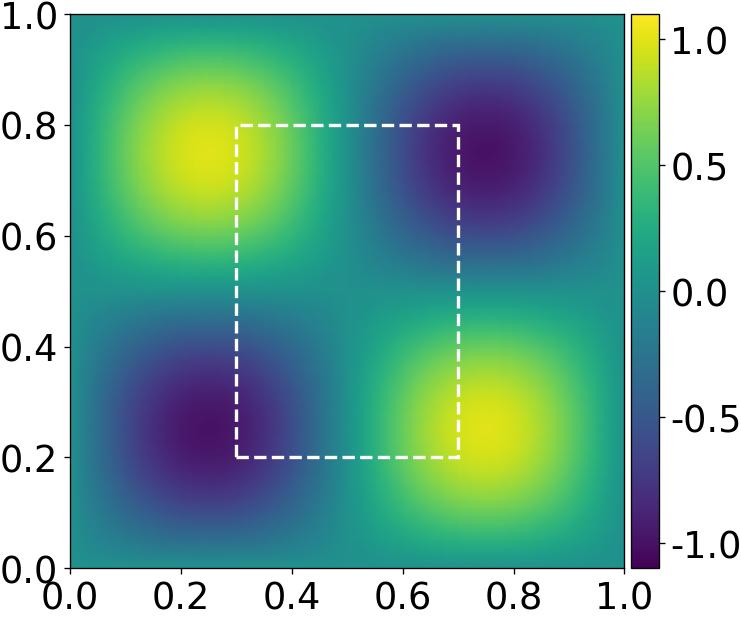}
    \caption{second component}
    \end{subfigure}
    \begin{subfigure}{.24\textwidth}
    \includegraphics[width=1.0\linewidth]{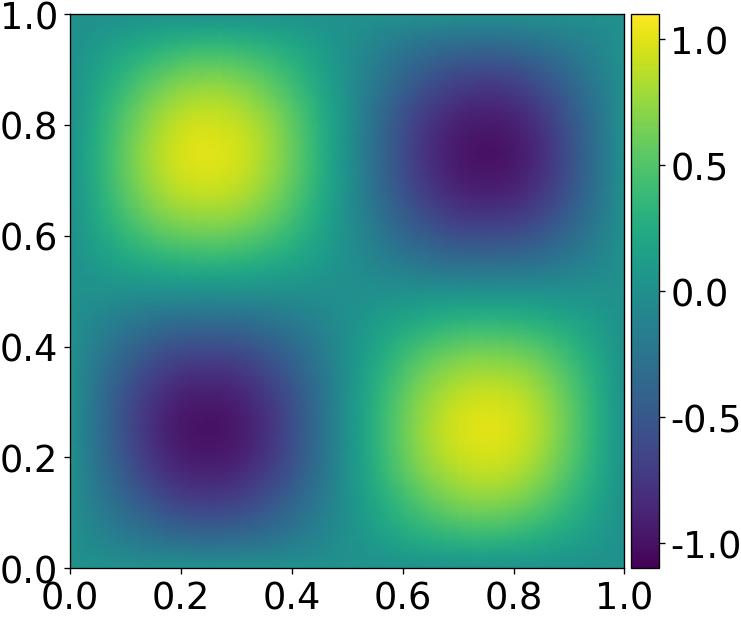}
    \caption{rel. error $1,6\%$}
    \end{subfigure}
    \begin{subfigure}{.24\textwidth}
    \includegraphics[width=1.0\linewidth]{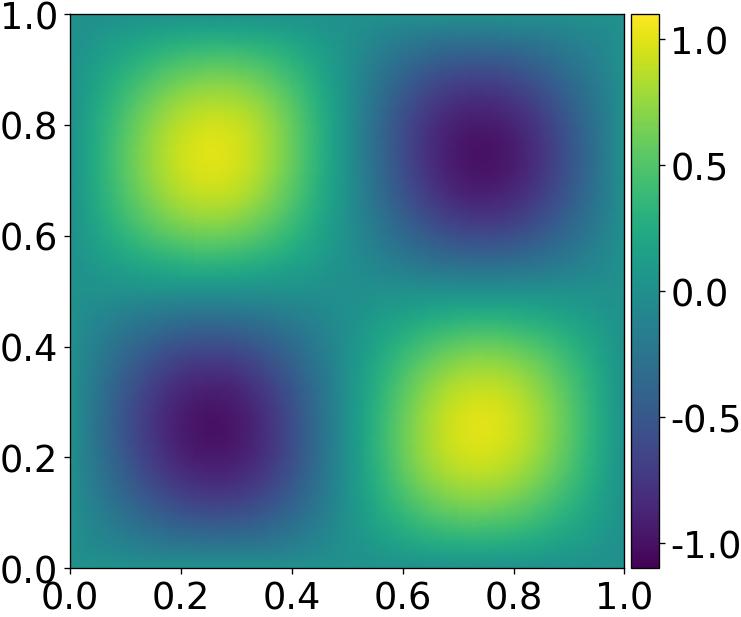}
    \caption{rel. error $4,0\%$}
    \end{subfigure}
    % \begin{subfigure}{.24\textwidth}
    % \includegraphics[width=1.0\linewidth]{numerical-result/2D/source_rec_2_comp_ex01_2D_f1_sigma3_obs_1_2nd_comp.jpg}
    % \caption{rel. error $9,7\%$}
    % \end{subfigure}
    \\
    \begin{subfigure}{.24\textwidth}
    \includegraphics[width=1.0\linewidth]{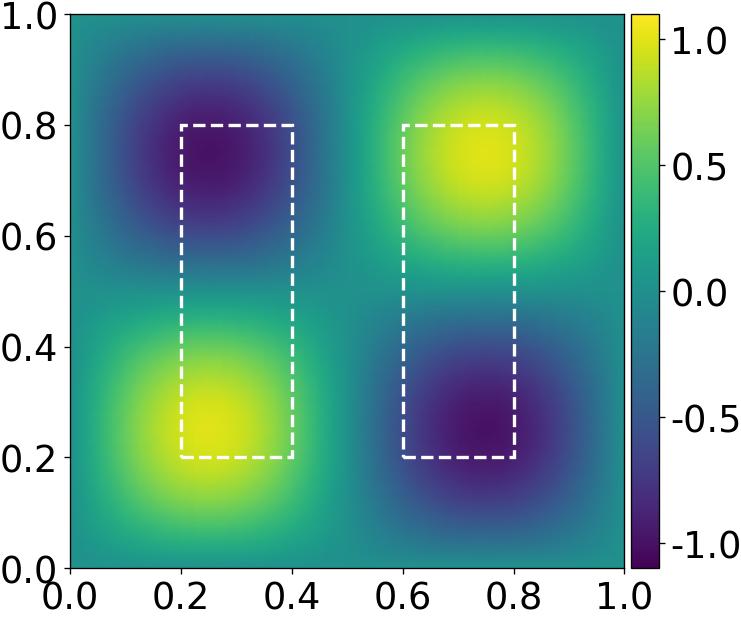}
    \caption{first component}
    \end{subfigure}
    \begin{subfigure}{.24\textwidth}
    \includegraphics[width=1.0\linewidth]{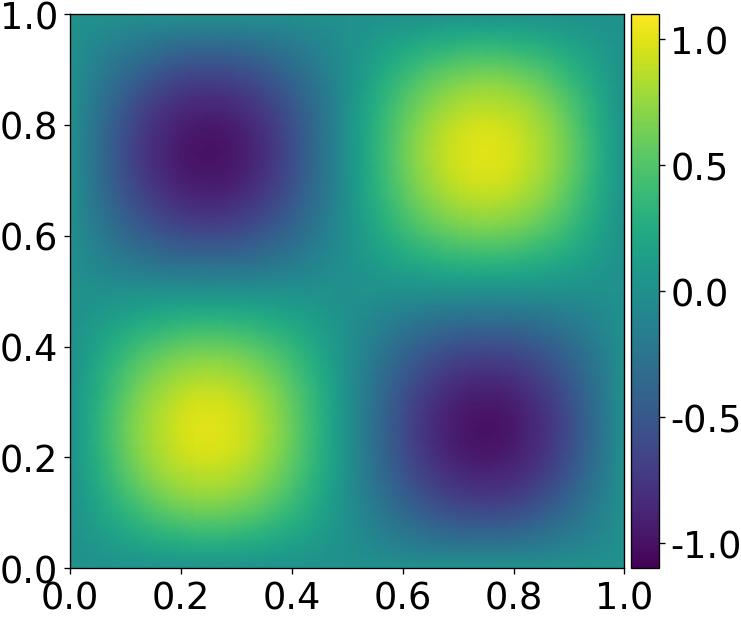}
    \caption{rel. error $1,0\%$}
    \end{subfigure}
    \begin{subfigure}{.24\textwidth}
    \includegraphics[width=1.0\linewidth]{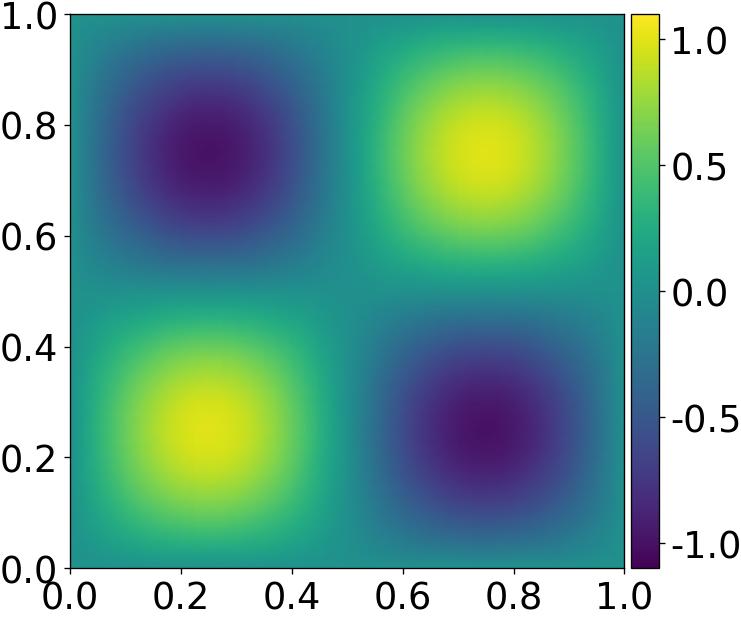}
    \caption{rel. error $0.9\%$}
    \end{subfigure}
    % \begin{subfigure}{.24\textwidth}
    % \includegraphics[width=1.0\linewidth]{numerical-result/2D/source_rec_1_comp_ex01_2D_f1_sigma3_obs_2_2nd_comp.jpg}
    % \caption{rel. error $100,0\%$}
    % \end{subfigure}
    \\
    \begin{subfigure}{.24\textwidth}
    \includegraphics[width=1.0\linewidth]{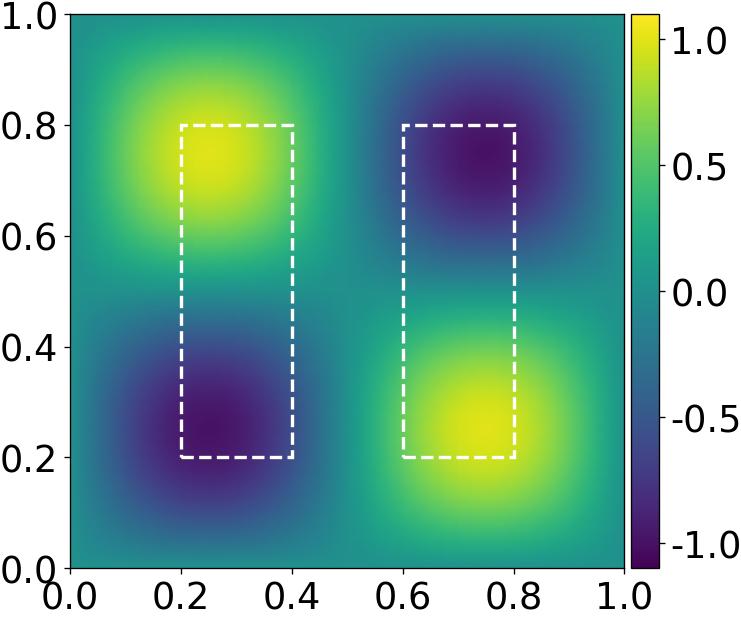}
    \caption{second component}
    \end{subfigure}
    \begin{subfigure}{.24\textwidth}
    \includegraphics[width=1.0\linewidth]{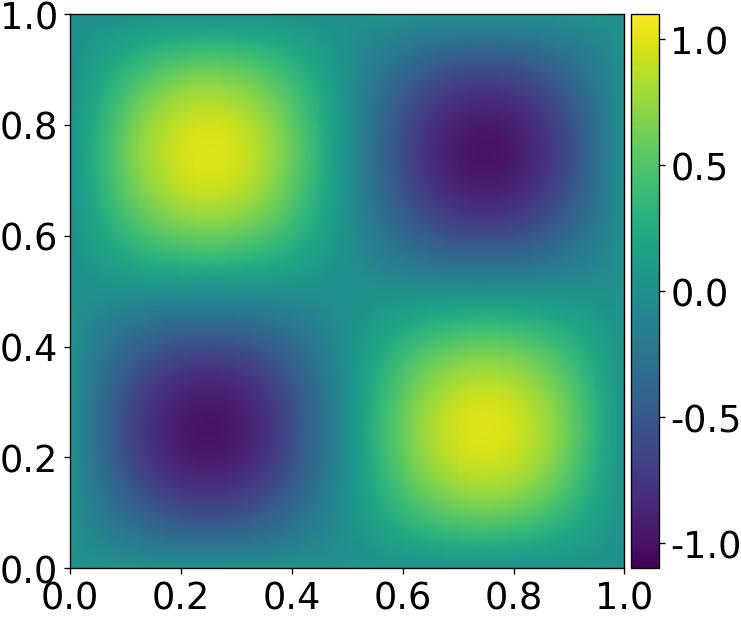}
    \caption{rel. error $0,8\%$}
    \end{subfigure}
    \begin{subfigure}{.24\textwidth}
    \includegraphics[width=1.0\linewidth]{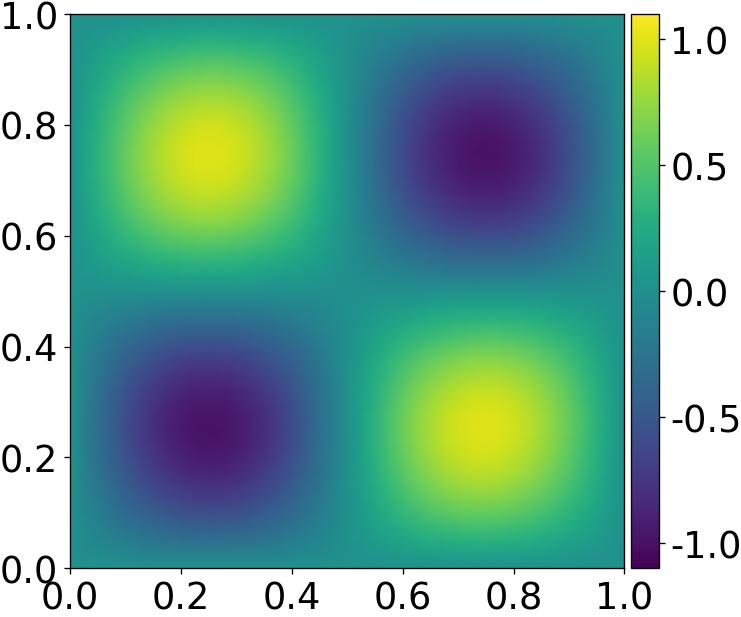}
    \caption{rel. error $1,7\%$}
    \end{subfigure}
    % \begin{subfigure}{.24\textwidth}
    % \includegraphics[width=1.0\linewidth]{numerical-result/2D/source_rec_2_comp_ex01_2D_f1_sigma3_obs_2_2nd_comp.jpg}
    % \caption{rel. error $3,4\%$}
    % \end{subfigure}
    \caption{Source reconstruction for $F_3(x,y)$ using a coupling matrix given by  $q_{11} = 1$, $q_{12} = 4$, $q_{21} = 0$, and $q_{22} = 1$. %The first two rows show reconstructions with measurements on  $\mathcal{O}_1$ (top), and the last two rows show reconstructions with measurements on $\mathcal{O}_2$ (bottom). 
    The first column displays the sources and the locations of the observation domains, indicated by a white dashed line. The second column displays the reconstructions from observing both components. The third column shows the reconstruction obtained by observing only the first component.}
    \label{fig2D_module_comparisom_ex1d11}
\end{figure}

\begin{figure}[htp]
    \centering
    \begin{subfigure}{.24\textwidth}
    \includegraphics[width=1.0\linewidth]{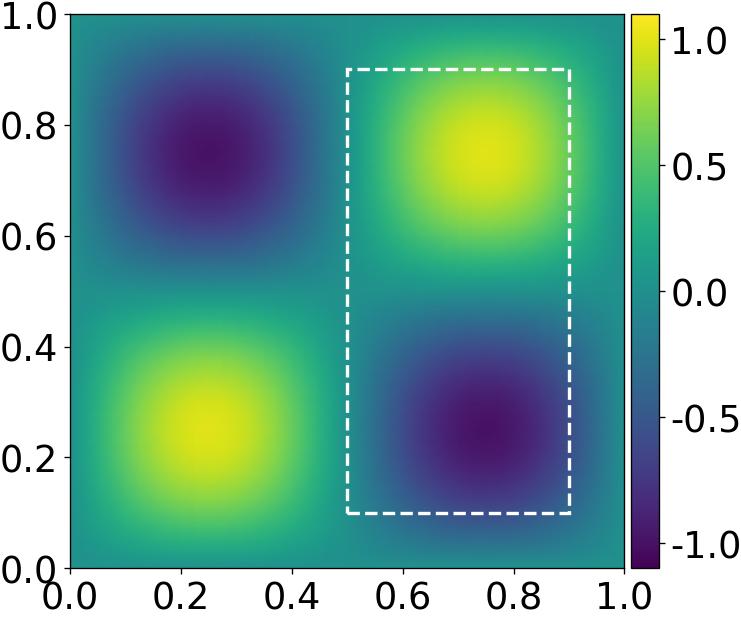}
    \caption{1st component}
    \end{subfigure}
    \begin{subfigure}{.24\textwidth}
    \includegraphics[width=1.0\linewidth]{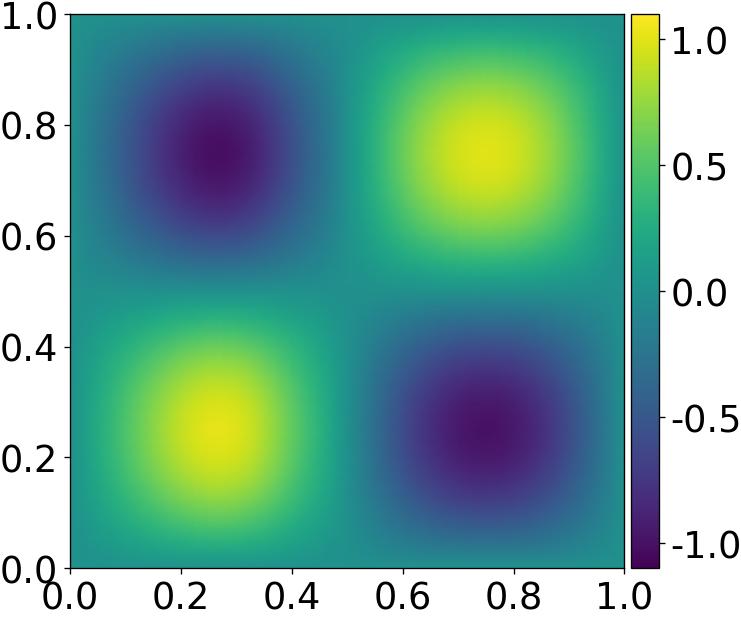}
    \caption{rel. error $5,5\%$}
    \end{subfigure}
    \begin{subfigure}{.24\textwidth}
    \includegraphics[width=1.0\linewidth]{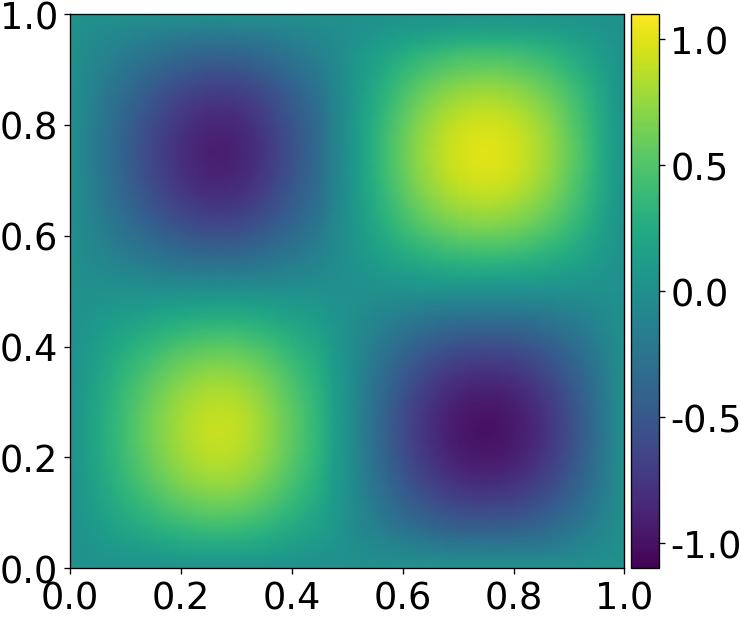}
    \caption{rel. error $9,9\%$}
    \end{subfigure}
    \begin{subfigure}{.24\textwidth}
    \includegraphics[width=1.0\linewidth]{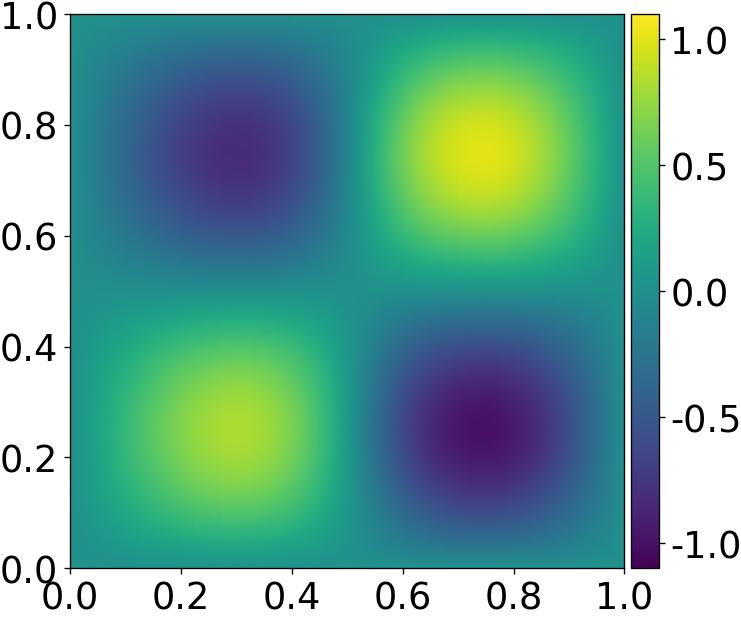}
    \caption{rel. error $17,5\%$}
    \end{subfigure}
    \begin{subfigure}{.24\textwidth}
    \includegraphics[width=1.0\linewidth]{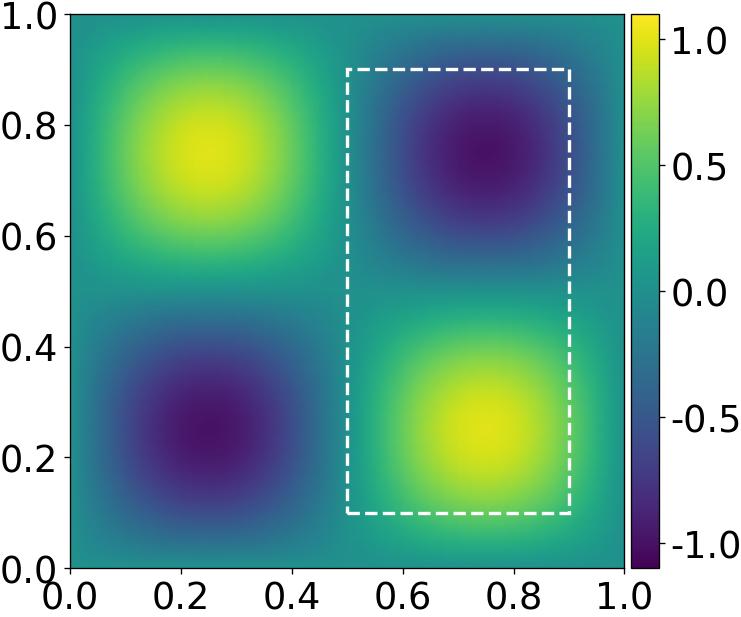}
    \caption{2nd component}
    \end{subfigure}
    \begin{subfigure}{.24\textwidth}
    \includegraphics[width=1.0\linewidth]{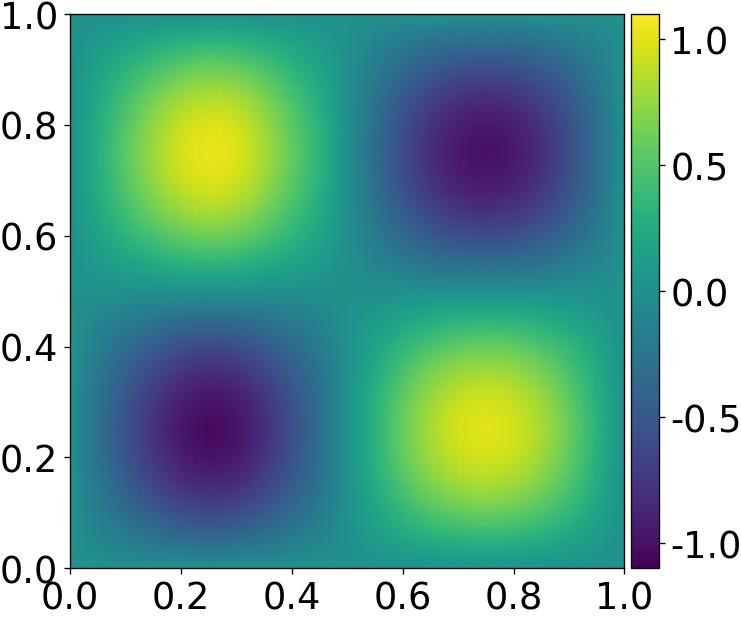}
    \caption{rel. error $2,8\%$}
    \end{subfigure}
    \begin{subfigure}{.24\textwidth}
    \includegraphics[width=1.0\linewidth]{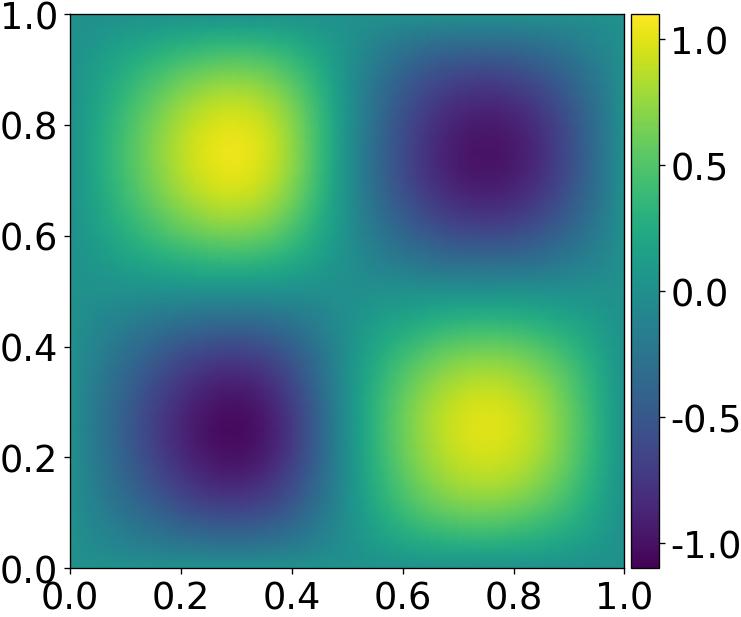}
    \caption{rel. error $10,6\%$}
    \end{subfigure}
    \begin{subfigure}{.24\textwidth}
    \includegraphics[width=1.0\linewidth]{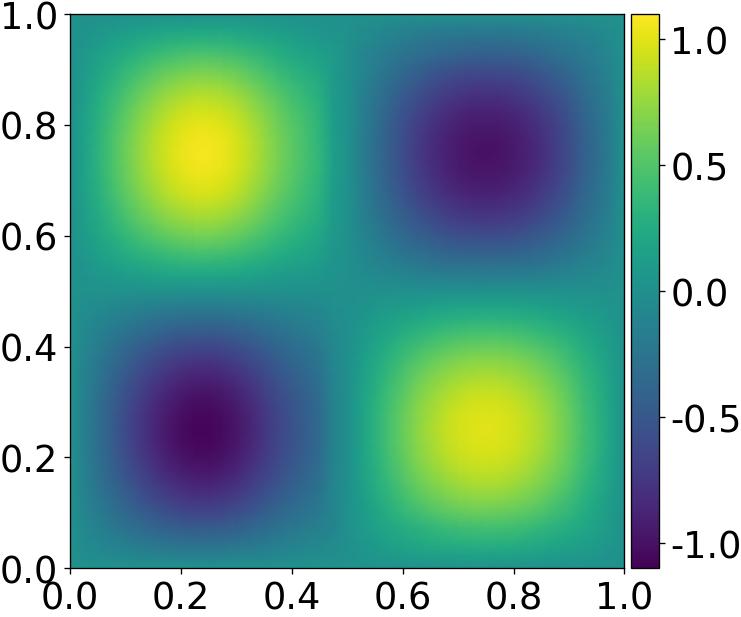}
    \caption{rel. error $4,3\%$}
    \end{subfigure}
    \begin{subfigure}{.24\textwidth}
    \includegraphics[width=1.0\linewidth]{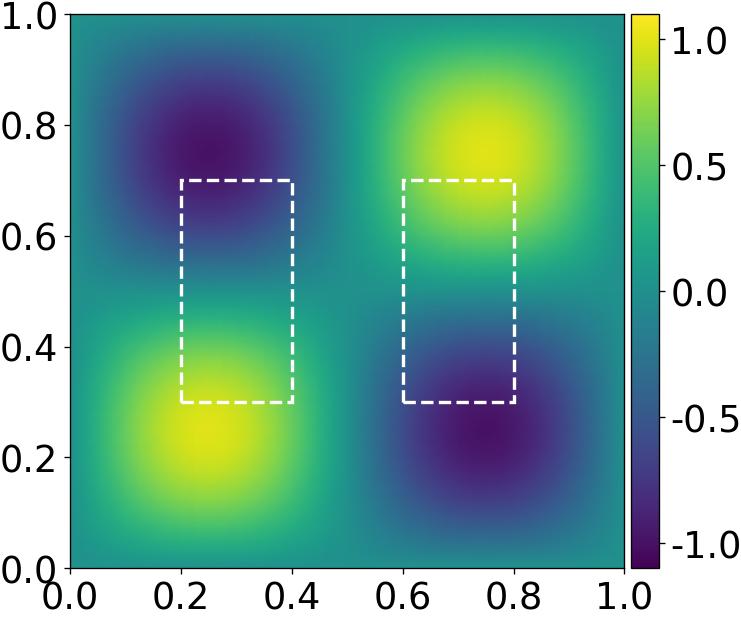}
    \caption{1st component}
    \end{subfigure}
    \begin{subfigure}{.24\textwidth}
    \includegraphics[width=1.0\linewidth]{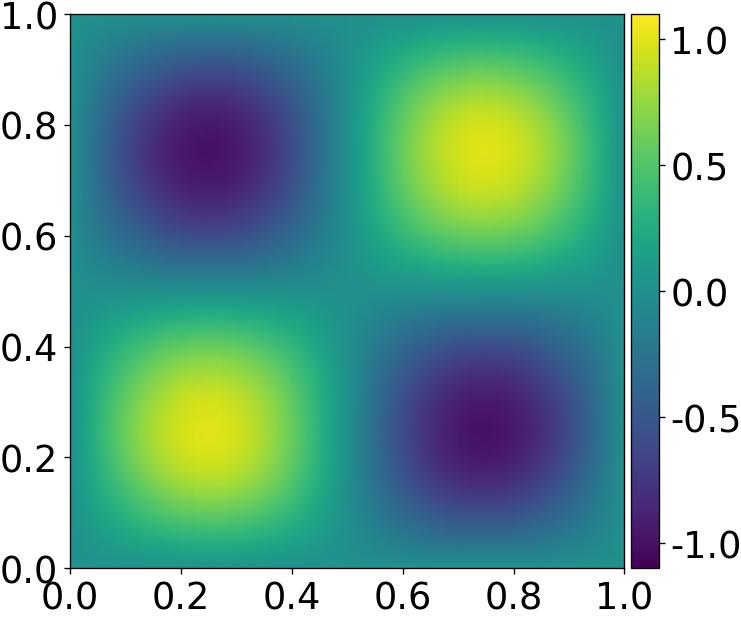}
    \caption{rel. error $1,5\%$}
    \end{subfigure}
    \begin{subfigure}{.24\textwidth}
    \includegraphics[width=1.0\linewidth]{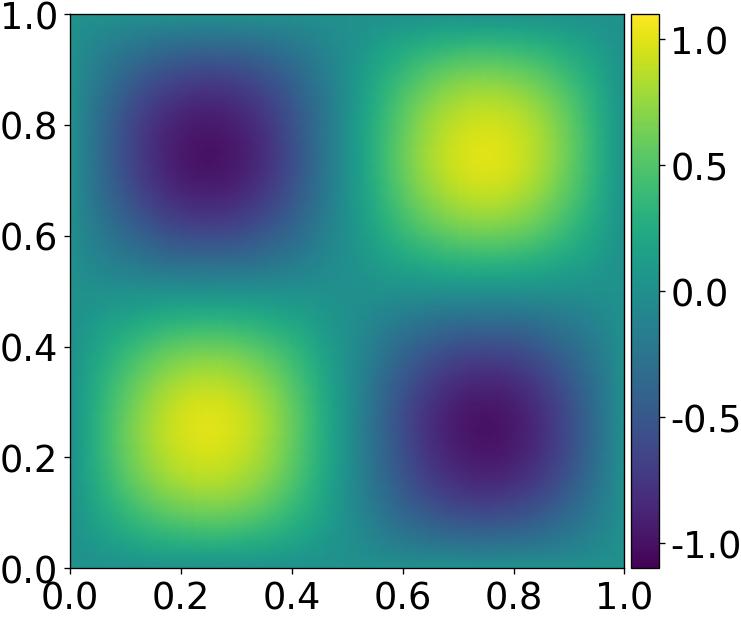}
    \caption{rel. error $1,3\%$}
    \end{subfigure}
    \begin{subfigure}{.24\textwidth}
    \includegraphics[width=1.0\linewidth]{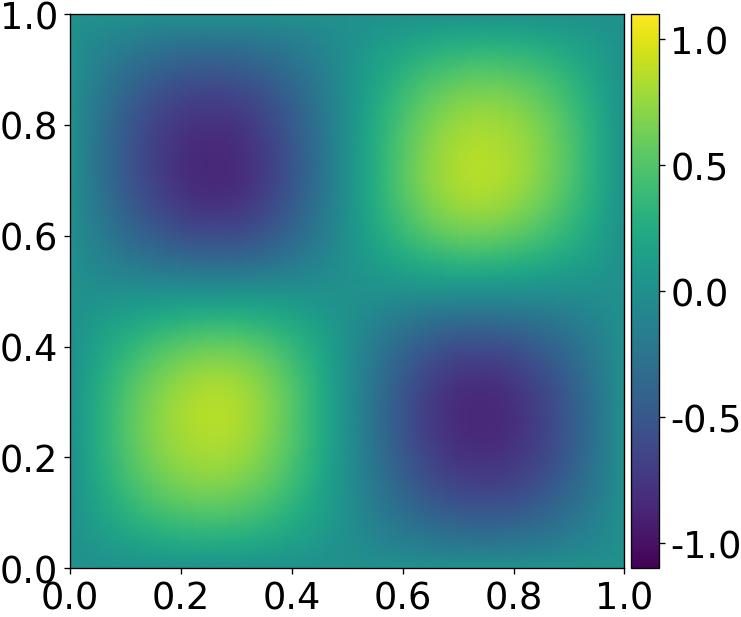}
    \caption{rel. error $13,4\%$}
    \end{subfigure}
    \begin{subfigure}{.24\textwidth}
    \includegraphics[width=1.0\linewidth]{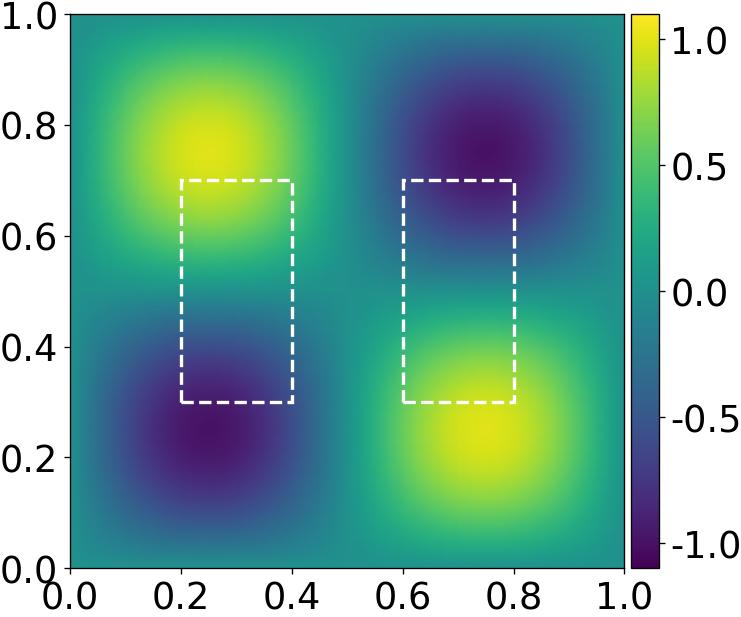}
    \caption{2nd component}
    \end{subfigure}
    \begin{subfigure}{.24\textwidth}
    \includegraphics[width=1.0\linewidth]{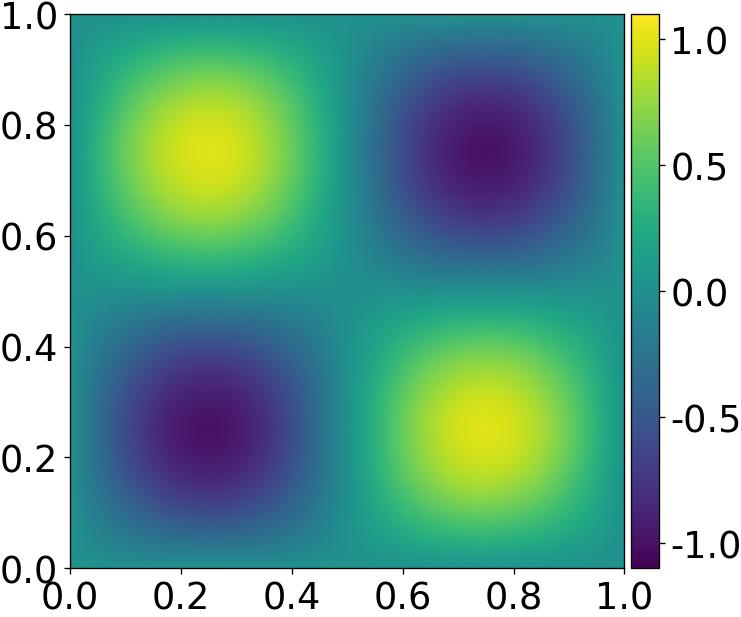}
    \caption{rel. error $0,9\%$}
    \end{subfigure}
    \begin{subfigure}{.24\textwidth}
    \includegraphics[width=1.0\linewidth]{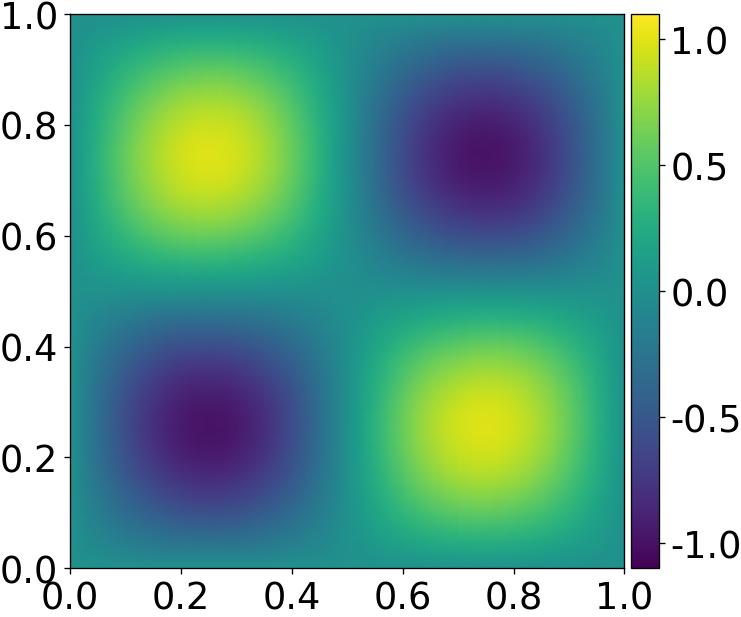}
    \caption{rel. error $2,7\%$}
    \end{subfigure}
    \begin{subfigure}{.24\textwidth}
    \includegraphics[width=1.0\linewidth]{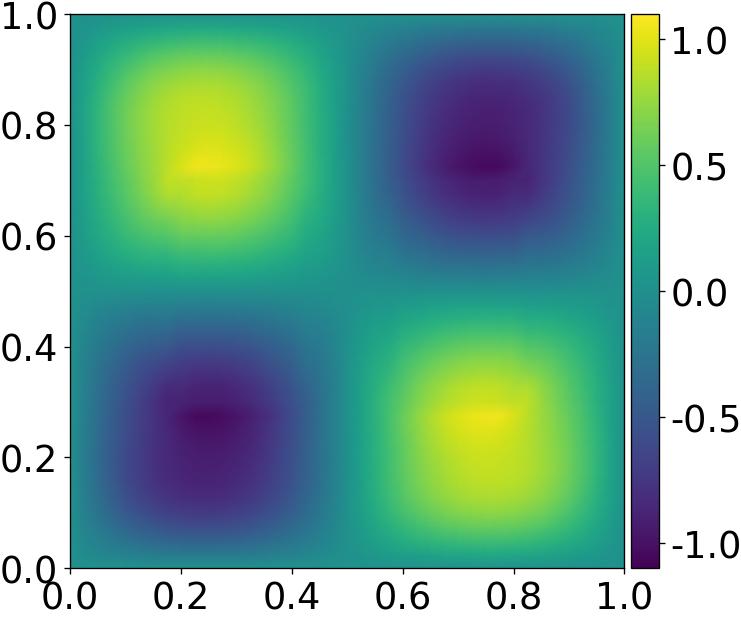}
    \caption{rel. error $10,0\%$}
    \end{subfigure}
    \caption{Source reconstruction for $F_3(x,y)$ using the coupling matrix  $q_{11}=q_{22}=0$, $q_{12} = 4$ and $q_{21}= 2$. %The first two rows show reconstructions with measurements on  $\mathcal{O}_1$ (top), and the last two rows show reconstructions with measurements on $\mathcal{O}_2$ (bottom). 
    The first column displays the sources and the locations of the observation domains, indicated by a white dashed line. The second, third, and fourth columns display the reconstructions with both components, the first component, and the second component, respectively.}
   \label{fig2D_module_comparisom_ex2d1}
\end{figure}

Figures~\ref{fig2D_module_comparisom_ex1d11} and~\ref{fig2D_module_comparisom_ex2d1} illustrate the numerical reconstruction results obtained for the above source. Figure~\ref{fig2D_module_comparisom_ex1d11} shows reconstructions from a constant coupling term given by $q_{11} = 1$, $q_{12} = 4$, $q_{21} = 0$, and $q_{22} = 1$.
%$$Q=\begin{pmatrix} 1 & 4 \\ 0 & 1\end{pmatrix}.$$
The two top rows present results based on observations within the domain $\mathcal{O}_3 = (0.3, 0.7)\times(0.2,0.8)$, while the two bottom rows show reconstructions derived from measurements in the domain $\mathcal{O}_4 = ((0.2, 0.4)\cup(0.6,0.8))\times(0.2,0.8)$. 
Figure~\ref{fig2D_module_comparisom_ex2d1} shows reconstructions from a more general constant coupling term
using the coupling matrix  $q_{11}=q_{22}=0$, $q_{12}(x,y) = 4$ and $q_{21}(x,y) = 2$,
%$$Q=\begin{pmatrix} 0 & 2 \\ 4 & 0\end{pmatrix},$$ 
with $\mathcal{O}_5 = (0.5, 0.9)\times(0.1,0.9)$ in the two top rows and $\mathcal{O}_6 = ((0.2, 0.4)\cup(0.6,0.8))\times(0.3,0.7)$ in the two bottom rows. Note that the last choice for $Q$ describes a setting more general than presented in Lemma~\ref{lema.controlQ}.

The relative error of each reconstruction is computed by component. The error goes down in dependency of the observation domain $\mathcal{O}\subset\Omega$. To be precise,  when the measure of the set $\Omega\backslash\mathcal{O}$ tends to zero, i.e., $\mu(\Omega\backslash\mathcal{O})\to 0$, the relative error decreases. However, for the cases where $\mu(\Omega\backslash\mathcal{O})\gg 0$, weaker reconstruction results are obtained. 
This leads to the interesting question of investigating geometric conditions to determine the optimal size of the observation domain $\mathcal{O}$ in order to obtain stronger results.

\subsection{Implementation notes}\label{sec:impl_details}
Through this numerical section, we employ a minimization-based formulation that circumvents the need to repeatedly solve large-scale controllability problems and Volterra systems, while still preserving the essential structural concepts of the theoretical framework. 
This methodological shift is due to the fact that implementing the reconstruction formulas on realistic meshes with practical time resolutions becomes computationally prohibitive. In what follows, we provide a detailed analysis of the computational cost associated with the controllability problems and the Volterra systems:
\begin{enumerate}
    \item In the null controllability problems, for each $\tau\in(0,T)$ and each eigenfunction indexed by $k$, we must solve a controllability system with initial data $\Psi_k^{0}$ (constructed from Laplace eigenfunctions). Implementing this via the Hilbert Uniqueness Method (HUM) requires both time discretization and space discretization. For a system with $n$ equations over a mesh with $M$ nodes, the cost is roughly $O(M\times n\times L)$ per $k$, where $L$  is the number of iterations required for convergence. Note that this grows quickly with mesh refinement or larger $n$.
    \item For the Volterra equations, once the control functions $U^{(\tau)}$ are computed, $n$ copies of the Volterra equations (first or second kind, depending on $\sigma(0)$) must be solved for each component and for each $k$. Even with efficient quadrature (e.g., trapezoidal rule), the cost is approximately $O(M\times N \times n\times (\ell + 1))$ per $k$, where $N$ is the number of time steps and $\ell$ indexes the time discretization in the integral operator.
\end{enumerate}
The functional defined in \eqref{minproblem.functional} provides a framework for analyzing the inverse source problem via an optimization-based approach. Its construction is consistent with the Lipschitz-type stability property established in Appendix~\ref{appendix.lipschitz}, thereby ensuring the theoretical soundness of the formulation. Thus, the selection of this functional is motivated by the inequality stated in Appendix~\ref{appendix.lipschitz}, which is associated with the Lipschitz stability result for the inverse source problem presented in Problem~\ref{Problem1}.

%% file: 03_conclusions.tex
\section{Conclusions}\label{conclusions}

In this paper, we introduced and analyzed inverse problems aimed at recovering spatial source distributions in coupled parabolic systems. We addressed two types of coupling: one involving a potential matrix with constant components and another involving a potential matrix with spatially dependent entries
%, and the other involving coupling through the main operator. 
A significant contribution of our work is the development of source reconstruction algorithms based on local internal measurements from one scalar component of the system's state. To the best of our knowledge, there are no previous studies on reconstruction algorithms for these sources.

One of the primary advantages of our approach is the minimal requirement for the measured state variables needed to reconstruct the source coefficients. Nevertheless, several related theoretical questions remain unresolved, particularly those concerning stability and uniqueness of the inverse source problem using the methodology via spectral analysis, control null results, and Volterra equations. This lack of resolution is not unexpected; even in the scalar case analyzed in~\cite{2013GOT}, achieving stability required stringent assumptions on the temporal derivative $\sigma'(t)$. Specifically, a stability result was obtained when $\sigma'(t)$  decreases exponentially in $(0,T)$, i.e., $\sigma'(t)\approx e^{-1/t}$ (see~\cite[Proof of Theorem 1.3, Step 4]{2013GOT}). 
This particular choice of $\sigma'(t)$  is closely linked to the null controllability problem since the control function cost has an exponential dependence of the form $e^{1/T}$ (see~\cite[Lemma 1.1]{2013GOT}). Consequently, establishing stability via our methodology would similarly require precise control-cost estimates for coupled parabolic systems. Moreover, such estimates must simultaneously allow appropriate bounds on $\sigma(t)$ and $\sigma'(t)$ while ensuring $\sigma(T)\neq 0$. Currently, however, obtaining these estimates remains an open problem, as the control function cost involves an exponential term of the form $e^{T+1/T}$. To provide a positive answer regarding the stability property, we have included a Lipschitz-type stability result for Problem \ref{Problem1} via Carleman inequalities in the Appendix \ref{appendix.lipschitz}.

From a theoretical perspective, our study provides reconstruction algorithms from local observations of only one scalar state; this problem has not been previously addressed in the existing literature. Our methodology intricately combines spectral analysis, null controllability properties, and Volterra equations (see Figure~\ref{fig:diagram}). 
In addition, we presented a spectral analysis for $2\times2$ systems with space-dependent coefficients (see Problem~\ref{Problem2}), which required employing Riesz bases due to the skew-adjoint nature of the main operator. 
This scenario contrasts with previous works, where self-adjoint operators allowed simpler Fourier-based representations.  
It is worth mentioning that the spectral study for a system of $2 \times 2$ equations with a skew-adjoint matrix operator is more challenging and involves concepts that are outside the scope of this work. 

From a numerical point of view, given the high computational cost of implementing reconstruction formulas via controllability and Volterra systems, we proposed a minimization-based approach. This method avoids repeatedly solving large-scale control problems and integral equations, while still preserving the core ideas of our theoretical framework. As such, it offers a more practical and scalable alternative for different space and time resolutions. Then, the inverse source problems were tackled using a descent method for an appropriate functional cost, including measures from the state and its first derivative in local subdomains (in line with the theoretical results). The parabolic systems were solved in 1D and 2D using finite element methods. We also presented a numerical analysis of the regularization term. We also compare the source reconstruction between local observations from one scalar component and all components of the state in every numerical experiment.
%presented in this work.

%% file: 04_appendix.tex
%%%%%%%%%%%%%%%%%%%%%%%%%%%%%%%%%%%%%%%%%%%%%%%%%%
% APPENDIX %%%%%%%%%%%%%%%%%%%%%%%%%%%%%%%%%%%%%%%
\appendix

\section{Lipschitz stability via Carleman estimates}\label{appendix.lipschitz}
	In this appendix, we prove a Lipschitz-type stability associated with Problem \ref{Problem1}. The strategy is based on the Bukhgeim-Klibanov method \cite{BK81}. Thus, we will use a global Carleman inequality satisfied by the solutions of 
    \begin{equation}\label{sys.linearforCarleman}
   		\left\{
    	\begin{array}{llll}
        	\partial_t Y-\Delta Y+QY=H, &\mbox{in} &\Omega\times(0,T),\\
         	Y=0, &\mbox{in} &\partial\Omega\times(0,T),\\
			Y(\,\cdot\, ,0)=0, &\mbox{in}&\Omega,        
    	\end{array}\right.
	\end{equation}	
	where $H=(h_1,\dots, h_n)^t\in L^2(0,T;L^2(\Omega)^n)$ and $Q\in L^\infty(\Omega)^{n^2}$ is the coupling matrix defined as in \eqref{def.matrixQandB} and \eqref{cond.matrixQ}.

	\begin{lemma}\label{lema.carleman}
		Let $\mathcal{O}\subset\Omega \subset \mathbb{R}^d$ be a nonempty open subset and $r\leq 3$. Then, there exists 
		a function $\alpha_0\in C^2(\overline{\Omega})$ (depending only on $\Omega$ and $\mathcal{O}$) and 
		three positive constants $s_1=s_1(\Omega,\mathcal{O}, T, \|q_{ij}\|_{\infty})$, 
		$C=C(\Omega,\mathcal{O}, T, \|q_{ij}\|_{\infty})$, and $\ell>3$ such that, for any $s\geq s_1$, the solution $Y$ of the problem~\eqref{sys.linearforCarleman} 
		satisfies the inequality
		\begin{equation}\label{ine.carleman}\begin{array}{lll}
			\sum\limits_{i=1}^{n}\mathcal{I}(3(n+1-i),y_i)& \leq C\Biggl(s^{\ell}\displaystyle\iint\limits_{\mathcal{O}\times(0,T)}
			e^{-2s\alpha}\xi(t)^{\ell}|y_n|^2\,dx\,dt\\
			&\hspace{1cm}+s^{r-3}\sum\limits_{i=1}^n\,\displaystyle\iint\limits_{\Omega\times(0,T)}e^{-2s\alpha}\xi(t)^{r-3}|h_i|^2\,dx\,dt\Biggr),
		\end{array}
		\end{equation}
		where 
		\begin{equation}\label{eq.termglobal.carleman}
			\mathcal{I}(r,z)=\iint\limits_{\Omega\times(0,T)}e^{-2s\alpha}\Biggl( (s\xi)^{r-4}|\partial_t z|^2
			+(s\xi)^{r-2}|\nabla z|^2+(s\xi)^{r}|z|^2\Biggr)\,dx\,dt,
		\end{equation}
		and the functions $\alpha$ and $\xi$ are defined as 
		\begin{equation*}
			\alpha(x,t)=\frac{\alpha_0(x)}{t(T-t)},\, (x,t)\in\Omega\times(0,T),\quad \xi(t)=(t(T-t))^{-1}, t\in (0,T).
		\end{equation*}
	\end{lemma}
	\begin{remark}
        The proof of Lemma \ref{lema.carleman} can be found in~\cite[Theorem 1.1]{2010Burgos}. We note that the first term of $\mathcal{I}(r,z)$, as well as the last term in the right-hand side of~\eqref{ine.carleman}, do not appear explicitly in \cite[Theorem 1.1]{2010Burgos}. However, by following the arguments therein, these additional terms can be incorporated without difficulty. Moreover, although Lemma \ref{lema.carleman} is valid for all 
		$r\in\mathbb{R}$, we restrict ourselves to the case $r\leq 3$ in order to prove the Lipschitz stability result associated with our inverse source problem.
	\end{remark}
    \begin{theorem}[Lipschitz stability]\label{teo.stability}
        Let $T'=\frac{T}{2}$. Suppose that $\sigma\in W^{1,\infty}(0,T)$ and that there exists a positive constant $\gamma_{0}$ such that $|\sigma(T')|\geq \gamma_{0}>0$. Assume further that $Y(\,\cdot\, , T')=\tilde{Y}(\,\cdot\, , T')$ in $\Omega$, where $Y$ and $\tilde{Y}$ denote solutions of \eqref{sys.linearforCarleman} corresponding to the sources $\sigma(t)F(x)$ and $\sigma(t)\tilde{F}(x)$, respectively. Then, there exists a positive constant $C=C(\Omega,\mathcal{O}, T', T, \|q_{ij}\|_{\infty})$ such that the following estimate holds:
		\begin{equation}\label{ine.main.stability}
			\|F-\tilde{F}\|_{L^2(\Omega)^n}\leq C\|\partial_t y_n-\partial_t\tilde{y}_n\|_{L^2(0,T;L^2(\mathcal{O}))}.	
		\end{equation}
	\end{theorem}
	\begin{proof}
		Let $Y$ and $\tilde{Y}$ be solutions of \eqref{sys.linearforCarleman} associated to the sources $\sigma F$ and $\sigma\tilde{F}$, respectively. By defining $U:=Y-\tilde{Y}$ and $Z:=\partial_t U$, we obtain the system:
		\begin{equation}\label{sys.aux1.carleman}
   		\left\{
    	\begin{array}{llll}
        	\partial_t Z-\Delta Z+QZ=(\partial_t\sigma)(F-\tilde{F}), &\mbox{in} &\Omega\times(0,T),\\
         	Z=0, &\mbox{on} &\partial\Omega\times(0,T).
    	\end{array}\right.
		\end{equation}
		Using that $\sigma\in W^{1,\infty}(0,T)$ and $|\sigma(T')|\geq \gamma_{0}>0$, and by applying the Carleman estimate \eqref{ine.carleman} to system 
		\eqref{sys.aux1.carleman}, we obtain
		\begin{equation}\label{sys.aux2.carleman}
		\begin{array}{lll}
		\displaystyle\sum\limits_{i=1}^{n}\mathcal{I}(3(n+1-i),z_i)&\leq Cs^{\ell}\displaystyle\iint\limits_{\mathcal{O}\times(0,T)}
			e^{-2s\alpha}\xi(t)^{\ell}|z_n|^2\,dx\,dt\\
			&\hspace{1cm}+Cs^{r-3}\displaystyle\sum\limits_{i=1}^n\,
			\displaystyle\iint\limits_{\Omega\times(0,T)}e^{-2s\alpha}\xi(t)^{r-3}|f_i-\tilde{f}_i|^2\,dx\,dt
		\end{array}
		\end{equation}
		where $C$ is a positive constant independent of $F$ and $\tilde{F}$.
		
		On the other hand, since $e^{-2s\alpha(x,0)}=0$ for $x\in\overline{\Omega}$, we 
        have for $z_i\in H^1(0,T;L^2(\Omega))$ the following estimate
		\begin{equation}\label{sys.auxxx.carleman}\small{
		\begin{array}{lll}
			&s^{r-2}\displaystyle\int\limits_{\Omega}\xi^{r-3}(T')|z_i(x,T')|^2e^{-2s\alpha(x,T')}\,dx\\
			&=s\displaystyle\int\limits_{0}^{T'}\frac{\partial}{\partial t}\Biggl(\int\limits_{\Omega}(s\xi)^{r-3}|z_i|^2
				e^{-2s\alpha}\,dx\Biggr)\,dt\\
			&=s^{r-2}\displaystyle\int\limits_{\Omega}\int\limits_{0}^{T'}
				\Biggl( (\partial_t\xi^{r-3})|z_i|^2+2\xi^{r-3}z_i\partial_{t}z_i-2s\xi^{r-3}(\partial_t\alpha)|z_i|^2\Biggr)
				e^{-2s\alpha}\,dx\,dt\\
			&\leq C\displaystyle\iint\limits_{\Omega\times(0,T')}
            \Biggl((s\xi)^{r-2}|z_i|^2
			+(s\xi)^{r-4}|\partial_{t}z_i|^2+s^{r}\xi^{r-2}|z_i|^2+(s\xi)^{r-1}|z_i|^2\Biggr)
			e^{-2s\alpha}\,dx\,dt.		
		\end{array}}
		\end{equation}
 		Here, we used the fact that $|\partial_t\alpha(x,t)|\leq C\xi^2(t)$ and $|\partial_t\xi(t)|\leq C\xi^2(t)$
 		for $(x,t)\in \Omega\times(0,T)$, and Young's inequality 
 		(i.e., $ab\leq \frac{a^p}{p}+\frac{b^q}{q},\,\frac{1}{p}+\frac{1}{q}=1, a,b>0$) with 
 		$a=s^{\frac{r}{2}}\xi^{\frac{r-2}{2}}|z_i|,\, b=(s\xi)^{\frac{r-4}{2}}|\partial_{t}z_i|$ and $p=q=2$.
 		
 		Given $Y(\,\cdot\, ,T')=\tilde{Y}(\,\cdot\, , T')$, we have 
 		$Z(\,\cdot\, ,T')=\partial_{t}\sigma(T')(F(\cdot)-\tilde{F}(\cdot))$. Since 
        $|\sigma(T')|\geq \gamma_{0}>0$, using \eqref{sys.aux2.carleman} and \eqref{sys.auxxx.carleman}, we obtain 
 		\begin{equation*}
 		\begin{array}{lll}
 			&s^{r-2}\sum\limits_{i=1}^n\,\displaystyle\int\limits_{\Omega}
 			|f_i(x)-\tilde{f}_i(x)|^2e^{-2s\alpha(x,T')}\,dx\\
 			&\leq Cs^{r-2}\displaystyle\int\limits_{\Omega}\xi^{r-3}(T')
 			|\partial_{t}\sigma(T')(F(x)-\tilde{F}(x))|^2e^{-2s\alpha(x,T')}\,dx\\
 			&\leq C s^{r-2}\sum\limits_{i=1}^{n}\displaystyle\int\limits_{\Omega}\xi^{r-3}(T')|z_i(x,T')|^2
 			e^{-2s\alpha(x,T')}\,dx\\
 			&\leq Cs^{\ell}\displaystyle\iint\limits_{\mathcal{O}\times(0,T)}e^{-2s\alpha}\xi(t)^{\ell}|z_n|^2\,dx\,dt\\
			&\quad+Cs^{r-3}\sum\limits_{i=1}^n\,
			\displaystyle\iint\limits_{\Omega\times(0,T)}e^{-2s\alpha}\xi(t)^{r-3}|f_i(x)-\tilde{f}_i(x)|^2\,dx\,dt. 	
 		\end{array}
 		\end{equation*}
		Given that $\alpha(x,T')\leq \alpha(x,t)$ for $(x,t)\in \Omega\times(0,T)$ and $r\leq 3$, we can absorb the second term on the right-hand side by the left-hand side by taking $s>1$. Therefore, the proof of Theorem \ref{teo.stability} is complete.	
	\end{proof}
%\section{NEW FIGURES}
%\input{new_figures}